\documentclass[12pt,reqno]{amsart}

\usepackage{amsmath}
\usepackage{amssymb}
\usepackage{fullpage}
\usepackage{algorithm}
\usepackage{algpseudocode}
\usepackage{graphics}
\usepackage[colorlinks,citecolor=blue]{hyperref}
\usepackage{booktabs}
\usepackage{cite}
\usepackage{enumitem}
\usepackage{tikz,pgfplots}
\usetikzlibrary{shapes.multipart,backgrounds,shapes,fit}
\usetikzlibrary{arrows,decorations.markings,shapes.arrows,arrows.meta}
\usepackage{stix}
\usepackage{comment}

\newcommand{\R}{\mathbb{R}}
\renewcommand{\SS}{\mathbb{S}}
\newcommand{\TT}{\mathbb{T}}

\DeclareMathOperator{\aut}{\operatorname{Aut}}
\DeclareMathOperator{\auttr}{\operatorname{Aut_\triangle}}

\DeclareMathOperator{\tr}{\operatorname{Tr}}
\DeclareMathOperator{\Tr}{\operatorname{Tr}}

\DeclareMathOperator{\inte}{\operatorname{int}}

\DeclareMathOperator{\spam}{\operatorname{span}}

\newcommand{\Srank}{r}
\DeclareMathOperator{\SC}{\mathrm{SC}}

\newcommand{\cA}{\mathcal{A}}
\newcommand{\cB}{\mathcal{B}}
\newcommand{\cC}{\mathcal{C}}
\newcommand{\cH}{\mathcal{H}}
\newcommand{\cL}{\mathcal{L}}

\newcommand{\cP}{\mathcal{P}}

\newcommand{\cV}{\mathcal{V}}
\newcommand{\cW}{\mathcal{W}}

\newtheorem{thm}{Theorem}[section]

\newtheorem{dfn}[thm]{Definition}

\theoremstyle{remark}
\newtheorem{examp}[thm]{Example} 

\newtheorem{theorem}{Theorem}[section]
\newtheorem{proposition}[theorem]{Proposition}

\newcommand{\iprod}[2]{\left\langle {#1}, {#2} \right\rangle}

\newcommand{\adj}{\mathrm{adj}}
\newcommand{\madj}{\mathrm{adj}^+}
\newcommand{\ladj}{\mathrm{adj}^-}

\newcommand{\mL}{\mathcal{L}}

\newcommand{\ch}{\mathrm{ch}}
\newcounter{algorithmctr}[section]
\renewcommand{\thealgorithmctr}{\thesection.\arabic{algorithmctr}}
\newenvironment{algdesc}[1]%
    {\refstepcounter{algorithmctr}\begin{list}{}{%
        \setlength{\rightmargin}{.05\linewidth}%
        \setlength{\leftmargin}{.05\linewidth}}%
        \item[]{\setlength{\parskip}{0ex}\bigskip\par%
         \nopagebreak%
         \underline{{\bf Algorithm \thealgorithmctr.} \emph{#1.}}}}%
    {{\setlength{\parskip}{-1ex}\nopagebreak\smallskip\par} \end{list}}

\title[Linear Optimization over Homogeneous Matrix Cones]
{Linear Optimization over Homogeneous Matrix Cones}

\author{Levent Tun\c{c}el \and Lieven Vandenberghe}
\thanks{Levent Tun\c{c}el: Department of Combinatorics and Optimization, 
Faculty of Mathematics, University of Waterloo, Waterloo, Ontario N2L 3G1,
Canada (e-mail: ltuncel@uwaterloo.ca).
\\
Lieven Vandenberghe: Department of Electrical and Computer
Engineering, UCLA, Los Angeles, CA 90095, USA (e-mail: vandenbe@ucla.edu). 
}

\date{October 31, 2022}

\begin{document}

\label{firstpage}
\maketitle

\begin{abstract}
A convex cone is homogeneous if its automorphism group acts 
transitively on the interior of the cone.  
Cones that are homogeneous and self-dual are called symmetric. 
Conic optimization problems over symmetric cones have been extensively 
studied in convex optimization, in particular in the
literature on interior-point algorithms, and as the foundation of  
modeling tools for convex optimization.
In this paper, we consider the less well-studied conic optimization 
problems over cones that are homogeneous but not necessarily 
self-dual.  

We start with cones of positive semidefinite symmetric 
matrices with a given sparsity pattern.  Homogeneous cones in this
class are characterized by nested block-arrow sparsity patterns, a
subset of the chordal sparsity patterns.  
Chordal sparsity guarantees that positive define matrices in the cone 
have zero-fill Cholesky factorizations.  
The stronger properties that make the cone homogeneous
guarantee that the inverse Cholesky factors have the same zero-fill
pattern.  We describe transitive subsets of the cone automorphism groups, 
and important properties of the composition of log-det barriers with the 
automorphisms.

Next, we consider extensions to linear slices of the positive semidefinite 
cone, and review conditions that make such cones homogeneous.  
An important example is the matrix norm cone, the epigraph of a 
quadratic-over-linear matrix function.
The properties of homogeneous sparse matrix cones are shown to extend to 
this more general class of homogeneous matrix cones.
 
In the third part we give an overview of the
algebraic theory of homogeneous cones due to Vinberg and Rothaus.  
A fundamental consequence of this theory is that every
homogeneous cone admits a spectrahedral (linear matrix inequality) 
representation.   
 
We conclude by discussing the role of homogeneous
structure in primal--dual symmetric interior-point methods.
We make a contrast with the well-developed algorithms for symmetric cones
that exploit the strong properties of self-scaled barriers,
and with symmetric primal--dual methods for general convex cones.
\end{abstract}

\section{Introduction}
\label{sec:intro}

The conic programming framework has been used extensively in 
the development of convex optimization theory, applications, algorithms, 
and modeling \cite{NN1994,BeN:01,BoV:04}.  
As with any type of optimization problem, a fundamental 
step in a successful treatment of large-scale conic programs is the 
identification and efficient exploitation of special structure.
In this paper, we discuss convex cones represented as slices of the 
positive semidefinite cone, i.e., as intersections
\begin{equation} \label{e-K}
K = \mathcal V \cap \SS^N_+
\end{equation}
of $\SS^N_+$ (the cone of symmetric positive semidefinite $N$-by-$N$ 
matrices) and a subspace $\mathcal V$, and we examine the special 
structure of $\mathcal V$ that makes $K$ a \emph{homogeneous} convex cone.
A convex cone is homogeneous if for every pair of points in its
interior there exists an automorphism of the cone that maps one point 
to the other.

Inequalities with respect to slices of the positive semidefinite
cone arise in nonsymmetric formulations of semidefinite programming
problems.
Consider a semidefinite program (SDP) in inequality form 
\begin{subequations}\label{e-sdp}
\begin{equation} 
\begin{array}[t]{ll}
\mbox{minimize} & c^\top y \\
\mbox{subject to} & \sum\limits_{i=1}^m y_i A_i +X = B \\
            &  X \succeq 0
\end{array}
\end{equation}
and its dual problem,
\begin{equation} 
\begin{array}[t]{ll}
\mbox{maximize} & {-\iprod{B}{S}} \\
\mbox{subject to} & \iprod{A_i}{S} + c_i = 0, \;\; i=1,\ldots,m \\
            &  S \succeq 0.
\end{array}
\end{equation}
\end{subequations}
The primal variables are $y\in\R^m$, $X\in\SS^N$.  The dual variable 
is $S\in\SS^N$.  
The inequalities $X\succeq 0$, $S\succeq 0$ mean that $X, S \in \SS^N_+$.
The positive semidefinite matrix cone $\SS^N_+$ is a 
\emph{symmetric cone,} i.e., self-dual and homogeneous, and the special 
properties of symmetric cones are key to the design and 
implementation of primal--dual interior-point algorithms for semidefinite 
optimization.

If the matrices $A_1$, \ldots, $A_m$, $B$ all belong to a subspace 
$\mathcal V$ of $\SS^N$, the problems~(\ref{e-sdp}) are equivalent to 
the pair of conic optimization problems
\begin{subequations}\label{e-sdp-not-symmetric}
\begin{equation} 
\begin{array}[t]{ll}
\mbox{minimize} & c^\top y \\
\mbox{subject to} & \sum\limits_{i=1}^m y_i A_i +X = B \\
                  &  X \in K
\end{array}
\end{equation} 
and
\begin{equation} 
\begin{array}[t]{ll}
\mbox{maximize} & {-\iprod{B}{S}} \\
\mbox{subject to} & \iprod{A_i}{S} + c_i = 0, \;\; i=1,\ldots,m \\
   &  S \in K^*
\end{array}
\end{equation}
\end{subequations}
where $K$ is defined in~(\ref{e-K}), $K^*$ is the dual
cone of $K$, and the variables $X, S$ are matrices in~$\mathcal V$.   
The formulation~(\ref{e-sdp-not-symmetric}) is of interest for
large-scale algorithm development because the subspace $\mathcal V$ 
can be of much lower dimension than $\SS^N$, possibly as low as the
dimension of the span of the coefficient matrices
$\bar{\mathcal  V} = \spam{\{A_1, \ldots, A_m, B\}}$.
However, the efficiency of algorithms for handling the conic inequalities 
with respect to $K$ and $K^*$ depends on more properties of 
$\mathcal V$ than just the dimension, and this may require embedding 
$\bar{\mathcal V}$ in a higher-dimensional subspace.
The standard choice in current primal--dual interior-point methods is to 
embed $\bar{\mathcal  V}$ in a space of block-diagonal matrices with
dense diagonal blocks.
For this choice of $\mathcal V$, the cone $K$ is symmetric.
For almost all other subspaces $\mathcal V$, the cone $K$ and its
dual $K^*$ are not equal; hence they are not symmetric cones. 
(The exceptions are semidefinite representations of the small
number of symmetric cones, for example, direct products of
second order cones.)
However $K$ and $K^*$ may still be homogeneous.
Homogeneous convex cones were algebraically classified in the 1960s by 
Vinberg \cite{Vin:65} and are the subject of a large literature
in algebra and statistics~\cite{LetacMassam2007,AnW:04,BHM:11,KhR:11}.
The conditions for a matrix cone of the form~(\ref{e-K}) to be 
homogeneous have been studied by Letac and Massam 
\cite{LetacMassam2007} and Ishi \cite{Ishi2013,Ishi2015}. 
Homogeneous cones have several important properties in common with 
symmetric cones.  One can note, for example, that their definition 
contains two fundamental concepts in 
primal--dual interior-point algorithms for optimization over 
symmetric cones.  The automorphisms of a cone
(invertible linear transformations that leave the cone invariant) are the 
\emph{scalings} used in interior-point methods, for example,
the positive diagonal scalings of the nonnegative orthant in 
algorithms for linear programming. 
The second property, that the automorphisms act transitively in the 
interior of the cone, implies that any given pair of primal and dual 
iterates can mapped to the same point by a cone automorphism, as we 
will discuss in Section~\ref{s-barriers}.
Hence, homogeneous cones are a natural subject of study in
conic optimization.  However, with some notable exceptions
\cite{Guler1996,GulerTuncel1998,TruongTuncel2004,Chua2009}, work on
algorithms for homogeneous conic optimization appears to be quite limited. 
It is the purpose of this article to describe 
properties of homogeneous matrix cones that are useful in 
algorithms for optimization problems of the 
form~(\ref{e-sdp-not-symmetric}).
We also discuss specific examples and structural properties that may be 
useful for optimization modeling tools.

In Sections~\ref{s-graphs}--\ref{s-barriers} we first consider 
matrix subspaces $\mathcal V$ defined by sparsity patterns.   
If the coefficient matrices $A_1$, \ldots, $A_m$, $B$ 
in  problem~(\ref{e-sdp}) have a common (aggregate) sparsity pattern 
then the subspace $\mathcal V$ 
in~(\ref{e-sdp-not-symmetric}) can be defined as the set of symmetric
$N$ by $N$ matrices with that pattern, or any extension of the
aggregate sparsity pattern.
The primal cone $K$ is the cone of positive semidefinite
matrices with a given sparsity pattern;
the dual cone $K^*$ is the cone of symmetric matrices
with the same sparsity pattern that have a positive semidefinite 
completion.
The non-symmetric conic formulation~(\ref{e-sdp-not-symmetric}) has 
been studied in recent approaches to exploit sparsity in sparse 
semidefinite 
optimization \cite{FKMN2001,BYZ:00,ADV:10,SRV:04,Bur:03}.
Table~\ref{t-cones} summarizes the definitions that relate this paper
to existing literature on semidefinite programming.
\begin{table}
\begin{tabular}{c@{\hskip 2em}c@{\hskip 2em}c@{\hskip 2em}c}  \toprule
Sparsity pattern & Linear algebra & Convex cone \\\midrule
dense & spectral theory & symmetric \\
& & & \\
homogeneous chordal 
& zero-fill Cholesky factor 
& homogeneous  \\
& and inverse factor &   \\
& & & \\
chordal & zero-fill Cholesky factor & slice of PSD cone  \\ 
& & & \\
general & sparse Cholesky factor & slice of PSD cone  \\
\bottomrule \\*[.5ex]
\end{tabular}
\caption{Four classes of sparse positive semidefinite matrix cones, 
classified by type of sparsity, the linear algebra tools available for 
their analysis, and fundamental properties of the cones.}
\label{t-cones}
\end{table}
It distinguishes sparse positive semidefinite matrix cones by type of 
sparsity.
At the top level, we have the dense positive semidefinite cones
(i.e., without any restriction on the sparsity pattern).
The dense positive semidefinite cone is symmetric (self-dual and
homogeneous). Symmetric primal--dual algorithms for them rely heavily on 
eigenvalue and generalized  eigenvalue decompositions of symmetric
positive semidefinite matrices (for example, for computing
the matrix geometric mean, or for joint diagonalization of positive 
definite matrices).
At the lowest level of the table we have the positive semidefinite matrix
cones with a general, unstructured sparsity pattern.
They form lower-dimensional slices of the positive semidefinite cone.  
Such cones are convex, but not homogeneous or self-dual. 
Implementations of non-symmetric interior-point algorithms for these
cones, for example, dual barrier algorithms \cite{BYZ:00},  
benefit from the possibility of computing sparse Cholesky factors, 
using fill-reducing ordering heuristics. 
Level three in the table is occupied by 
the positive semidefinite matrices with chordal sparsity patterns.
Chordal sparsity has been studied intensively in sparse semidefinite 
optimization (see \cite{VaA:15,ZFP:21} for recent surveys).
The chordal structure can be exploited to formulate efficient algorithms 
for key computations needed in semidefinite programming 
algorithms, such as the evaluation of primal and dual barrier functions 
and their derivatives,
and finding maximum-determinant or minimum-rank positive semidefinite 
completions  \cite{GrT:84,AHMR:88,GJSW:84}.
All these algorithms can be derived from the basic property that 
positive semidefinite matrices with a chordal sparsity patterns have a 
zero-fill Cholesky factorization.
The second row of the table is the focus of 
Sections~\ref{s-graphs}--\ref{s-barriers} of this paper.
The sparsity patterns that are referred to here as ``homogeneous chordal''
define matrix cones that are homogeneous but not necessarily symmetric.
These sparsity patterns have been characterized by 
Letac and Massam \cite[Theorem 2.2]{LetacMassam2007} and
Ishi \cite[Theorem A]{Ishi2013}. 
As we will discuss in Sections~\ref{s-graphs} and~\ref{s-cone}, 
they are block-arrow sparsity patterns and recursive 
generalizations of block-arrow structures.
They form a subset of the chordal patterns, with the additional
useful property that the inverse Cholesky factor has the same, 
zero-fill, sparsity pattern as the Cholesky factor itself.

Note that any class of semidefinite programming problems on a higher level 
in the table includes the lower ones.
One can always extend, at no loss of generality, a general sparsity
pattern to make it chordal, or a chordal pattern to make it 
homogeneous chordal, 
or a homogeneous chordal sparsity pattern to make it dense.
However, there is an obvious trade-off.
The higher levels come with stronger results and more powerful
techniques from linear algebra, and with more efficient primal, dual,
or primal--dual conic optimization algorithms.  
They also embed the optimization problem in 
higher-dimensional spaces and exploit less of the detailed structure
in the sparsity pattern.

The three sections on homogeneous sparse matrix cones are organized
as follows.  Section~\ref{s-graphs} is a survey of 
results and algorithms from sparse matrix and graph theory related 
to chordal and homogeneous chordal sparsity patterns.
In Section~\ref{s-cone} we show that the positive semidefinite
cone with a homogeneous chordal pattern and the associated dual cone
are homogeneous.  We establish a transitive subset of the automorphism
group constructed from congruences with sparse lower-triangular matrices.
In Section~\ref{s-barriers} we derive implications for the log-det barrier 
function and its conjugate.  We show that the Hessians of the logarithmic
barrier functions can be factorized as a composition of a cone 
automorphism and its adjoint.  
This leads to a generalization of the Nesterov--Todd scaling point for 
symmetric cones.

In Section~\ref{s-hom-matrix-cones} we then turn to more general
homogeneous slices of the positive semidefite matrix cone,
with subspaces $\mathcal V$ that can be defined by other linear relations
than the sparsity pattern.
The properties of $\mathcal V$ that make the cone~(\ref{e-K})
homogeneous are described by Ishi \cite{Ishi2015}.
The results in this section will parallel the properties of homogeneous
sparse matrix cones.  In particular, Cholesky factors and inverse 
Cholesky factors inherit the structure of the subspace $\mathcal V$.

Section~\ref{sec:homogeneous-cones} reviews the general, algebraic 
classifications of all homogeneous cones and connects these theories to 
the earlier sections.
An important result is that every homogeneous cone has 
a semidefinite representation, i.e., is linearly isomorphic to 
a slice of the positive semidefnite cone.

We conclude the paper with a survey of recent work on
interior-point methods for nonsymmetric conic optimization, and
point out the potential benefits of exploiting the special
properties of homogeneous cones (Section~\ref{s-ipm}).
The two appendices contain background material from graph theory
and algorithmic details.

The paper is primarily intended as a survey.  Its main contributions 
are the following.
\begin{itemize}
\item
We identify a class of conic optimization problems (based on homogeneous 
sparse matrix cones, called \emph{homogeneous chordal cones})
which lie strictly between SDPs and homogeneous cone programming problems 
(in the context of the set of convex cones $K$ allowed in the optimization 
problems~\eqref{e-sdp-not-symmetric}).
In this context, the class of convex optimization problems over 
homogeneous chordal cones provides a generalization of second order cone 
programming that has important computational advantages over 
semidefinite programming.

\item
We build on results from convex optimization and analysis, graph theory, data structures and algorithms, sparse matrix computation and theory, 
abstract algebra and show how to perform fundamental linear algebra 
operations in an efficient way for many
families of algorithms for our class of conic optimization problems.

\item
We show how to compute primal and dual scalings that are automorphisms 
of the underlying cones and in doing so we solve an open problem about 
the existence of automorphism based primal--dual scalings for pairs of 
interior-points in homogeneous cones and in their duals.

\item We extend the results from to homogeneous sparse matrix cones 
 to homogeneous matrix cones defined by slices of the positive
 semidefinite cone.  Constraints of this type are important in
 semidefinite representations of the spectral matrix norm and the trace
 norm.
\end{itemize}

\section{Homogeneous chordal sparsity} \label{s-graphs}
We denote by $\SS^N$ the space of $N$-by-$N$ symmetric matrices with
real entries, by $\SS^N_+$ the convex cone of positive semidefinite
matrices in $\SS^N$, and by $\SS^N_{++} := \inte(\SS^N_+)$ the cone
of positive definite matrices in $\SS^N$.   
For $X,Y\in \SS^N$, the inequalities $X\succeq Y$ and $X\succ Y$ mean 
that $X-Y\in \SS^N_+$ and $X-Y\in \SS^N_{++}$, respectively.
The standard trace inner product is used for~$\SS^N$:
\[
\langle X, Y \rangle = \Tr(XY) = \sum_{i=1}^N\sum_{j=1}^N X_{ij}Y_{ij}.
\]
The set of $N$-by-$N$ lower-triangular matrices with real entries is 
denoted by $\TT^N$.

\subsection{Sparse matrices}
An $N$-by-$N$ symmetric sparsity pattern is represented by a 
simple undirected graph $G = (V,E)$ with vertex set 
$V = \{1,2,\ldots,N\}$ and edge set $E$.  
An edge connecting vertices $i$ and $j$ is denoted by $\{i,j\}$.
A matrix $X\in\SS^N$ is said to have the sparsity pattern $E$
if $X_{ij} =  X_{ji} = 0$ whenever $i\neq j$ and $\{i,j\} \not\in E$.
The diagonal entries and the entries indexed by $E$ are called the
\emph{nonzeros} in the pattern.  The other entries (indexed by the 
complement of $E$) are the \emph{zeros}.
The set of symmetric $N$-by-$N$  matrices with sparsity pattern $E$
is denoted by $\SS^N_E$:
\[
\SS^N_E := \{ X\in \SS^N : \, X_{ij} = X_{ji} = 0 \mbox{\ if $i\neq j$
and $\{i,j\} \not\in E$}\}.
\]
We use $\Pi_E$ to denote orthogonal projection on $\SS^N_E$.
For $X\in \SS^N$, the matrix $\Pi_E(X)$ is the matrix in $\SS^N_E$ with
nonzero entries given by $(\Pi_E(X))_{ij} = X_{ij}$ if $i=j$ or 
if $i\neq j$ and $\{i,j\}\in E$. 

The cone of positive semidefinite matrices in $\SS^N_E$ is the intersection
\begin{equation} \label{e-psd-cone}
 \SS^N_E \cap \SS^N_+  = \{ X\in \SS^N_E : \, X\succeq 0\}.
\end{equation}
This cone is closed, convex, and pointed.  It also has nonempty interior
relative to $\SS^N_E$ (it includes the identity matrix $I$),
so it is a \emph{regular} (or \emph{proper}) cone.
The cone of matrices in $\SS^N_E$  that have a positive 
semidefinite completion is the projection of $\SS^N_+$ on $\SS^N_E$.
We denote this set by 
\begin{equation} \label{e-psdc-cone}
 \Pi_E(\SS^N_+) = \{ \Pi_E(Y) : \, Y \succeq 0\}.
\end{equation}
The cone $\Pi_E(\SS^N_+)$ is clearly convex, pointed, and has nonempty 
interior relative to~$\SS^N_E$.
Closedness follows from the fact that if $\Pi_E(Y) = 0$ and $Y\succeq 0$
then $Y=0$.  Hence, the positive semidefinite completable 
cone $\Pi_E(\SS^N_+)$ is also regular.
The two cones $\SS^N_E \cap \SS^N_+$ and $\Pi_E(\SS^N_+)$ are duals
of each other under the trace inner product in the space $\SS^N_E$.

The graph $(V,E)$ can also be used to describe the sparsity pattern 
of lower-triangular matrices.  We say $L\in \TT^N$ has sparsity pattern 
$E$ if $L+L^\top \in \SS^N_E$.  The notation
\[
\TT^N_E = \{ L \in \TT^N : \, L+L^\top \in \SS^N_E\}
\]
will be used for this set.  

We define the Cholesky factorization of a positive 
definite matrix $X$ as a decomposition
\begin{equation} \label{e-chol}
 PXP^\top = LL^\top
\end{equation}
where $P$ is a permutation matrix and  $L$ is lower-triangular with
positive diagonal entries.  
In general, the factorization introduces fill in the sparsity pattern of 
$PXP^\top$.  We say the sparsity pattern of $L$ is an \emph{extension} of 
the sparsity pattern of $PXP^\top$. 

\subsection{Chordal sparsity}
We now give a short overview of the properties of chordal graphs 
and chordal sparsity patterns that will be important in the discussion
of homogeneous chordal patterns in the next section.
The interested reader is referred to the 
surveys~\cite{VaA:15,BlP:93,Gol:04,ZFP:21} for more background on 
chordal graphs and their history.

An undirected graph $(V,E)$ is called \emph{chordal} if it does not contain
a cycle $C_k$ of length $k\geq 4$ as a node induced subgraph (from now on,
we will simply say \emph{induced graph} to mean \emph{node induced graph}).
A classical result states that a graph is chordal if and only if it has 
a \emph{perfect elimination ordering} \cite{FuG:65}. 
An ordering of the graph is a bijection $\sigma$ from $\{1,2,\ldots,|V|\}$
to the vertex set $V$.
An ordering $\sigma$ is a perfect elimination ordering if
\begin{equation} \label{e-peo2}
\left. \begin{array}{l}
\{u,v\} \in E, \; \{u,w\} \in E, \\
\sigma^{-1}(u) < \sigma^{-1}(v) < \sigma^{-1}(w) 
\end{array}\right\}
\qquad \Longrightarrow \qquad \{v,w\} \in E.
\end{equation}
In other words, the \emph{higher neigborhood} 
\[
\madj(u) := \{ v \in V : \, \{u,v\} \in E, \, 
 \sigma^{-1}(u) < \sigma^{-1}(v)\}
\]
of every vertex induces a complete subgraph of $G$:
\begin{equation} \label{e-peo}
v, w \in \madj(u) \quad \Longrightarrow \quad \{v,w\} \in E.  
\end{equation}
In sparse matrix language, a perfect elimination ordering of a sparsity 
pattern $E$ defines a  permutation matrix that yields a zero-fill
Cholesky factorization~(\ref{e-chol}), i.e., 
$P^\top(L+L^\top)P\in\SS^N_E$ if $X \in \SS^N_E$.

Efficient linear-time algorithms exist for testing chordality
of a graph and finding a perfect elimination ordering if one exists
\cite{RTL:76,TaY:84}.
For non-chordal graphs, the connection with the sparse Cholesky 
factorization~(\ref{e-chol}) suggests a practical heuristic for finding
efficient chordal extensions: apply a fill-reducing reordering 
to the sparsity pattern of $X$ and calculate the sparsity pattern of
the Cholesky factor $L$.  

\begin{figure}
\hspace*{\fill}
\begin{minipage}{.45\linewidth}
\hspace*{\fill}
\begin{tikzpicture}[scale=0.4]
    \tikzset{Bullet/.style = { shape = circle, minimum size = 3pt, 
        inner sep = 0pt, fill=black, draw=black, thick}}
    \foreach \i/\j in {
        1/0, 6/0, 7/0,
        4/1, 6/1, 7/1, 
        4/2, 7/2,
        5/3, 8/3,
        6/4, 7/4, 
        7/5, 8/5,
        7/6, 8/6,
        8/7 
    }{
        \node[Bullet] at (\j,-\i){};
        \node[Bullet] at (\i,-\j){};
      }
   \foreach \k in { 1, 2, 3, 4, 5, 6, 7, 8, 9} 
       \node[font=\footnotesize] at (\k-1,-\k+1){$\k$};
   \draw[thick] (-0.5,0.5) rectangle (8.5,-8.5);
\end{tikzpicture}
\hspace*{\fill}
\end{minipage}
\begin{minipage}{.45\linewidth}
\hspace*{\fill}
\begin{tikzpicture}[xscale=8.00, yscale = 8] 
   \tikzset{VertexStyle/.style = {shape = circle, draw, 
       minimum size = 12pt, inner sep = 1, fill=none}}
   \foreach \k/\x/\y in {
       3/0.2/0.1,
       4/0.3/0.2,
       5/0.1/0.2,
       6/0.3/0.3,
       9/0.2/0.5,
       1/0.0/0.0,
       2/0.0/0.1,
       7/0.1/0.3,
       8/0.2/0.4}
       \node[VertexStyle, font=\small](\k) at (\x,\y){$\k$};
   \foreach \i/\j in {1/2, 2/5, 3/5, 4/6, 5/7, 6/8, 7/8, 8/9} 
       \draw[thick] (\i)--(\j);
\end{tikzpicture}
\hspace*{\fill}
\end{minipage}
\hspace*{\fill}
\caption{\emph{Left.} A chordal graph with vertices 
$V=\{1,2,\ldots, 9\}$ and perfect elimination ordering $1, \ldots, 9$.  
The dots in the array represent the edges in the graph.  
\emph{Right.} Elimination tree.}
\label{f-tree}
\end{figure}
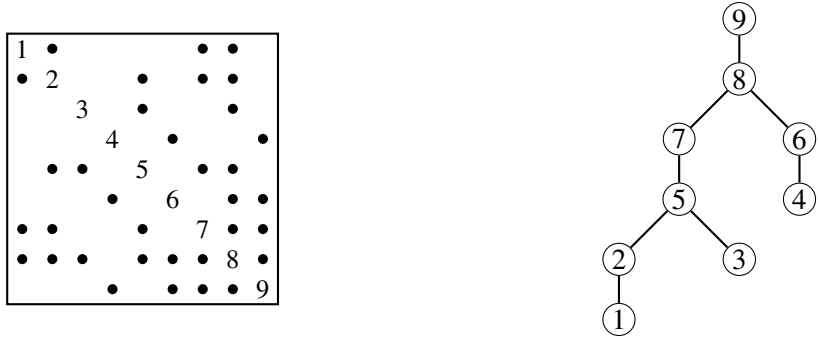

Elimination trees play an important role in sparse matrix algorithms,
such as the multifrontal algorithm for sparse Cholesky 
factorization~\cite{DuR:83,Liu:90}.
The elimination tree of a chordal graph $G$ with perfect elimination
ordering $\sigma$ is a tree (or a forest if the graph is not connected),
with vertex set $V$.
The parent $p(u)$ of a non-root vertex $u$ in the tree is the first 
element of $\madj(u)$.  The perfect elimination 
property~(\ref{e-peo}) holds if and only if 
\begin{equation} \label{e-tree-peo}
 \madj(u) \subseteq \{p(u)\} \cup \madj(p(u))
\end{equation}
for all non-root vertices $u$.
Figure~\ref{f-tree} shows an example.

It is useful to note that the elimination tree provides a summary 
of the graph, but is not an equivalent representation.  
For example, from the elimination tree 
in Figure~\ref{f-tree} and the property~(\ref{e-tree-peo}), we can 
conclude that 
vertex $6$ is not adjacent to vertex $1$; however, the information in the 
elimination tree does not allow us to decide whether vertex $5$ is 
adjacent to vertex $1$ or not.

\subsection{Homogeneous chordal sparsity}
\begin{figure}
\hspace*{\fill}
\hspace*{\fill}
\begin{tikzpicture}[scale=1.0,baseline]
   \centering
   \tikzset{Vertex/.style = {shape = circle, draw, 
       minimum size = 10pt, inner sep = 1, fill=none}}
   \foreach \k / \x / \y in {
       1/0.0/ 1.0, 2/1.0/ 1.0, 3/1.0/ 0.0, 4/0.0/ 0.0}
       \node[Vertex, font=\small](\k) at (\x,\y){};
   \foreach \i/\j in {1/2, 2/3, 3/4, 4/1} \draw (\i)--(\j);
\end{tikzpicture}
\hspace*{\fill}
\begin{tikzpicture}[scale=1.0,baseline]
   \centering
   \tikzset{Vertex/.style = {shape = circle, draw, 
       minimum size = 10pt, inner sep = 1, fill=none}}
   \foreach \k / \x / \y in {
       1/0.0/ 0.5,
       2/1.0/ 0.5,
       3/2.0/ 0.5,
       4/3.0/ 0.5}
       \node[Vertex, font=\small](\k) at (\x,\y){};
   \foreach \i/\j in {1/2, 2/3, 3/4} \draw (\i)--(\j);
\end{tikzpicture}
\hspace*{\fill}
\hspace*{\fill}
\caption{$C_4$ and $P_4$ are forbidden induced subgraphs in 
a homogeneous chordal graph.}
\label{f-subgraphs}
\end{figure}
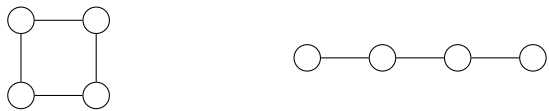

We define a \emph{homogeneous chordal graph} as an undirected graph 
that does not contain $C_4$ (a cycle of length four)  or $P_4$
(a path formed by three edges on four vertices) as induced subgraphs.
These forbidden subgraphs are shown in Figure~\ref{f-subgraphs}.
It is clear from the definition that a homogeneous chordal graph 
does not contain any induced cycle $C_k$ of length $k \geq 5$; 
so, all homogeneous chordal graphs are chordal.

Homogeneous chordal graphs were first studied by Wolk 
\cite{Wol:62,Wol:65}, who called them \emph{D-graphs}.  
Golumbic proposed the more commonly used term
\emph{trivially perfect} graphs \cite{Gol:78}.
They are known as \emph{homogeneous graphs} 
in the statistics literature on Gaussian graphical models 
\cite{LetacMassam2007,KhR:12}.
Other names include
\emph{quasi-threshold graphs} \cite{YCC:96},
\emph{co-chordal graphs} \cite{KhR:12}, and 
\emph{chordal co-graphs}\footnote{Graphs that do not contain $P_4$ are 
also known as \emph{co-graphs}
(complement reducible graphs), $D^*$-graphs, or hereditary Dacey graphs 
(due to a connection to work on orthomodular lattices).
So the homogeneous chordal graphs are the chordal co-graphs.}.
Our motivation for the name \emph{homogeneous chordal graphs}
will become clear in Section~\ref{s-cone}.

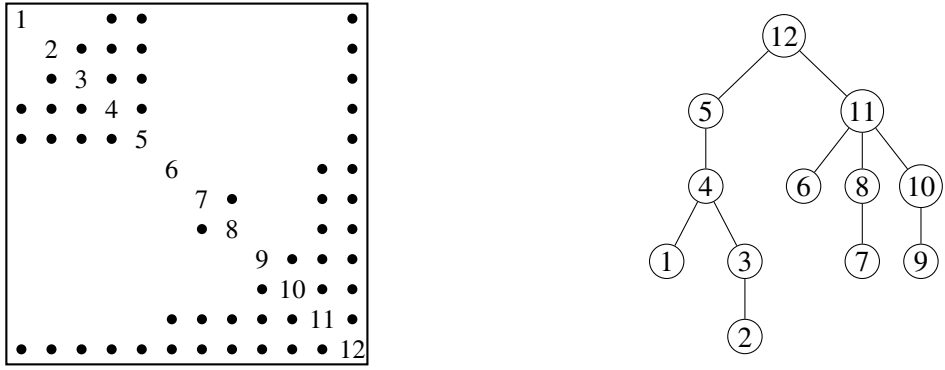
\begin{figure}
\hspace*{\fill}
\begin{minipage}{.45\linewidth}
\centering
\begin{tikzpicture}[scale=0.4]
   \tikzset{Bullet/.style = {
       shape = circle, minimum size = 3pt, inner sep = 0pt, fill=black, 
       draw=black, thick}}
   \foreach \i/\j in { 
        1/  5,
        1/  4,
        1/ 12,
        2/  3,
        2/  5,
        2/  4,
        2/ 12,
        3/  5,
        3/  4,
        3/ 12,
        5/  4,
        5/ 12,
        4/ 12,
        6/ 11,
        6/ 12,
        7/  8,
        7/ 11,
        7/ 12,
        8/ 11,
        8/ 12,
       10/  9,
       10/ 11,
       10/ 12,
        9/ 11,
        9/ 12,
       11/ 12
   } { 
       \node[Bullet] at (\j-1,-\i+1){}; 
       \node[Bullet] at (\i-1,-\j+1){}; 
   }
   \foreach \i in { 1, 2, 3, 5, 4, 6, 7, 8, 9, 10, 11, 12} 
       \node[font=\footnotesize] at (\i-1,-\i+1) {$\i$};  
   \draw[thick] (-0.5,0.5) rectangle (11.5,-11.5);
\end{tikzpicture}
\end{minipage}
\hspace*{\fill}
\begin{minipage}{.45\linewidth}
\centering
\begin{tikzpicture}[xscale=1.3, yscale=1.0]
   \tikzset{VertexStyle/.style = {shape = circle, draw, 
       minimum size = 13pt, inner sep = 1, fill=none}}
   \foreach \v/\k/\x/\y in {
       2/ 1/-1.2/1,
       1/ 2/-0.4/0,
       9/ 3/-0.4/1,
       6/ 4/-0.8/2,
      12/ 5/-0.8/3,
       8/12/ 0  /4,
       5/11/ 0.8/3,
      11/ 6/ 0.2/2,
       7/ 8/ 0.8/2,
       4/ 7/ 0.8/1,
      10/10/ 1.4/2, 
       3/ 9/ 1.4/1 }
      \node[VertexStyle, font=\normalsize](\v) at (\x,\y){\small $\k$};
   \foreach \i/\j in {2/6, 1/9, 9/6, 6/12, 12/8, 5/8, 11/5, 7/5, 4/7,
      10/5, 3/10} \draw (\i)--(\j);
\end{tikzpicture}
\end{minipage}
\hspace*{\fill}
\caption{The homogeneous chordal graph on the left is the
comparability graph of the tree on the right.  This tree is also the
elimination tree for the perfect elimination ordering $1, \ldots, 12$.} 
\label{f-arrow}
\end{figure}

Wolk \cite{Wol:62,Wol:65} showed that the absence of $P_4$ and $C_4$
characterizes the comparability graphs of rooted forests:  
a graph $G=(V,E)$ is a homogeneous chordal graph 
if and only if there exists a rooted forest with vertex
set $V$ and such that $\{v,w\} \in E$ if and only if $v$ is an ancestor 
of $w$ or $w$ is an ancestor of $v$ in the forest (in which case
we call $v$ and $w$ \emph{comparable} vertices). 
As a key step in his proof, he also established the important property 
that every connected component of a homogeneous chordal graph has a 
\emph{universal vertex}, i.e., a vertex adjacent to all other vertices
in the same connected component \cite[page 18]{Wol:62}.
This leads to a useful recursive characterization \cite{YCC:96}.  
Every homogeneous chordal graph can be constructed starting from 
a single-vertex graph by a repeated application of the following 
two operations.
\begin{itemize} \label{p-recursive-D-graphs}
\item \emph{Disjoint union.}  If $(V_1,E_1)$ and $(V_2,E_2)$
are homogeneous chordal graphs and $V_1 \cap V_2 = \emptyset$,
then $(V_1\cup V_2, E_1\cup E_2)$ 
is a homogeneous chordal graph.
\item \emph{Addition of a universal vertex.} If $(V,E)$ is a homogeneous
chordal graph and $w\not\in V$, then 
$\left(V \cup \{w\}, E \cup \{\{w,v\} : \, v \in V\}\right)$
is a homogeneous chordal graph.
\end{itemize}
These two operations have a simple interpretation for graphs that
describe sparsity patterns.  By making a disjoint union we construct
a sparsity pattern of size $N_1+N_2$ as a block-diagonal pattern
with diagonal blocks of size $N_1$  and $N_2$ (up to a symmetric 
reordering). 
Adding a universal vertex to a sparsity pattern of size $N_1$-by-$ N_1$
corresponds to adding a dense row and column to define a pattern
of size $(N_1+1) $-by-$ (N_1+1)$.
By repeating the two operations we construct a nested block-arrow
pattern (up to a symmetric reordering).  
Figure~\ref{f-arrow} shows an example.

Chu \cite{Chu:08} presents a linear-time algorithm for recognizing
homogeneous chordal graphs. 
The algorithm, described in detail in Appendix~\ref{s-app-graphs},
is an instance of the
\emph{Lexicographic Breadth First Search} (LBFS) algorithm that was first 
developed for testing chordality \cite{RTL:76} and later extended
for testing a variety of other graph properties \cite{Corneil2004,HMPV:00}.
Chu's algorithm also produces a perfect elimination ordering 
and an elimination tree.
The perfect elimination ordering $\sigma$ produced by the LBFS algorithm
has the following property, in addition to~(\ref{e-peo2}):
\begin{equation} \label{e-peo3}
\left.
\begin{array}{l}
\{u,v\} \in E, \; \{v,w\} \in E  \\
\sigma^{-1}(u) < \sigma^{-1}(v) < \sigma^{-1}(w) 
\end{array}\right\}
\qquad \Longrightarrow \qquad \{u,w\} \in E.
\end{equation}
Combined with~(\ref{e-peo2}) this
implies that two vertices are adjacent in the graph 
if and only if they form an ancestor--descendant pair
in the elimination tree: the homogeneous chordal 
graph is the comparability graph of the elimination tree.  
We will call a perfect elimination ordering that 
satisfies~(\ref{e-peo3}) a \emph{trivially perfect elimination ordering}.
For a trivially perfect elimination ordering, 
property~(\ref{e-tree-peo}) can be strengthened to
\begin{equation} \label{e-tree-peo-tp}
 \madj(u) = \{p(u)\} \cup \madj(p(u)).
\end{equation}
Hence, in contrast to general chordal patterns, a homogeneous chordal 
graph is completely characterized by an elimination tree.
This is illustrated in Figure~\ref{f-arrow}.
Here the numerical ordering is a trivially perfect elimination ordering 
of the homogeneous chordal graph on the left.
Each vertex in this graph is adjacent to all its ancestors and
descendants in the elimination tree.
The ordering in this example is also a \emph{postordering,} i.e., 
if $\sigma^{-1}(v)=j$ and $v$ has $k$ descendants in the elimination tree, 
then the descendants are numbered $j-1$, \ldots, $j-k$.
The postordering property holds for all trivially perfect 
elimination orderings computed by LBFS (see Appendix~\ref{s-app-graphs}).

\begin{figure}
\hspace*{\fill} 
\begin{minipage}{.2\linewidth}
\centering

\hspace*{\fill}
\begin{tikzpicture}[scale=0.5]
   \tikzset{Bullet/.style = {
       shape = circle, minimum size = 3pt, inner sep = 0pt, fill=black, 
       draw=black, thick}}
   \foreach \i/\j in { 
        1/3,
        2/3
   } { 
       \node[Bullet] at (\j-1,-\i+1){}; 
       \node[Bullet] at (\i-1,-\j+1){}; 
   }
   \foreach \i\j in { 1/1, 2/2, 3/3 }
       \node[font=\footnotesize] at (\i-1,-\i+1) {$\j$};  
   \draw[thick] (-0.4,0.4) rectangle (2.4,-2.4);
\end{tikzpicture}
\end{minipage}
\hspace*{\fill}
\begin{minipage}{.2\linewidth}
\centering
\begin{tikzpicture}[xscale=1.1, yscale=1.1, label distance = -.2em]
   \tikzset{VertexStyle/.style = {shape = circle, draw, 
       minimum size = 12pt, inner sep = 1, fill=none}}
   \foreach \v/\k/\x/\y in {
       1/1/0.0/0.0,
       2/2/1.0/0.0,
       3/3/0.5/1.0 }
      \node[VertexStyle, font=\normalsize](\v) 
          [label = right: \scriptsize $\k$]
          at (\x,\y){\small $\v$};
   \foreach \i/\j in {1/3, 2/3} \draw (\i)--(\j);
\end{tikzpicture}
\end{minipage}
\hspace*{\fill} \hspace*{\fill} 
\hspace*{\fill} \hspace*{\fill}
\begin{minipage}{.2\linewidth}
\centering
\begin{tikzpicture}[scale=0.5]
   \tikzset{Bullet/.style = {
       shape = circle, minimum size = 3pt, inner sep = 0pt, fill=black, 
       draw=black, thick}}
   \foreach \i/\j in { 
        1/2,
        2/3
   } { 
       \node[Bullet] at (\j-1,-\i+1){}; 
       \node[Bullet] at (\i-1,-\j+1){}; 
   }
   \foreach \i\j in { 1/1, 2/3, 3/2 }
       \node[font=\footnotesize] at (\i-1,-\i+1) {$\j$};  
   \draw[thick] (-0.4,0.4) rectangle (2.4,-2.4);
\end{tikzpicture}
\end{minipage}
\begin{minipage}{.2\linewidth}
\centering
\begin{tikzpicture}[xscale=1.1, yscale=1.0, label distance = -.2em]
   \tikzset{VertexStyle/.style = {shape = circle, draw, 
       minimum size = 12pt, inner sep = 1, fill=none}}
   \foreach \v/\k/\x/\y in {
       1/ 1/0.0/0.0,
       2/ 3/0.0/1.0,
       3/ 2/0.0/2.0 }
      \node[VertexStyle, font=\normalsize](\v) 
          [label = right: \scriptsize $\v$]
          at (\x,\y){\small $\k$};
   \foreach \i/\j in {1/2, 2/3} \draw (\i)--(\j);
\end{tikzpicture}
\end{minipage}
\hspace*{\fill}\hspace*{\fill}
\caption{Two perfect elimination orderings of a homogeneous
chordal graph and the corresponding elimination trees.  
The number next to node $v$ in the elimination trees is $\sigma^{-1}(v)$,
the position of $v$ in the ordering.
The ordering on the left is a trivially perfect elimination ordering.
The ordering on the right is a perfect elimination ordering, but is
not trivially perfect.} 
\label{f-vinberg}
\end{figure}
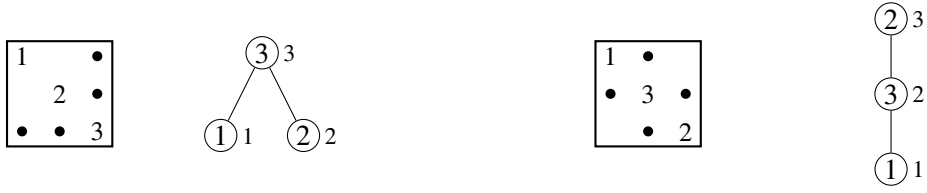

Note that not every perfect elimination ordering of a 
homogeneous chordal graph satisfies~(\ref{e-tree-peo-tp}).
Figure~\ref{f-vinberg} shows the smallest non-trivial (not dense
and not diagonal) sparsity pattern.
The figure shows two perfect elimination orderings and the corresponding 
elimination trees.
The first ordering is trivially perfect. The second ordering is not
because 
\[
\madj(1) = \{3\} \neq \{p(1)\} \cup \madj(p(1)) = \{2, 3\}.
\]

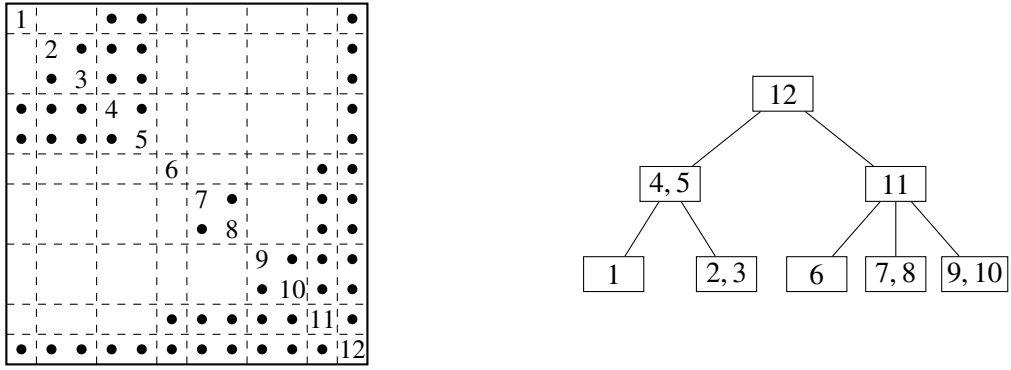
\begin{figure}
\hspace*{\fill}
\begin{minipage}{.45\linewidth}
\centering
\begin{tikzpicture}[scale=0.4]
   \tikzset{Bullet/.style = {
       shape = circle, minimum size = 3pt, inner sep = 0pt, fill=black, 
       draw=black, thick}}
   \foreach \i/\j/\k/\l in { 
        1/  5/ 2/12, 
        1/  4/ 2/ 6, 
        1/ 12/ 2/ 8, 
        2/  3/ 1/ 9, 
        2/  5/ 1/12,
        2/  4/ 1/ 6, 
        2/ 12/ 1/ 8, 
        3/  5/ 9/12,
        3/  4/ 9/ 6, 
        3/ 12/ 9/ 8, 
        5/  4/12/ 6, 
        5/ 12/12/ 8, 
        4/ 12/ 6/ 8, 
        6/ 11/11/ 5,  
        6/ 12/11/ 8, 
        7/  8/ 4/ 7, 
        7/ 11/ 4/ 5, 
        7/ 12/ 4/ 8, 
        8/ 11/ 7/ 5,
        8/ 12/ 7/ 8,
       10/  9/10/ 3,
       10/ 11/10/ 5, 
       10/ 12/10/ 8, 
        9/ 11/ 3/ 5, 
        9/ 12/ 3/ 8, 
       11/ 12/ 5/ 8
   } { 
       \node[Bullet] at (\j-1,-\i+1){}; 
       \node[Bullet] at (\i-1,-\j+1){}; 
   }
   \foreach \i\k in {
       1/2, 2/1, 3/9, 5/12, 4/6, 6/11, 7/4, 8/7, 9/3, 10/10, 11/5, 12/8} 
       \node[font=\footnotesize] at (\i-1,-\i+1) {$\i$};  
   \draw[thick] (-0.5,0.5) rectangle (11.5,-11.5);

   \draw[dashed] (0.5,-11.5) -- (0.5, 0.5);
   \draw[dashed] (-0.5, -0.5) -- (11.5, -0.5);

   \draw[dashed] (2.5,-11.5) -- (2.5, 0.5);
   \draw[dashed] (-0.5, -2.5) -- (11.5, -2.5);

   \draw[dashed] (4.5,-11.5) -- (4.5, 0.5);
   \draw[dashed] (-0.5, -4.5) -- (11.5, -4.5);

   \draw[dashed] (5.5,-11.5) -- (5.5, 0.5);
   \draw[dashed] (-0.5, -5.5) -- (11.5, -5.5);

   \draw[dashed] (7.5,-11.5) -- (7.5, 0.5);
   \draw[dashed] (-0.5, -7.5) -- (11.5, -7.5);

   \draw[dashed] (9.5,-11.5) -- (9.5, 0.5);
   \draw[dashed] (-0.5, -9.5) -- (11.5, -9.5);

   \draw[dashed] (10.5,-11.5) -- (10.5, 0.5);
   \draw[dashed] (-0.5, -10.5) -- (11.5, -10.5);

\end{tikzpicture}
\end{minipage}
\hspace*{\fill}
\begin{minipage}{.45\linewidth}
\centering
\begin{tikzpicture}[xscale=1.5, yscale=1.2]
   \tikzset{VertexStyle/.style = {shape = rectangle, draw, 
       minimum width = 0.8cm, minimum height = 1.1em, inner sep = 2}}
   \foreach \v/\k/\x/\y in {
       1/ {1}    / 0.0/ 0.0,
       2/ {2, 3} / 1.0/ 0.0,
       4/ {4, 5} / 0.5/ 1.0,
       6/ {6}      / 1.8/ 0.0,
       7/ {7, 8}   / 2.5/ 0.0,
       9/ {9, 10}  / 3.2/ 0.0,
      11/ {11}     / 2.5/  1, 
      12/ {12}     / 1.5/  2}
      \node[VertexStyle, font=\small](\v) at (\x,\y){\small $\k$};
   \foreach \i/\j in {1/4, 2/4, 6/11, 7/11, 9/11, 11/12, 4/12} 
      \draw (\i)--(\j);
\end{tikzpicture}
\end{minipage}
\hspace*{\fill}
\caption{Fundamental supernode partition and supernodal 
elimination tree for the example in Figure~\ref{f-arrow}. The 
representative vertices are $1$, $2$, $4$, $6$, $7$, $9$, $11$, $12$.}
\label{f-supernodal}
\end{figure}

The elimination tree for a trivially perfect elimination ordering 
can be compressed into a \emph{supernodal} elimination tree, 
in which the nodes of the elimination tree are combined into larger 
supernodes.
Each supernode is associated with a \emph{representative} vertex.  
The representative vertices are the leaf nodes in the elimination tree
and all the nodes with more than one child.
The supernode with representative vertex $v$ contains the representative 
vertex $v$ itself plus the nodes in the elimination tree
between $v$ and the first ancestor $w$ that is also a 
representative vertex.  
In the supernodal elimination tree, the supernode with representative 
vertex $w$ is the parent of the supernode with representative~$v$.
The supernodes therefore form a partition of the vertex set.
Each supernode induces a complete subgraph. The vertices
in a supernode are adjacent to all vertices in the supernodes that
are its ancestors or descendants in the supernodal elimination tree.
The definitions are illustrated in Figure~\ref{f-supernodal} 
for the example in Figure~\ref{f-arrow}.  
Note that several other definitions of supernodes exist in the
sparse matrix literature. The supernodes as defined here are known as
\emph{fundamental supernodes} \cite{LNP:93}.

To conclude we summarize the properties of the
example in Figures~\ref{f-arrow} and \ref{f-supernodal} that generalize 
to arbitrary homogeneous chordal sparsity patterns in $\SS^N$.
After applying a symmetric reordering one can assume that the numerical 
ordering $1,2, \ldots, N$ is a trivially perfect elimination ordering
and a postordering.  A matrix with a homogeneous chordal sparsity
pattern will then have the form
\begin{equation} \label{e-X-supernodal}
 X = \left[\begin{array}{ccccc}
    X_{\beta_1\beta_1} & 0 & \cdots & 0  & X_{\beta_1\nu} \\
    0 & X_{\beta_2\beta_2} & \cdots & 0 & X_{\beta_2\nu} \\
    \vdots  & \vdots & \ddots & \vdots & \vdots \\
    0 & 0 &  \cdots & X_{\beta_k\beta_k} & X_{\beta_k\nu} \\
    X_{\nu\beta_1} & X_{\nu\beta_2} & \cdots & X_{\nu\beta_k}
    & X_{\nu\nu} \end{array}\right],
\end{equation}
where each of the diagonal blocks $X_{\beta_i\beta_i}$, 
for $i=1,\ldots,k$, has a block-arrow structure of the same form.  
If the sparsity pattern is not block-diagonal, $\nu$ is the supernode
at the root of the supernodal elimination tree.  Assume the root $\nu$
has $k$ children, denoted by $\nu_1$, \ldots, $\nu_k$.
Then the index set $\beta_i$ is the union of the supernode $\nu_i$ and 
its descendants in the supernodal elimination tree. 
The postordering property implies that each of these index sets 
$\beta_i$ contains consecutive indices, that precede the indices in $\nu$, 
so the matrices $X_{\beta_i\beta_i}$ are diagonal blocks.
Each of the matrices $X_{\beta_i\beta_i}$ has a homogeneous 
chordal sparsity pattern, with supernodal elimination tree given by the 
subtree rooted at $\nu_i$. 
We have assumed that the entire sparsity pattern is not block-diagonal
($\nu$ is not empty).
If it is block-diagonal, the associated sparsity graph is not connected, 
and the supernodal elimination tree is a forest with connected components
$\beta_1$, \ldots, $\beta_k$.

It is easily verified that the matrix~(\ref{e-X-supernodal}) is
positive definite, then its Cholesky factor in $X = LL^\top$ is
structured as
\begin{equation} \label{e-L-angular}
 L = \left[\begin{array}{ccccc}
    L_{\beta_1\beta_1} & 0 & \cdots & 0  & 0 \\
    0 & L_{\beta_2\beta_2} & \cdots & 0 & 0 \\
    \vdots  & \vdots & \ddots & \vdots & \vdots \\
    0 & 0 &  \cdots & L_{\beta_k\beta_k} & 0 \\
    L_{\nu\beta_1} & L_{\nu\beta_2} & \cdots & L_{\nu\beta_k}
    & L_{\nu\nu} \end{array}\right], 
\end{equation}
where each block $L_{\beta_i\beta_i}$ is the Cholesky factor
of $X_{\beta_i\beta_i}$ and therefore has a similar angular
sparsity pattern.

\subsection{Homogeneous chordal extension}
Homogeneous chordal patterns in the reordered form~(\ref{e-X-supernodal})
have a long history in many areas, including least squares 
fitting \cite[\S6.3]{Bjo:96} \cite{GoP:80}, decomposition methods in 
optimization~\cite{Las:02}, and graphical statistical 
models \cite{PeW:94,CDR:07,DrR:08,LetacMassam2007}.
The term \emph{nested block-angularity} is used in \cite[p.24]{Sau:72}.

They also arise naturally as \emph{extensions} of general unstructured
sparsity patterns, reordered using a nested dissection ordering
\cite{DER:17,GeL:81}.
Here, $\nu$ is the vertex separator in the first dissection step;
the other non-leaf supernodes are the separators in subsequent levels
of dissection.  Such a pattern is a homogeneous chordal pattern if at each 
level we treat the last block row and column in~(\ref{e-X-supernodal}) as 
dense, and also treat the principal blocks indexed by the leaves of the
supernodal elimination tree as dense.   
In applications to linear equations the matrix will have a 
large number of additional zeros within these blocks, so the actual
sparsity pattern is an unstructured sparsity pattern $E'$
(or a non-homogeneous chordal sparsity pattern if it is the filled 
pattern of a Cholesky factor), and the homogeneous chordal
pattern $E$ is an extension ($E' \subseteq E$).

When used in the nonsymmetric formulation~(\ref{e-sdp-not-symmetric})
of a sparse semidefinite program, a homogeneous chordal extension can 
be obtained by applying nested dissection to the aggregate sparsity 
pattern of $A_1, \ldots, A_m, B$.
If the homogeneous chordal extension is used to define $\mathcal V$, 
then, as we will see in the next section, the cone $K$
is a homogeneous convex cone.
The coefficient matrices $A_1, \ldots, A_m, B$ are sparse matrices
in $\mathcal V$, but their zeros within the homogeneous chordal 
pattern are not exploited in the definition of the cone $K$. 

Nested dissection ordering provides a heuristic for obtaining
homogeneous chordal extensions, with no guarantee of optimality. 
As it was proved by Yannakakis~\cite{Yannakakis1981},
given a sparsity pattern, it is NP-hard to compute the minimum number of
edges to add to make the underlying graph chordal. 
Analogously, El-Mallah and Colbourn~\cite{ElColbourn1988}
proved that given a sparsity pattern, it is NP-hard to find the smallest number of edges
to add to the graph to make it a co-graph (a graph that does not contain $P_4$ as an
induced subgraph). We can
show that given a sparsity pattern, it is NP-hard to find the largest
induced subgraph which is homogeneous chordal.

\begin{proposition}
Given a graph $G=(V,E)$ describing the sparsity pattern of a
symmetric matrix, it is NP-hard to compute the largest 
principal submatrix with homogeneous chordal sparsity pattern.
\end{proposition}

\begin{proof}
We use Theorem 3 of \cite{Bartholdi1982} (whose proof relies on 
Yannakakis's related results). This theorem states that 
given a square matrix $A$ with 0,1 entries, and a positive integer $k$,
it is NP-hard to decide whether $A$ has a $k$-by-$k$ principal 
submatrix satisfying property $\cP$, provided 
\begin{itemize}
\item
property $\cP$ is \emph{nontrivial} (meaning that it holds for infinitely
many 0,1 matrices and it fails for infinitely many 0,1 matrices);
\item
property $\cP$ holds for identity matrices;
\item
property $\cP$ is hereditary on principal submatrices.
\end{itemize}
Thus, it suffices for us to check that the property of homogeneous
chordal sparsity satisfies these required conditions. Using the excluded
induced subgraph characterization of homogeneous chordal graphs,
we note that identity matrices
correspond to empty (no edges) graphs which are homogeneous chordal;
sparsity patterns of principal submatrices correspond to induced 
subgraphs and if the original graph does not contain a $C_4$ or $P_4$ 
neither does any of its induced subgraphs. 
Finally, there are infinitely many graphs which do not contain a 
$C_4$ or $P_4$; as well, there are infinitely many graphs which
do contain either a $C_4$ or a $P_4$ (possibly both and many copies). 
Thus, homogeneous chordal sparsity satisfies the assumptions of Theorem 3 
of \cite{Bartholdi1982} and the underlying problem is NP-hard. 
\end{proof}

Therefore, one has to rely on heuristic algorithms in general 
(including polynomial time approximation algorithms for the 
minimum fill-in problems~\cite{NSS2000}), as in the
approaches used in applications of chordal extensions of sparsity 
patterns.

\section{Homogeneous sparse matrix cones} \label{s-cone}
We now apply the results of the previous section to derive properties
of the two matrix cones 
\begin{equation} \label{e-K-Kstar}
 K := \SS^N_E \cap \SS^N_+, \qquad K^* = \Pi_E(\SS^N_+).
\end{equation}
The cone $K$ is the cone of positive semidefinite matrices with 
sparsity pattern $E$.
The dual cone $K^*$ is the cone of positive semidefinite completable
matrices with sparsity pattern~$E$.  Note that $K\subseteq K^*$.
We assume that $E$ is a homogeneous chordal sparsity pattern and that the
numerical order $1,\ldots,N$ is a trivially perfect elimination ordering,
as in the example of Figure~\ref{f-arrow}.

The \emph{automorphism group} $\aut(K)$ of a regular cone $K$ is the set 
of nonsingular linear transformations that map $K$ to itself.  
A regular cone $K$ is called \emph{homogeneous} if for
every pair of points $x, y \in \inte(K)$ there exists an automorphism
of $K$ that maps $x$ to $y$.
So, a regular cone $K$ is homogeneous if and only if the automorphism 
group of $K$ acts transitively in the interior of $K$.
A subset $\mathcal  H \subseteq \aut(K)$ is a \emph{transitive subset}
of $\aut(K)$ if for every pair of points $x,y \in \inte(K)$ there exists 
an automorphism in $\mathcal H$ that maps $x$ to $y$.

Ishi \cite[Theorem A]{Ishi2013} proves that the sparse matrix
cones~(\ref{e-K-Kstar}) are homogeneous if and only if $E$ is a 
homogeneous chordal sparsity pattern.  
In this section we describe transitive subsets of the primal and dual 
automorphism groups.

\subsection{Computations with sparse triangular matrices}
\label{s-linalg}
The properties of homogeneous chordal sparse matrices that will 
be needed follow from four facts presented in the
next theorem.

\begin{theorem}
\label{thm:1.1}\ 
Let $E$ be a homogeneous chordal sparsity pattern in $\SS^N$,
with trivially perfect elimination ordering $1,\ldots, N$, 
and assume $L\in\TT^N_E$.
\begin{enumerate}[label = \arabic{enumi}.]
\item If $\tilde L \in \TT^N_E$, then $L\tilde L \in \TT^N_E$.
\item If $L$ is nonsingular, then $L^{-1} \in \TT^N_E$.
\item If $X\in \SS^N_E$, then $LXL^\top \in \SS^N_E$.
\item If $Y\in \SS^N$, then 
 $\Pi_E(L^\top YL) = \Pi_E(L^\top \Pi_E(Y)L)$.
\end{enumerate}
\end{theorem}

The second property appears in \cite{KhR:12}.
None of the four properties holds for general chordal sparsity
patterns, as can be seen by considering the example of a tridiagonal 
pattern, which is chordal but not homogeneous if $N\geq 4$. 
We also note the assumption of a trivially perfect elimination ordering.
In the example on the right in Figure~\ref{f-vinberg}, the
ordering $\sigma(1) = 1$, $\sigma(2)=3$, $\sigma(3)=2$ is a perfect 
elimination ordering and results in a zero-fill bidiagonal Cholesky 
factor.  However the inverse Cholesky factor will generally have a 
nonzero entry in position $1,2$.

\begin{proof}
To simplify the notation we denote the set $\madj(i)$ by $\alpha_i$.   
This is the set of row indices of the lower-triangular nonzeros in 
column~$i$. 
The set $\{i\} \cup \alpha_i$ is denoted by $\bar\alpha_i$.
If the order of the elements in $\alpha_i$ and $\bar\alpha_i$ matters, 
it is assumed that they are sorted in increasing order.
In this notation, the property~(\ref{e-tree-peo-tp}) can be expressed as 
\begin{equation} \label{e-col-tp}
    \alpha_i = \bar\alpha_{p(i)} \quad 
 \mbox{for all $i$,}
\end{equation}
where we interpret $\bar\alpha_{p(i)}$ as the empty  set
if $i$ is a root of the elimination tree.
In the example of Figure~\ref{f-arrow}, 
$\alpha_3 = \{4,5,12\}$, $\bar\alpha_3 = \{3,4,5,12\}$, and $p(3)=4$.

To prove the first property, we examine the sparsity pattern of 
$L\tilde L$.  The $ij$ element, with $i\geq j$, is
\[
(L\tilde L)_{ij} = \sum_{k=j}^N L_{ik} \tilde L_{kj} =  
 \sum_{k\in\bar \alpha_j} L_{ik} \tilde L_{kj}.
\]
The simplification in the second expression follows
because $\tilde L_{kj} = 0$ for $k\not\in \bar\alpha_j$.
Since $L_{ik}$ is zero if $i\not \in \bar \alpha_k$,
we have $(L\tilde L)_{ij} = 0$ for 
$i\not \in \bigcup_{k\in \bar\alpha_j} \bar\alpha_k$.
It follows from~(\ref{e-col-tp}) that 
$\bigcup_{k\in \bar\alpha_j} \bar\alpha_k = \bar \alpha_j$.
We conclude that the nonzeros of column $j$ of $L\tilde L$
are in the positions indexed by $\bar\alpha_j$, i.e., $L\tilde L
\in \TT^N_E$.

For property~2, we consider the forward substitution method 
for computing column $k$ of $L^{-1}$.  
To solve $Lx =e_k$, where $k$ is the the $k$th unit vector,
we set $x = e_k$  and run the iteration
\[
 \left[\begin{array}{c} x_j \\ x_{\alpha_j} \end{array}\right]
 :=  \left[\begin{array}{cc} 1/L_{jj} & 0 \\
 -L_{\alpha_jj}/L_{jj} & I \end{array}\right]
 \left[\begin{array}{c} x_j \\ x_{\alpha_j} 
 \end{array}\right], \quad  j =  k, k+1, \ldots, N.
\]
Since initially $x = e_k$, and $\alpha_j$ is the set of ancestors of 
vertex $j$ in the elimination tree, 
the iteration only modifies entries of $x$
on the path between $k$ and the root of the tree.
In other words, the iteration  can be simplified as
\[
 \left[\begin{array}{c} x_j \\ x_{\alpha_j} \end{array}\right]
 :=  \left[\begin{array}{cc} 1/L_{jj} & 0 \\
 -L_{\alpha_jj}/L_{jj} & I \end{array}\right]
 \left[\begin{array}{c} x_j \\ x_{\alpha_j} 
 \end{array}\right], \quad  j =  k, p(k), p^2(k), \ldots,
\]
where $p^2(k) = p(p(k))$, et cetera, i.e., we iterate over 
$j\in\bar \alpha_k$ in ascending order.
After completing the iteration, the nonzeros of $x$ are in the positions
indexed by $\bar\alpha_k$.  Therefore $L^{-1} \in\TT^N_E$.    

Next we prove property~3.  Consider the following expression for the
lower-triangular entry of $LXL^\top$ in position $ij$, with $i>j$:
\begin{equation} \label{e-prop3-pf}
(LXL^\top)_{ij}
= \sum_{k=1}^N \left(L_{ik} L_{kj} X_{kk}
 + \sum_{l\in\alpha_k} X_{lk} (L_{il}L_{jk} + L_{ik} L_{jl})\right).
\end{equation}
Suppose $i\not\in\alpha_j$, i.e.,
$i$ is not an ancestor of $j$ in the elimination tree.
We show that $(LXL^\top)_{ij} = 0$.
The first term in the sum~(\ref{e-prop3-pf}) is zero
because $L_{ik}L_{kj}\neq 0$ only if $i\in\bar\alpha_k$
and $k\in\bar\alpha_j$, which implies $i$ is on the path from
vertex $j$ to the root.
The second term is zero because $L_{il}L_{jk} \neq 0$
implies $i\in\bar\alpha_l\subset \bar \alpha_k$ and 
$j\in\bar\alpha_k$, so $i$ and $j$ are both on the path from vertex $k$
to the root, and since $i>j$, vertex $i$ is an ancestor of $j$.
Similarly, the last term is zero  because $L_{ik} L_{jl} \neq 0$
implies that $i\in\bar\alpha_k$ and $j\in\bar\alpha_l \subset \alpha_k$,
so $i$ and $j$ are both on the path from vertex $k$ to the root
and $i$ is an ancestor of $j$.

The last property in the list follows from the 3rd property.
It is sufficient to show that $\Pi_E(LYL^\top) = 0$ 
whenever $\Pi_E(Y) =0$. 
To see this, we choose any $X\in\SS^N_E$ and note that
\[
\Tr(X \Pi_E(L^\top YL)) = \Tr(XL^\top YL) = \Tr(LXL^\top Y) = 0
\]
because $LXL^\top \in\SS^N_E$ by property~3 and $\Pi_E(Y) = 0$.
\end{proof}

The properties in Theorem~\ref{thm:1.1} are also easily verified by 
induction for a pattern in the postordered block-matrix 
form~(\ref{e-X-supernodal}).  
To verify property~2, we note that if $L$ in~(\ref{e-L-angular}) is 
invertible, its inverse is 
\[
 L^{-1} = \left[\begin{array}{ccccc}
    L_{\beta_1\beta_1}^{-1} & 0 & \cdots & 0  & 0 \\
    0 & L_{\beta_2\beta_2}^{-1} & \cdots & 0 & 0 \\
    \vdots  & \vdots & \ddots & \vdots & \vdots \\
    0 & 0 &  \cdots & L_{\beta_k\beta_k}^{-1} & 0 \\
    -L_{\nu\nu}^{-1}L_{\nu\beta_1}L_{\beta_1\beta_1}^{-1} & 
    -L_{\nu\nu}^{-1}L_{\nu\beta_2}L_{\beta_2\beta_2}^{-1} & \cdots & 
    -L_{\nu\nu}^{-1}L_{\nu\beta_k}L_{\beta_k\beta_k}^{-1}
    & L_{\nu\nu}^{-1} \end{array}\right],
\]
and it is clear that $L^{-1}$ has the same sparsity pattern as $L$.

\subsection{Primal cone automorphisms} \label{s-primal-auto}
We now show that the linear transformations of the form
\begin{equation} \label{e-mL}
   \mL(X) = LXL^\top,
\end{equation}
with nonsingular $L\in\TT^N_E$, form a transitive subset of $\aut(K)$.
Property~3 in Theorem~\ref{thm:1.1} shows that $\mL(X) \in \SS^N_E$ 
for $X\in\SS^N_E$.  Since $L^{-1}\in \TT^N_E$ (by property~2), 
the same is true for the inverse mapping $\mL^{-1}(X) = L^{-1}XL^{-\top}$.
The two transformations $\mL$ and $\mL^{-1}$ preserve positive 
definiteness, so they are automorphisms for~$K$.
To show that the transformations $\mL$ form a transitive subset, 
we show that for every pair of matrices $X_1,X_2\in \inte(K$) there 
exists a nonsingular $L\in \TT^N_E$ such that $LX_1L^\top = X_2$.
Let $L_1, L_2 \in \TT^N_E$ be the triangular factors 
in the Cholesky factorizations $X_1 = L_1L_1^\top$ and $X_2=L_2L_2^\top$. 
The matrix $L = L_2 L_1^{-1}$ is nonsingular and 
in $\TT^N_E$ (by the first two properties in Theorem~\ref{thm:1.1}).
The automorphism $\mL$ defined by $L$ maps $X_1$ to $X_2$:
\[
 \mL(X_1) = LX_1L^\top = LL_1L_1^{\top}L^\top = L_2L_2^\top = X_2.
\]
We will use the notation $\auttr(K)$
for the transitive subset of $\aut(K)$ containing the 
transformations of the form~(\ref{e-mL}) with nonsingular $L \in \TT^N_E$.

\subsection{Dual cone automorphisms} \label{s-dual-auto}
The adjoint of $\mL$ is the linear mapping from $\SS^N_E$ to $\SS^N_E$
that satisfies $\langle \mL^*(S), X\rangle = \langle S, \mL(X)\rangle$
for all $S,X\in \SS^N_E$.  Since we use the trace inner product,
\[
 \langle S, \mL(X)\rangle = \Tr(SLXL^\top) = \Tr(L^\top SLX)  = 
 \langle \Pi_E(L^\top SL), X\rangle,
\]
so the adjoint is given by
\begin{equation} \label{e-mL*}
 \mL^*(S) = \Pi_E(L^\top SL).
\end{equation}
The projection in the expression $\Pi_E(L^\top SL)$ cannot be omitted 
because, unlike for the forward mapping $LXL^\top$, the product 
$L^\top SL$ is not necessarily in $\SS^N_E$. 

The linear transformations of the form $\mL^*$, where 
$\mL\in\auttr(K)$, form a transitive subset of $\aut(K^*)$.
The fact that $\mL^*$ is an automorphism of $K^*$ follows directly
from being the adjoint of an automorphism of $K$:
\begin{eqnarray*}
S \in K^* 
& \Longleftrightarrow & 
\langle S, X \rangle \geq 0 \quad \mbox{for all $X\in K$} \\
& \Longleftrightarrow & 
\langle S, \mL(X) \rangle \geq 0 \quad \mbox{for all $X\in K$} \\
& \Longleftrightarrow & 
\langle \mL^*(S), X \rangle \geq 0 \quad \mbox{for all $X\in K$} \\
& \Longleftrightarrow & 
\mL^*(S) \in K^*.
\end{eqnarray*}
On line~2 we use the fact that $\mL$ is an automorphism of $K$.
Next we prove that the mappings $\mL^*$  form a transitive subset 
of $\aut(K^*)$, by showing how for every $S_1, S_2\in\inte(K^*)$
one can find $L$ such that $\mL^*(S_1) = S_2$.
We use a classical result from the theory of
positive definite matrix completions,
stating that for every $S\in\inte(K^*)$ there exists an 
$X\in \inte(K)$ that satisfies $\Pi_E(X^{-1}) = S$
\cite{GroneJohnsonSa1984}.  
The matrix $X$ is the inverse of the maximum-determinant positive 
definite completion, i.e., the unique solution $Y$ of the convex 
optimization problem
\begin{equation} \label{e-maxdet}
\begin{array}{ll}
 \mbox{minimize} & -\ln\det(Y) \\
 \mbox{subject to} & \Pi_E(Y) = S 
 \end{array}
\end{equation}
over $Y\in\SS^N_{++}$.
The optimality conditions for this problem,
\[
 Y^{-1} = X \succ 0, \qquad \Pi_E(Y) = S, 
\]
where $X\in\SS^N_E$ is a multiplier for the equality constraint 
of~(\ref{e-maxdet}), show that $\Pi_E(X^{-1})= S$.
Now consider two matrices $S_1, S_2\in \inte(K^*)$.
To construct an automorphism $\mL^*$ (of $K^*$) that maps $S_1$ to $S_2$,
we compute the matrices $X_1, X_2 \in \inte(K)$ that satisfy
$\Pi_E(X_1^{-1}) = S_1$, $\Pi_E(X_2^{-1}) = S_2$. 
Let $L_1, L_2 \in\TT^N_E$ be the Cholesky factors of $X_1$ and $X_2$,
and define $L = L_1L_2^{-1}$.  Then
\begin{eqnarray*}
\mL^*(S_1) 
& = & \Pi_E(L^\top S_1L) \\
& = & \Pi_E(L^\top \Pi_E(L_1^{-\top}L_1^{-1})L) \\
& = & \Pi_E(L^\top L_1^{-\top}L_1^{-1}L) \\
& = & \Pi_E(L_2^{-\top}L_2^{-1}) \\
& = & S_2.
\end{eqnarray*}
On line~3 we apply property~4 in Theorem~\ref{thm:1.1}.

\subsection{Matrix inverse}
The inverse of a positive definite matrix $X\in \inte(K)$ can be
factorized as $X^{-1} = RR^\top$ where the upper-triangular
matrix $R = L^{-\top}$  is sparse and satisfies $R^\top \in \TT^N_E$.
Suppose the pattern is in the postordered block-matrix 
form~(\ref{e-X-supernodal}).
Then 
\begin{eqnarray*}
R & = & \left[\begin{array}{ccccc}
    R_{\beta_1\beta_1} & 0 & \cdots & 0  & R_{\beta_1\nu} \\
    0 & R_{\beta_2\beta_2} & \cdots & 0 &  R_{\beta_2\nu} \\
    \vdots  & \vdots & \ddots & \vdots & \vdots \\
    0 & 0 &  \cdots & R_{\beta_k\beta_k} & R_{\beta_k\nu} \\
    0 & 0 & \cdots & 0 & R_{\nu\nu} \end{array}\right]  \\
 & = & \left[\begin{array}{ccccc}
    L^{-\top}_{\beta_1\beta_1} & 0 & \cdots & 0  & 
    -L_{\beta_1\beta_1}^{-1}L_{\nu\beta_1}^\top L_{\nu\nu}^{-1} \\
    0 & L_{\beta_2\beta_2}^{-\top} & \cdots & 0 &  
    -L_{\beta_2\beta_2}^{-1}L_{\nu\beta_2}^\top L_{\nu\nu}^{-1} \\
    \vdots  & \vdots & \ddots & \vdots & \vdots \\
    0 & 0 &  \cdots & L_{\beta_k\beta_k}^{-\top} & 
    -L_{\beta_k\beta_k}^{-1}L_{\nu\beta_k}^\top L_{\nu\nu}^{-1} \\
    0 & 0 & \cdots & 0 & L_{\nu\nu}^{-\top} \end{array}\right]
\end{eqnarray*}
and $X^{-1}$ is the sum of a block-diagonal and a low-rank matrix
\begin{eqnarray*}
X^{-1} & = & 
\left[\begin{array}{cccc}
X_{\beta_1\beta_1}^{-1} & \cdots & 0 & 0 \\
\vdots & \ddots & \vdots & \vdots \\
0 & \cdots & X_{\beta_k\beta_k}^{-1} & 0 \\
0 & \cdots & 0 & 0 \end{array}\right] + 
\left[\begin{array}{c}
R_{\beta_1\nu} \\ \vdots \\ R_{\beta_k\nu} \\ 
R_{\nu\nu} \end{array}\right] \left[\begin{array}{c}
R_{\beta_1\nu} \\ \vdots \\ R_{\beta_k\nu} \\ R_{\nu\nu} 
\end{array}\right]^\top.
\end{eqnarray*}
Moreover each diagonal block $X_{\beta_i\beta_i}^{-1}$ has a similar
block-diagonal plus low-rank structure.

Conversely, consider a block-diagonal plus low-rank matrix
\[
 Y =
\left[\begin{array}{cccc}
Y_{\beta_1\beta_1} & \cdots & 0 & 0 \\
\vdots & \ddots & \vdots & \vdots \\
0 & \cdots & Y_{\beta_k\beta_k} & 0 \\
0 & \cdots & 0 & 0 \end{array}\right] + 
\left[\begin{array}{c}
W_{\beta_1\nu} \\ \vdots \\ W_{\beta_k\nu} \\ W_{\nu\nu} 
 \end{array}\right] 
\left[\begin{array}{c}
W_{\beta_1\nu} \\ \vdots \\ W_{\beta_k}\nu \\ W_{\nu\nu} 
\end{array}\right]^\top
\]
where the matrices $Y_{\beta_1\beta_1}$,
\ldots, $Y_{\beta_k\beta_k}$ are positive definite, and $W_{\nu\nu}$
is invertible.
Then the inverse is a block-arrow matrix
\[
Y^{-1} = 
\left[\begin{array}{cccc}
Y_{\beta_1\beta_1}^{-1} & \cdots & 0 & 
    -Y_{\beta_1\beta_1}^{-1}W_{\beta_1\nu} W_{\nu\nu}^{-1} \\
\vdots & \ddots & \vdots & \vdots \\
0 & \cdots & Y_{\beta_k\beta_k}^{-1} & 
    -Y_{\beta_k\beta_k}^{-1} W_{\beta_k\nu} W_{\nu\nu}^{-1} \\
-W_{\nu\nu}^{-1}W_{\beta_1\nu}^{\top}Y_{\beta_1\beta_1}^{-1} & \cdots & 
-W_{\nu\nu}^{-1}W_{\beta_k\nu}^{\top}Y_{\beta_k\beta_k}^{-1} & 
W_{\nu\nu}^{-\top} S W_{\nu\nu}^{-1} \end{array}\right]
\]
where 
\[
S = I + \sum_{i=1,\ldots,k} 
W_{\beta_i\nu}^TY_{\beta_i\beta_i}^{-1}W_{\beta_i\nu}.
\]

\section{Logarithmic barriers} \label{s-barriers}
The function $-\ln\det(X)$ for symmetric positive definite $X$ has 
important applications in statistics, machine learning, information 
theory, and semidefinite optimization.
Here, we restrict the function to the symmetric matrices with a given
homogeneous chordal sparsity pattern $E$.  We denote this function by
$F: \SS^{N}_E \to (-\infty, +\infty]$,
\begin{equation} \label{e-log-det-barrier}
F(X) := \left\{ \begin{array}{rl} -\ln\det(X), 
& \mbox{ if }X \in \inte(K)\\
+\infty, & \mbox{ otherwise,}
\end{array}
\right.
\end{equation}
where $K$ is the primal cone in~(\ref{e-K-Kstar}),
and refer to $F$ as the \emph{logarithmic barrier} for $K$.

The gradient and Hessian of $F$ (as a function on $\SS^N_E$) 
at $X\in\inte(K)$ are given by
\begin{equation} \label{e-F-derivatives} 
 F'(X) = -\Pi_E(X^{-1}), \qquad
 F''(X;Y) = \Pi_E(X^{-1}YX^{-1}).
\end{equation}
Here $F''(X;Y)$ denotes the directional derivative of $F'$ at 
$X$ in the direction $Y\in \SS^N_E$, i.e.,
\[
 F''(X;Y) = \left. \frac{d}{d\alpha} F'(X+\alpha Y) \right|_{\alpha = 0}.
\]
The \emph{conjugate barrier} of $F$ is defined as
\[
 F_*(S) =  \sup_{X\in\inte(K)} \left\{-\langle S,X\rangle - F(X)\right\}
\]
and has domain $\inte(K^*)$.  This is the logarithmic barrier for $K^*$.
The maximizer in the optimization problem in the definition
is the positive definite solution $\hat X$ of the nonlinear equation 
\[
 F'(X) = -\Pi_E(X^{-1}) = -S,
\]
with variable $X\in\SS^N_E$.  
The inverse $\hat X^{-1}$ of the solution is the maximum-determinant 
positive definite completion of $S$.  From 
$\hat X$ we obtain the function value 
$F_*(S) = - F(\hat X) - N$ and the derivatives 
\begin{equation} \label{e-F*-derivatives} 
 F_*'(S) = -\hat X, \qquad F_*''(S) = F''(\hat X)^{-1}.
\end{equation}

In this section, we derive some interesting properties of compositions
of $F$ and $F_*$ with the cone automorphisms~(\ref{e-mL}) 
and~(\ref{e-mL*}), respectively.

\subsection{Composition with primal cone automorphism}
As in Section~\ref{s-cone}, we assume that the numerical order
is a trivially perfect elimination ordering for $E$.  Clearly,
\begin{equation} \label{e-F-scaling}
 F(\mL(X)) = F(LXL^T) = F(X) + F(LL^\top) 
\end{equation}
for all $X\in\inte(K)$  and nonsingular $L\in\TT_E$.   
Differentiating the left- and right-hand sides with respect to $X$ 
shows that 
\begin{equation} \label{e-composition}
 F'(\cL(X)) = \cL^{-*}(F'(X)), \qquad
 F''(\cL(X)) = \cL^{-*} \circ F''(X) \circ \cL^{-1}
\end{equation}
for all $X\in\inte(K)$ and nonsingular $L\in \TT_E$.
These properties can also be verified from the 
definitions~(\ref{e-F-derivatives}) and Theorem~\ref{thm:1.1}.
For the gradient,
\begin{eqnarray*}
F'(\cL(X))
 & = &  -\Pi_E{((LXL^\top)^{-1})} \\
 & = &  -\Pi_E{(L^{-\top}X^{-1}L^{-1})} \\
 & = &  -\Pi_E{(L^{-\top}\Pi_E(X^{-1})L^{-1})} \\
 & = &  \cL^{-*}(F'(X)).
\end{eqnarray*}
On line 3 we use property~4 in Theorem~\ref{thm:1.1}.
The result for the Hessian follows similarly from
\begin{eqnarray*}
F''(\cL(X); Y)
 & = & \Pi_E{\left((LXL^{\top})^{-1} Y (LXL^{\top})^{-1}\right)} \\
 & = & \Pi_E{(L^{-\top}X^{-1}L^{-1} Y L^{-\top} X^{-1}L^{-1})}\\
 & = & \Pi_E{(L^{-\top} \Pi_E{(X^{-1}L^{-1} Y L^{-\top} X^{-1})} L^{-1})}\\
 & = & \cL^{-*}(F''(X; \cL^{-1}(Y)))
\end{eqnarray*}
for every $Y\in \SS_E$.

\subsection{Composition with dual cone automorphism}
Similar properties hold for the dual barrier.  Using~(\ref{e-F-scaling})
in the definition of the dual barrier, we find that
\begin{eqnarray*}
 F_*(S) 
& = & \sup_X{ \{\langle -S,X\rangle - F(X)\}} \\
& = & \sup_X{ \{\langle -S,\mathcal L(X)\rangle - F(\mathcal L(X))\}} \\
& = & \sup_X{ \{\langle -\mathcal L^*(S),X\rangle - F(X)\}} - F(LL^T) \\
& = & F_*(\mathcal L^*(S)) - F(LL^T).
\end{eqnarray*} 
Hence, $F_*(\mathcal L^*(S)) = F_*(S) + F(LL^T)$
for all $S\in \inte(K^*)$ and nonsingular $L\in\TT_E$.
Differentiating with respect to $S$ shows that
\begin{equation} \label{e-composition-dual}
 F'_*(\cL^*(S)) = \cL^{-1}(F'_*(S)), 
 \qquad
 F''_*(\cL^*(S)) = \cL^{-1} \circ F''_*(S) \circ \cL^{-*}.
\end{equation}
To verify these properties directly, we note that, by definition,
\begin{eqnarray*}
\hat X = -F'_*(S) & \Longleftrightarrow &
\Pi_E{(\hat X^{-1})} = S, \\
\hat Y = -F'_*(\cL^*(S))  & \Longleftrightarrow &
\Pi_E{(\hat Y^{-1})} = \Pi_E{(L^\top SL)}.
\end{eqnarray*}
Combining the two properties, we obtain
\[
\Pi_E{(\hat Y^{-1})} 
= \Pi_E{(L^\top \Pi_E(\hat X^{-1}) L)} 
= \Pi_E{(L^\top \hat X^{-1} L)}.
\]
Since the maximum-determinant positive definite completion is unique,
we conclude that
\[
-F'_*(\cL^*(S)) = \hat Y = L^{-1}\hat XL^{-\top} = -\cL^{-1}(F'_*(S)).
\]
The Hessian property in~(\ref{e-composition-dual}) follows from
\begin{eqnarray*}
 F''_*(\cL^*(S))
& = & F''(\hat Y)^{-1} \\
& = & F''(\cL^{-1}(\hat X))^{-1} \\
& = & (\cL^* \circ F''(\hat X) \circ \cL)^{-1} \\
& = & \cL^{-1} \circ F''(\hat X)^{-1} \circ \cL^{-*} \\
& = & \cL^{-1} \circ F_*''(S) \circ \cL^{-*}.
\end{eqnarray*}

\subsection{Hessian factorization}\label{sec:4.3}
An important consequence of the second relation in~(\ref{e-composition})
is that the Hessian of $F$ at any point $X\in\inte(K)$ can be factored as
\begin{equation} \label{e-primal-Hessian-fact}
 F''(X) = \mL^{-*} \circ \mL^{-1},
\end{equation}
where $\mL\in\auttr(K)$, namely the automorphism
that maps the identity matrix $I$ to $\mL(I) = X$
(and defined by the Cholesky factor of $X$).
Similary, from~(\ref{e-composition-dual}), the Hessian of $F_*$ at
any point $S\in\inte(K^*)$ admits a factorization
\[
 F''_*(S) = \mL \circ \mL^*,
\]
where $\mL^*$ is the dual cone automorphism that maps $S$ to $\mL^*(S)= I$.

In~\cite{NT1997} (also see ~\cite[Theorem 3.1]{Tuncel2001}), it is
shown that for every $X\in\inte(K)$ and $S\in\inte(K^*)$ 
there exists a unique $W\in \inte(K)$  that satisfies
\[
 F''(W;X) = S,
\]
where $F''(W;X)$ is the directional derivative of $F'$ at $W$ in
the direction $X$. 
The matrix $W$ is the solution
of the convex optimization problem
\[
 \mbox{minimize} \quad {-\langle F'(W), X\rangle + \langle S, W\rangle}
\]
with variable $W$.  
By factorizing $F''(W)$ as $F''(W) = \mL^{-*} \circ \mL^{-1}$, we 
obtain the following theorem.  
\begin{theorem}
\label{thm:asymmetric-factor}
For every pair of interior points $X \in \inte(K)$ and $S \in \inte(K^*)$,
there exists a unique $\cL\in\auttr(K)$ 
which satisfies 
\[
 \cL^{-1}(X) = \cL^{*}(S),
\]
i.e., there exists a nonsingular $L\in \TT_E$ such that
$L^{-1} XL^{-\top} = \Pi_E(L^{\top} SL)$.
\end{theorem}

Theorem~\ref{thm:asymmetric-factor} can be generalized to all homogeneous 
cones (see the discussion following Theorem~\ref{thm:6.3}).  
Efficient computation of the matrix $W$ is a topic of current research.

A closely related result on convex cones is discussed 
in \cite{Tun:98}.  Theorem 4.2 of \cite{Tun:98} states that if there
exists a subset $G\subseteq \aut(K)$ such that for every
$x \in \inte(K)$ and $s\in\inte(K^*)$ there exists a self-adjoint
$\mathcal D\in G$ which satisfies
\[
\mathcal D^{-1}(x) = \mathcal D(s),
\]
then $K$ must be a symmetric cone (homogeneous and self-dual).
Theorem \ref{thm:asymmetric-factor} does not contradict 
Theorem 4.2 of \cite{Tun:98} because the automorphism $\mL$ 
in Theorem \ref{thm:asymmetric-factor} is not self-adjoint.

\section{Homogeneous matrix cones} \label{s-hom-matrix-cones}
As an extension of~(\ref{e-K-Kstar}) we now consider slices of the 
positive semidefinite cone
\begin{equation} \label{e-K-V}
 K := \mathcal V \cap \SS^N_+,
\end{equation}
where $\mathcal V$ is a subspace of $\SS^N$.
It is clear that $K$ is a closed, pointed, and convex cone.
We will assume that $\mathcal V \cap \SS^N_{++}$ is nonempty,
so $K$ has nonempty interior (relative to $\mathcal V$).
The corresponding dual cone (in the subspace $\mathcal V$)
is given by 
\begin{equation} \label{e-K*-V}
K^* =  \Pi_\mathcal V(\SS^N_+),
\end{equation}
where $\Pi_\mathcal V$ denotes Euclidean projection on $\mathcal V$.
To see this, we first note that the cone $\Pi_\mathcal V(\SS^N_+)$
is closed.  This follows from \cite[theorem 9.1]{Roc:70}
and the fact that if $\Pi_\mathcal V(Y)= 0$ and $Y\succeq 0$ then 
$Y$ must be zero, because $\Pi_\mathcal V(Y)=0$ implies that
$\Tr(YX) = 0$ for all $X\in \mathcal V$ and, by assumption, 
$\mathcal V$ contains positive definite matrices.
Next, it is easily verified that the dual of the cone 
$\Pi_\mathcal V(\SS^N_+)$ is given by the cone $K$ defined 
in~(\ref{e-K-V}):
\begin{eqnarray*}
\left(\Pi_{\mathcal V}(\SS^N_+)\right)^*
& = &
\{ X\in\mathcal V : \Tr(SX) \geq 0 \mbox{\ for all 
$S\in\Pi_\mathcal V(\SS^N_+)$}\} \\
& = & 
\{ X\in\mathcal V : \Tr(YX) \geq 0 \mbox{\ for all $Y \in \SS^N_+$}\} \\
& = & \mathcal V \cap \SS^N_+.
\end{eqnarray*}
Hence $K = (\Pi_\mathcal V(\SS^N_+))^*$.
Since $\Pi_\mathcal V(\SS^N_+)$ is closed,
we have 
\[
K^* = \left(\Pi_\mathcal V(\SS^N_+)\right)^{**} = \Pi_\mathcal V(\SS^N_+).
\]
We conclude that $K$ and $K^*$ form a dual pair of regular cones.
We also note that $K\subseteq K^*$.

Ishi \cite{Ishi2015} presents conditions on $\mathcal V$ that imply that
the cone $K$ defined in~(\ref{e-K-V}) is 
homogeneous.
Suppose that after a suitable reordering, the matrices $X\in\mathcal V$
can be partitioned as $r$-by-$r$ block matrices
\begin{equation} \label{e-X-partitioned}
X = 
 \left[\begin{array}{ccccc}
   X_{11} & X_{21}^\top & X_{31}^\top & \cdots & X_{r1}^\top \\*[.5ex]
   X_{21} & X_{22} & X_{32}^\top & \cdots & X_{r2}^\top \\*[.5ex]
   X_{31} & X_{32} & X_{33} & \cdots & X_{r3}^\top \\*[.5ex]
   \vdots & \vdots & \vdots & \ddots & \vdots \\*[.5ex]
   X_{r1} & X_{r2} & X_{r3} & \cdots &  X_{rr}\end{array} \right],
\end{equation}
with blocks $X_{ij}$ of size $N_i$-by-$ N_j$, and that
\begin{equation} \label{e-V-def}
 \mathcal V := \left\{ X \in \SS^N : X_{ij} \in \mathcal V_{ij}, \;
 i \in \{1,\ldots, r\}, j \in\{1,\ldots,i\} \right\}
\end{equation}
where $\mathcal V_{ii}$ is a subspace of $\SS^{N_i}$ and, for $i\neq j$, 
$\mathcal V_{ij}$ is a subspace of $\R^{N_i \times N_j}$.
For $j>i$ we define $\mathcal V_{ij} = \{ U^\top : U\in \mathcal V_{ji}\}$.
Suppose the subspaces $\mathcal  V_{ij}$ satisfy the following properties.
\begin{enumerate}[label=\textbf{P\arabic*.}, ref=P\arabic*]
\item \label{p-ishi1}
The diagonal blocks are multiples of the identity:
$\mathcal V_{ii} = \{\alpha I : \alpha \in \R\}$ for $i=1,\ldots, r$.
\item \label{p-ishi2}
 The lower-triangular blocks have orthogonal rows of equal norm:
 if $i>j$ and $A\in \mathcal V_{ij}$, then $AA^\top$ is a multiple of the
 identity.  
\item \label{p-ishi3}
 If $i>j>k$, then the subspaces $\mathcal V_{ij}$, $\mathcal V_{jk}$, 
 $\mathcal V_{ik}$ are related as follows:
 \[
 A \in \mathcal V_{ik}, \; B \in \mathcal V_{jk} \quad
   \Longrightarrow  \quad AB^\top \in \mathcal V_{ij}. 
 \]
\item \label{p-ishi4}
 If $i>j>k$, then the subspaces $\mathcal V_{ij}$,
 $\mathcal V_{jk}$, $\mathcal V_{ik}$ 
 are related as follows:
 \[
 A \in \mathcal V_{ij}, \; B \in \mathcal V_{jk}  \quad
  \Longrightarrow \quad AB\in \mathcal V_{ik}. 
 \]
\end{enumerate}
Ishi \cite[theorem 3]{Ishi2015} shows that the cone $K$ is homogeneous.
Sections~\ref{s-homogeneous-chol}--\ref{s-homogeneous-aut}
will explain this in more detail.  

Property~\ref{p-ishi1} implies that $I \in \mathcal V$, so 
$V\cap \SS^N_{++} \neq \emptyset$, as assumed at the beginning
of this section.  
A useful equivalent form of~\ref{p-ishi2} is the following:
if $i>j$ and $B, C\in\mathcal V_{ij}$, then $BC^T + CB^T$ is a multiple 
of the identity.  This follows from~\ref{p-ishi2} applied to $A=B+C$ and, 
conversely, clearly implies~\ref{p-ishi2} if we take $B=C=A$.

In the next sections we use the following notation for the set of 
lower-triangular matrices with $L+L^\top\in\mathcal V$:
\begin{equation} \label{e-T-def}
\mathcal T := \left\{ L\in\TT^N : L_{ij} \in \mathcal V_{ij}, \;
 i \in \left\{1,\ldots, r \right\}, \; j \in \left\{ 1, \ldots, i\right\} 
 \right\}. 
\end{equation}
Here $L_{ij}$ refers to the $N_i$-by-$N_j$ submatrix of $L$, 
partitioned as in~(\ref{e-X-partitioned}).

\subsection{Examples} \label{s-hom-cones-examples}

\subsubsection*{Homogeneous sparse matrix cones.}
The homogeneous sparse matrix cones of 
Sections~\ref{s-cone}--\ref{s-barriers} 
are a special case with $\mathcal V := \SS^N_E$.
Suppose $E$ is a homogeneous chordal sparsity pattern and that
the numerical order $1,\ldots, N$ is a trivially perfect
elimination ordering.  Define $r:=N$, $N_1 :=\cdots := N_r :=1$, and
\[
\mathcal V_{ij} := \left\{\begin{array}{ll}
\{0\} & \mbox{$i\neq j$ and $\{i,j\}\not\in E$} \\
\R & \mbox{otherwise.} \end{array}\right.
\]
Properties~\ref{p-ishi1} and~\ref{p-ishi2} hold trivially, 
since $N_i=1$ for all $i$.  Property \ref{p-ishi3} reduces to
\[
i > j > k, \quad \{i,k\} \in E, \quad \{j,k\} \in E 
\qquad  \Longrightarrow \qquad 
\{i,j\} \in E.
\]
This is the property~(\ref{e-peo2}) of a perfect elimination ordering 
of a chordal graph.
Property~\ref{p-ishi4} is
\[
i > j > k, \quad \{i,j\}\in E, \quad \{j,k\}\in E 
\qquad \Longrightarrow \qquad 
\{i,k\} \in E.
\]
This is the additional property~(\ref{e-peo3}) of a trivially perfect
elimination ordering.

\subsubsection*{Block-sparsity.}
As an extension, we can define a block-sparsity pattern for a
matrix partitioned as in~(\ref{e-X-partitioned}) as an 
undirected graph with vertex set $V = \{1,2,\ldots,r\}$
and edge set
\[
 E = \{ \{i,j\} : i\neq j, \mathcal V_{ij} \neq \{0\}\}.
\]
Properties~\ref{p-ishi3} and~\ref{p-ishi4} imply (among other conditions
on the subspaces) that the graph
$(V,E)$ represents an $r$-by-$r$ homogeneous chordal sparsity pattern
with trivially perfect ordering $1,\ldots, r$.
As an example, Properties~\ref{p-ishi1}--\ref{p-ishi4} are
satisfied by the subspace $\mathcal V$ of matrices of the form
\[
 \left[\begin{array}{ccc}
   \alpha I & 0 & u \\
   0 & \beta I & v \\
   u^\top & v^\top & \gamma 
 \end{array}\right]
\]
with $u\in\R^{N_1}$, $v\in\R^{N_2}$, and $\alpha, \beta,\gamma \in \R$.
The corresponding homogeneous chordal sparsity pattern is the 
$3$-by-$3$ pattern of Figure~\ref{f-vinberg}.

\subsubsection*{Rotated quadratic cone.}
The subspace
\[
\mathcal V = \left\{ \left[\begin{array}{cc} 
  \alpha I & u \\ u^\top & \beta  \end{array}\right] : 
  \alpha, \beta \in\R, \; u\in\R^{N-1}\right\}
\]
is a special case with $r=2$,  $N_1 = N-1$, $N_2=1$,
and $\mathcal V_{12} = \R^{N-1}$.
The cone $K = \mathcal V \cap \SS^N_+$ is linearly isomorphic
to the cone 
\begin{eqnarray}
\mathcal Q_\mathrm r 
& = & \left\{ (\alpha, \beta, u) \in \R \times\R \times \R^{N-1} :
 \alpha, \beta \geq 0, \; \alpha\beta \geq u^\top u \right\} 
\label{e-rot-quad}\\
& = & \left\{ (\alpha, \beta, u) \in \R \times\R \times \R^{N-1} :
\left[\begin{array}{cc} 
  \alpha I & u \\ u^\top & \beta \end{array}\right] \succeq 0 \right\}.
\nonumber
\end{eqnarray}
The cone $\mathcal Q_\mathrm r$ is known as the \emph{rotated quadratic 
cone} and is a symmetric cone.
It can be used to represent the second order cone
\[
\mathcal Q = \{ (t, y) \in \R\times \R^N : \|y\|_2 \leq t \}
\] 
as
\[
\mathcal Q = \{ (t, y) \in \R\times \R^N :
 (t+y_1, t-y_1, \bar y) \in Q_\mathrm r\}
\]
where $\bar y  = (y_2,\ldots,y_N)$.

\subsubsection*{Non-sparse example.}
Define $\mathcal V$ as the set of matrices of the form
\[
\left[\begin{array}{ccccc} 
  \alpha I & 0 & u_1 & -u_2 & v_1\\*[.3ex]
  0 & \alpha I & u_2 & u_1 & v_2 \\*[.3ex]
  u_1^\top & u_2^\top & \beta & 0 & w_1 \\*[.3ex]
 -u_2^\top & u_1^\top & 0 & \beta & w_2 \\*[.3ex]
 v_1^\top & v_2^\top & w_1 & w_2 & \gamma
\end{array}\right] 
\]
with $\alpha, \beta,\gamma, w_1, w_2\in\R$ and 
$u_1, u_2, v_1, v_2 \in\R^M$.
This is a special case with
$N=2M+3$, $r=3$, $N_1 = 2M$, $N_2 =2$, $N_3=1$, and
\[
\mathcal V_{12} =
\left\{ \left[\begin{array}{cc}
    u_1 & -u_2 \\ u_2 & u_1 \end{array}\right] :
    u_1, u_2\in\R^M \right\}, \qquad
\mathcal V_{13} = \R^{2M}, \qquad
\mathcal V_{23} = \R^2. 
\]

\subsubsection*{Matrix norm cone.}
Define $\mathcal V$ as 
\[
\mathcal V = \left\{ \left[\begin{array}{cc} 
  \alpha I & U^\top \\ U & V \end{array}\right] : t\in\R, \;
   U\in\R^{K\times L}, \; V \in\SS^K\right\}.
\]
This is a special case of Ishi's general structure, with
$r = K+1$, $N_1=L$, $N_2=\cdots = N_r =1$. 
The off-diagonal subspaces are defined as
\[
\mathcal V_{ij} = 
\left\{\begin{array}{ll}
 \R^L & i=1, \; j \in \{2,\ldots, r \}  \\
 \R & i \in \{2,\ldots, r\}, \; j \in \{2,\ldots, i-1\}.  
 \end{array}\right.
\]
The cone $K = \mathcal V \cap \SS^N_+$ is known as the 
\emph{matrix norm cone} and is important for trace norm 
minimization problems \cite{KarimiT2019,KarimiT2020}.

\subsubsection*{Sparse matrix norm cone.}
The matrix norm cones and homogeneous sparse matrix cones can be combined 
in a new class of homogeneous matrix cones.  Define $\mathcal V$ as 
\[
\mathcal V := \left\{ \left[\begin{array}{cc} 
  \alpha I & U^\top \\ U & V \end{array}\right] : \alpha \in\R, \;
   U\in\mathcal U, \; V \in\SS_E^K \right\},
\]
where $\mathcal U$ is a subspace of $\R^{K\times L}$
with the property that for every $U\in \mathcal U$,  the product
$UU^\top \in\SS^K_E$.
An example is the set of positive semidefinite
matrices of the form
\[
\left[\begin{array}{ccccc|ccc}
 \alpha & 0 & 0 & 0 & 0 & u_1 & 0 & u_6 \\
 0 & \alpha & 0 & 0 & 0 & 0 & u_4 & u_7 \\
 0 & 0 & \alpha & 0 & 0 & 0 & u_5 & u_8 \\
 0 & 0 & 0 & \alpha & 0 & u_2 & 0 & u_9 \\
 0 & 0 & 0 & 0 & \alpha & u_3 & 0 & u_{10} \\ \hline
 u_1 & 0 & 0 & u_2 & u_3  & v_1 & 0 & v_2 \\
 0  & u_4 & u_5  & 0 & 0 & 0 & v_3 & v_4 \\
 u_6 & u_7 & u_8  & u_9 & u_{10} & v_2 & v_4 & v_5 
\end{array}\right].
\]

\subsection{Cholesky factorization} \label{s-homogeneous-chol}
In this section we assume that $\mathcal V$ 
satisfies~\ref{p-ishi1}, \ref{p-ishi2}, \ref{p-ishi3} but not 
necessarily \ref{p-ishi4}. 
We show that every positive definite matrix 
$X\in \mathcal V \cap \SS^N_{++}$ has a Cholesky factorization
$X = LL^\top$ where $L\in\mathcal T$.
This is the counterpart of the zero-fill Cholesky factorization
of positive definite matrices with chordal sparsity patterns. 

\begin{proof}
The proof is by induction on $r$.
For $r=1$, we have $\mathcal V = \{ \alpha I : \alpha\in\R\}$ and
the result is obvious, with $L=\sqrt{\alpha}I$ if $X=\alpha I$. 
We show that the result holds for $r=m$ if it holds for $r=m-1$.
Suppose $X$ is positive definite of the form
\[
 X = \left[\begin{array}{ccccc}
   \alpha_1I & X_{21}^\top & X_{31}^\top & \cdots & X_{r1}^\top \\*[.5ex]
   X_{21} & \alpha_2I & X_{32}^\top & \cdots & X_{r2}^\top \\*[.5ex]
   X_{31} & X_{32} & \alpha_3 I & \cdots & X_{r3}^\top \\*[.5ex]
   \vdots & \vdots & \vdots & \ddots & \vdots \\*[.5ex]
   X_{r1} & X_{r2} & X_{r3} & \cdots &  \alpha_r I \end{array} \right] 
\]
with $X_{ij} \in \mathcal V_{ij}$ for $i \in \{2, \ldots, r\}$, 
$j \in\{1,\ldots, i-1\}$, and that the subspaces $\mathcal V_{ij}$ 
satisfy~\ref{p-ishi2} and~\ref{p-ishi3}.
The matrix can be factored as
\begin{equation} \label{e-X-chol-pf}
X = 
 \left[\begin{array}{
 @{\hskip 0.2em}c
 @{\hskip 0.7em}c
 @{\hskip 0.7em}c
 @{\hskip 0.7em}c
 @{\hskip 0.7em}c
 @{\hskip 0.2em}  
 }
 L_{11} & 0  & 0  & \cdots & 0 \\*[.5ex]
 L_{21} & I  & 0  & \cdots & 0 \\*[.5ex]
 L_{31} & 0  & I  & \cdots & 0 \\*[.5ex]
 \vdots & \vdots  & \vdots & \ddots & \vdots \\*[.5ex]
 L_{r1} & 0 & 0   & \cdots &  I \end{array} \right]  
 \left[\begin{array}{
 @{\hskip 0.2em}c
 @{\hskip 0.7em}c
 @{\hskip 0.7em}c
 @{\hskip 0.7em}c
 @{\hskip 0.7em}c
 @{\hskip 0.2em}  
 }
 I & 0  & 0  & \cdots & 0 \\*[.5ex]
 0 & Y_{22}  & Y_{32}^\top  & \cdots & Y_{r2}^\top \\*[.5ex]
 0 & Y_{32}  & Y_{33}  & \cdots & Y_{r3}^\top  \\*[.5ex]
 \vdots & \vdots  & \vdots & \ddots & \vdots \\*[.5ex]
 0 & Y_{r2} & Y_{r3}   & \cdots &  Y_{rr} \end{array} \right] 
 \left[\begin{array}{
 @{\hskip 0.2em}c
 @{\hskip 0.7em}c
 @{\hskip 0.7em}c
 @{\hskip 0.7em}c
 @{\hskip 0.7em}c
 @{\hskip 0.2em}  
 }
 L_{11} & 0  & 0  & \cdots & 0 \\*[.5ex]
 L_{21} & I  & 0  & \cdots & 0 \\*[.5ex]
 L_{31} & 0  & I  & \cdots & 0 \\*[.5ex]
 \vdots & \vdots  & \vdots & \ddots & \vdots \\*[.5ex]
 L_{r1} & 0 & 0   & \cdots &  I \end{array} \right]^\top,
\end{equation}
where $L_{11} := \sqrt{\alpha_1}I \in \mathcal V_{11}$, and 
$L_{i1} := X_{i1}/\sqrt{\alpha_1} \in \mathcal V_{i1}$ for 
$i \in \{2,\ldots, r\}$.  The matrix $Y$ is the Schur complement
\[
 \left[\begin{array}{ccccc}
 Y_{22}  & Y_{32}^\top  & \cdots & Y_{r2}^\top \\*[.5ex]
 Y_{32}  & Y_{33}  & \cdots & Y_{r3}^\top  \\*[.5ex]
 \vdots  & \vdots & \ddots & \vdots \\*[.5ex]
 Y_{r2} & Y_{r3}   & \cdots &  Y_{rr} \end{array} \right] 
 = 
 \left[\begin{array}{ccccc}
 \alpha_2 I   & X_{32}^\top  & \cdots & X_{r2}^\top \\*[.5ex]
 X_{32}  & \alpha_3 I  & \cdots & X_{r3}^\top  \\*[.5ex]
 \vdots  & \vdots & \ddots & \vdots \\*[.5ex]
 X_{r2} & X_{r3}   & \cdots & \alpha_r I \end{array} \right] 
 - \frac{1}{\alpha_1} 
   \left[\begin{array}{c} X_{21} \\*[.5ex] X_{31} \\*[.5ex] 
   \vdots \\*[.5ex] X_{r1}
   \end{array}\right]
   \left[\begin{array}{c} X_{21} \\*[.5ex] X_{31} \\*[.5ex] 
   \vdots \\*[.5ex] X_{r1}
   \end{array}\right]^\top.
\]
By property~\ref{p-ishi2}, the diagonal blocks 
$Y_{ii} = \alpha_i I - (1/\alpha_1) X_{i1}X_{i1}^\top$
are multiples of the identity, so $Y_{ii} \in \mathcal V_{ii}$.
Property~\ref{p-ishi3} implies that  for $i>j$,
\[
Y_{ij} = X_{ij} - \frac{1}{\alpha_1} X_{i1}X_{j1}^\top
 \in \mathcal V_{ij}.
\]
Hence by the induction hypothesis, the matrix $Y$ can be factored as
\[
 \left[\begin{array}{ccccc}
 Y_{22}  & Y_{32}^\top  & \cdots & Y_{r2}^\top \\*[.5ex]
 Y_{32}  & Y_{33}  & \cdots & Y_{r3}^\top  \\*[.5ex]
 \vdots  & \vdots & \ddots & \vdots \\*[.5ex]
 Y_{r2} & Y_{r3}   & \cdots &  Y_{rr} \end{array} \right] 
= 
 \left[\begin{array}{ccccc}
 L_{22}  & 0  & \cdots & 0 \\*[.5ex]
 L_{32}  & L_{33}  & \cdots & 0  \\*[.5ex]
 \vdots  & \vdots & \ddots & \vdots \\*[.5ex]
 L_{r2} & L_{r3}   & \cdots &  L_{rr} \end{array} \right] 
 \left[\begin{array}{ccccc}
 L_{22}^\top  & L_{32}^\top  & \cdots & L_{r2}^\top \\*[.5ex]
 0  & L_{33}^\top  & \cdots & L_{r3}^\top  \\*[.5ex]
 \vdots  & \vdots & \ddots & \vdots \\*[.5ex]
 0 & 0   & \cdots &  L_{rr}^\top \end{array} \right] 
\]
with $L_{ij} \in \mathcal V_{ij}$.
Substituting the factorization of $Y$ in~(\ref{e-X-chol-pf}) gives
a Cholesky factorization $X=LL^\top$ with the desired properties.
\end{proof}

\subsection{Computations with triangular matrices} 
The following result generalizes Theorem~\ref{thm:1.1} to 
triangular matrices in a subspace with the structure
specified in Properties~\ref{p-ishi1}--\ref{p-ishi4}.
\begin{theorem}
\label{thm:V}\ 
Let $\mathcal T$ be the set of triangular matrices~(\ref{e-T-def}),
where the subspaces $\mathcal V_{ij}$ satisfy 
Properties \ref{p-ishi1}--\ref{p-ishi4}, and assume $L \in\mathcal T$.
\begin{enumerate}[label = \arabic{enumi}.]
\item If $\tilde L \in \mathcal T$, then $L\tilde L \in \mathcal T$.
\item If $L$ is nonsingular, then $L^{-1} \in \mathcal T$.
\item If $X\in\mathcal V$, then $LXL^\top \in \mathcal V$.
\item If $Y\in\mathcal V^\perp$, then $L^\top YL \in \mathcal V^\perp$. 
\end{enumerate}
\end{theorem}

\begin{proof}
Suppose $L, \tilde L \in \mathcal T$ are partitioned as
\[
L = \left[\begin{array}{ccccc}
   \alpha_1I & 0 & 0 & \cdots & 0 \\*[.5ex]
   L_{21} & \alpha_2I & 0 & \cdots & 0 \\*[.5ex]
   L_{31} & L_{32} & \alpha_3 I & \cdots & 0 \\*[.5ex]
   \vdots & \vdots & \vdots & \ddots & \vdots \\*[.5ex]
   L_{r1} & L_{r2} & L_{r3} & \cdots &  \alpha_r I \end{array} \right],
\qquad
\tilde L = \left[\begin{array}{ccccc}
   \tilde \alpha_1I & 0 & 0 & \cdots & 0 \\*[.5ex]
   \tilde L_{21} & \tilde \alpha_2I & 0 & \cdots & 0 \\*[.5ex]
   \tilde L_{31} & \tilde L_{32} & \tilde \alpha_3 I & \cdots & 0 \\*[.5ex]
   \vdots & \vdots & \vdots & \ddots & \vdots \\*[.5ex]
   \tilde L_{r1} & \tilde L_{r2} & \tilde L_{r3} & \cdots &  
   \tilde \alpha_r I \end{array} \right].
\]
The diagonal blocks in the product $L\tilde L$ are $(L\tilde L)_{ii}
= \alpha_i\tilde\alpha_i I \in \mathcal V_{ii}$.
For the lower-triangular off-diagonal blocks,
\[
(L\tilde L)_{ij} 
 = \tilde \alpha_j L_{ij} + \sum_{k=j+1}^{i-1}L_{ik}\tilde L_{kj}
   + \alpha_i \tilde L_{ij}.
\]
By assumption, the terms 
$L_{ij}$, $\tilde L_{ij}$ are in $\mathcal V_{ij}$.
The middle term in the expression on the right-hand side 
is in $\mathcal V_{ij}$ by property~\ref{p-ishi4}.
Hence, \ref{p-ishi1} and~\ref{p-ishi4} are sufficient to prove 
statement~1 of the theorem.

Part 2 is proved by induction on $r$.  For $r=1$ it is obvious, 
with $L=\alpha I$ and $L^{-1}=\alpha^{-1}I$.
Suppose the result holds for $r=m-1$ and consider a matrix 
$L\in\mathcal T$, partitioned in $r$-by-$r$ blocks as above.
By the induction hypothesis, the blocks in 
\[
 \left[\begin{array}{cccc}
  L_{22} & 0 & \cdots & 0 \\
  L_{32} & L_{33} & \cdots & 0 \\
  \vdots & \vdots & \ddots & \vdots \\
  L_{r2} & L_{r3} & \cdots & L_{rr} \end{array}\right]^{-1}
= 
 \left[\begin{array}{cccc}
  \alpha_2^{-1} I & 0 & \cdots & 0 \\
  G_{32} & \alpha_3^{-1}I & \cdots & 1 \\
  \vdots & \vdots & \ddots & \vdots \\
  G_{r2} & G_{r3} & \cdots & \alpha_r^{-1} I \end{array}\right]
\]
satisfy $G_{ij} \in \mathcal V_{ij}$ for $i>j>1$.
The inverse of $L$ is
\[
L^{-1} 
= \left[\begin{array}{ccccc}
  \alpha_1^{-1}I & 0 & 0 & \cdots & 0 \\
  G_{21} & \alpha_2^{-1} I & 0 & \cdots & 0 \\
  G_{31} & G_{32} & \alpha_3^{-1} I & \cdots & 0  \\
  \vdots        & \vdots & \vdots & \ddots & \vdots \\
  G_{r1} & G_{r2} & G_{r3} & \cdots & \alpha_r^{-1} I 
\end{array}\right]
\]
with
\begin{eqnarray*}
G_{21} & = & -(\alpha_1\alpha_2)^{-1} L_{21} \\
G_{31} & = & -\alpha_1^{-1} G_{32}L_{21} -(\alpha_1\alpha_3)^{-1} L_{31} \\
 & \vdots & \\
G_{r1} & = & -\alpha_1^{-1} (G_{r2} L_{21} + \cdots + G_{r,r-1}L_{r-1,1})
   - (\alpha_1\alpha_r)^{-1}L_{r1}.
\end{eqnarray*}
From the induction hypothesis ($G_{ij} \in\mathcal V_{ij}$ for $i>j>1$)
and Property~\ref{p-ishi4}, we see that $G_{i1}\in\mathcal V_{i1}$
for $i \in \{2,\ldots, r\}$.
Hence, statement~2 of the theorem follows from \ref{p-ishi1} 
and~\ref{p-ishi4}.

Next we show part 3.  
First consider the diagonal blocks of $LXL^\top$, 
\[
(LXL^\top)_{ii}
= \sum_{k=1}^i L_{ik}X_{kk}L_{ik}^\top + \sum_{k=2}^i \sum_{l=1}^{k-1} 
\left(L_{ik} X_{kl} L_{il}^\top  + 
L_{il} X_{kl}^\top L_{ik}^\top\right).
\]
Properties~\ref{p-ishi1}, \ref{p-ishi2}, and~\ref{p-ishi3} imply that
this is a product of the identity.  For the off-diagonal blocks
with $i>j$,
\[
(LXL^\top)_{ij}
= 
\sum_{k=1}^j L_{ik}X_{kk}L_{jk}^\top + 
\sum_{k=2}^j \sum_{l=1}^{k-1}  \left( L_{ik} X_{kl} L_{jl}^\top  + 
  L_{il} X_{kl}^\top L_{jk}^\top\right).
\] 
Properties \ref{p-ishi1}, \ref{p-ishi3}, \ref{p-ishi4} imply
that this is an element in $\mathcal V_{ij}$.

Part 4 is an immediate consequence of part 3.
Suppose $Y\in\mathcal V^\perp$.  For any $X\in\mathcal V$,
\[
\Tr(X L^\top YL) = \Tr(LXL^\top Y) =  0
\]
because $LXL^\top \in\mathcal V$ by part 3.
Therefore $L^\top YL \in\mathcal V^\perp$.
\end{proof}

\subsection{Primal and dual cone automorphisms} \label{s-homogeneous-aut}
Now consider the cones~(\ref{e-K-V}) and~(\ref{e-K*-V}), where
$\mathcal V$ is a  subspace that satisfies the four 
properties \ref{p-ishi1}--\ref{p-ishi4}.

The linear mappings $\mathcal L(X) = LXL^\top$ for nonsingular
$L\in \mathcal T$ form a transitive subset of $\aut(K)$. 
This is readily shown by extending the arguments in 
Section~\ref{s-primal-auto} using Theorem~\ref{thm:V} and 
the property that Cholesky factors of matrices in $\inte(K)$ are in 
$\mathcal T$ (see Section~\ref{s-homogeneous-chol}).
In the remainder of Section~\ref{s-hom-matrix-cones}, 
$\auttr(K)$ will be used to denote this
transitive subset of $\aut(K)$.  

The adjoints $\mathcal L^*(S) = \Pi_\mathcal V(L^\top SL)$
of mappings $\mathcal L \in \auttr(K)$ 
form a transitive subset of $\aut(K^*)$.   
The proof again parallels the arguments in Section~\ref{s-dual-auto}.
As noted in Section~\ref{s-dual-auto}, every adjoint of an 
automorphism of $K$ is an automorphism of $K^*$.
We also note that every $S\in \inte(K^*)$ can be expressed as
\[
S = \Pi_\mathcal V(X^{-1})
\]
where $X\in \inte(K)$.  The matrix $X$ is the solution of the
convex optimization problem 
\[
 \begin{array}{ll}
 \mbox{minimize} & \Tr(SX) -\ln\det(X) \\
 \end{array}
\]
with variable $X\in \mathcal V$.
The Cholesky factorization $X = LL^\top$, where $L\in \mathcal T$,
defines a mapping $\mathcal L\in \auttr(K)$ that satisfies
$\mathcal L^{-*}(I) =  \Pi_\mathcal V(L^{-\top} L^{-1}) = S$.
Therefore $\mathcal L^*(S) = I$.

Now consider any two matrices $S_1, S_2\in \inte(K^*)$.
Find $\mathcal L_1, \mathcal L_2 \in \auttr(K)$ that
satisfy $\mathcal L_1^*(S_1) = I$, $\mathcal L_2^*(S_2) = I$.
Then $\mathcal L  = \mathcal L_1 \circ \mathcal L_2^{-1}$
satisfies $\mathcal L^*(S_1) = S_2$ and 
$\mathcal L \in \auttr(K)$ by properties 1--3 in 
Theorem~\ref{thm:V}.
By establishing transitive subsets of $\aut(K)$ and $\aut(K^*)$ we have
shown that these cones are homogeneous.

Using our results, in particular Theorem~\ref{thm:V},
we can establish the following fact.
\begin{theorem}
\label{thm:6-squared}
Let $K$ be a homogeneous cone represented in $\SS^N$ as described in 
\eqref{e-K-V}. Then, for every $Z \in \cV \cap
\SS^N$, 
upon expressing $Z = L +L^{\top}$ for some $L \in \mathcal T$, we have
\[
L L^{\top} +\Pi_{\cV}\left(L^{\top} L\right) + L^2 +\left(L^{\top}\right)^2 \in K^*.
\]
\end{theorem}

Decompositions as the above have potential applications in linear 
and nonlinear complementarity problems over homogeneous cones and in the 
design of algorithms and theories utilizing Moreau decompositions, see for instance,
\cite{KongXiuTuncel2012}.

\subsection{Logarithmic barriers}
If $K=\mathcal V \cap \SS^N_+$ and $\mathcal V$ satisfies 
properties~\ref{p-ishi1}--\ref{p-ishi4},
then the log-det barrier 
\begin{equation} \label{e-log-det-barrier-hom}
F(X) := \left\{ \begin{array}{rl} -\ln\det(X), & \mbox{ if }
 X \in \inte(K) \\
+\infty, & \mbox{ otherwise}
\end{array}
\right.
\end{equation}
has the same scaling properties~(\ref{e-composition})
as the log-det barrier for a homogeneous sparse matrix cone.
The proof is exactly the same.  From~(\ref{e-F-scaling})  
and the fact that $\mathcal L$ is an automorphism, it follows
that~(\ref{e-composition}) holds for all $X\in\inte(K)$ and
all $\mathcal L \in \auttr(K)$.
Similarly, (\ref{e-composition-dual}) holds for all
$S\in\inte(K^*)$ and all $\mathcal L \in \auttr(K)$.

With the above definition of the barrier function $F$, 
Theorem~\ref{thm:asymmetric-factor} extends to all homogeneous matrix
cones $K$ discussed in this section.  This is stated in the following
theorem.
\begin{theorem}
\label{thm:6.3}
Let $K = \cV \cap\SS^N_+$, where $\cV$ is a subspace that satisfies
properties \ref{p-ishi1}--\ref{p-ishi4}.
Then, for every $X \in \inte(K)$ and $S \in \inte(K^*)$, there 
exists $\mathcal{L} \in \auttr(K)$ such that
\[
 \mathcal{L}^{-1}(X) = \mathcal{L}^* (S).
\]
\end{theorem}

Ishi has shown that every homogeneous cone can be represented in 
the form $K = \cV \cap\SS^N_+$, where $\cV$ is satisfies
properties \ref{p-ishi1}--\ref{p-ishi4}.
Theorem~\ref{thm:6.3} therefore shows that problem (2) on page 711 of 
\cite{Tun:98} is solvable for every homogeneous cone, and settles the 
open problem (i) (on page 714) of \cite{Tun:98}.

Next, we relate the above homogeneous matrix cones representation to 
algebraic classifications of all homogenous cones and explain why the 
results of this section apply to all homogeneous cones.

\section{Algebraic structure of homogeneous cones}
\label{sec:homogeneous-cones}

In the previous sections we discussed classes of homogeneous cones 
defined as linear slices of the positive semidefinite cone.
It turns out that every homogeneous cone can be expressed in this form.
As mentioned by Faybusovich \cite[p.214]{Faybusovich2002}
and Papp and Alizadeh \cite[p.1406]{PaA:13}, and worked out in detail
by Chua \cite{Chua2003}, this result is implicit in Vinberg's 
$T$-algebra based classification of homogeneous cones,
because Vinberg's results imply that every homogeneous cone is a 
``cone of squares'' for a suitable vector product.
Rothaus, announcing a similar result first in 1963, proved it using the 
inductive Siegel Domain based classification of homogeneous cones and 
convex cone duality \cite{Rothaus1963,Rothaus1966,Rothaus1968}.
Ishi's approach \cite{Ishi2013,Ishi2016,Ishi2015}, 
influenced in part by some recent work by Yamasaki and 
Nomura~\cite{YamasakiNomura2015},
brings Rothaus's Siegel Domain based inductive construction closer to 
more direct utilization of the $T$-algebra axioms. In this section, 
we discuss some of the results by Vinberg and Rothaus, and explain their
connections with the classes of homogeneous matrix cones described in 
Sections~\ref{s-cone}--\ref{s-hom-matrix-cones}.

It is useful to first clarify the meaning
of \emph{semidefinite representation} of a convex cone.
A convex cone $\cV \cap \SS^N_+$,  
where $\cV \subseteq \SS^N$ is a linear subspace, can be equivalently
represented as
\begin{equation} \label{e-lmi}
K=\left\{x \in \R^n : \sum_{i=1}^n x_i A_i \succeq 0 \right\},
\end{equation}
where $A_1, A_2, \ldots, A_n \in \SS^N$.
Given the subspace in the representation $\cV \cap \SS^N_+$,
we can pick a basis $A_1, A_2, \ldots, A_n \in \SS^N$ for $\cV$ to obtain 
the representation~(\ref{e-lmi}).
Given $A_1, A_2, \ldots, A_n \in \SS^N$ in the second 
representation, we define $\cV := \spam\{A_1, A_2, \ldots, A_n\}$
to obtain the former representation. 
The representation~(\ref{e-lmi}) is called
a \emph{linear matrix inequality (LMI)} or
\emph{spectrahedral} representation of the cone $K$.
In spectrahedral representations one typically
requires that $\cV \cap \SS^N_{++} \neq \emptyset$ 
(i.e., that there exists $\bar{x} \in \R^n$ such that
$\sum_{i=1}^m \bar{x}_i A_i \succ 0$).

Whenever a regular cone 
admits a spectrahedral representation with 
$\inte(K) = \cV \cap \SS^N_{++} \neq \emptyset$, the dual cone
in the space $\SS^N$, under the trace inner product, is given by
$\cV^{\perp} + \SS^N_+$.  (In general, a closure operation 
is needed on the right-hand side. However, 
it can be shown that the cone $\cV^\perp + \SS^N_+$
is closed if $\cV \cap \SS^N_{++} \neq \emptyset$, 
so the closure operation can be omitted.)
If we take the dual of $K=\cV \cap \SS^N_+$ with respect to the smaller 
space $\cV$, the dual cone is
 \[
 K^* = \Pi_{\cV}\left(\cV^{\perp}+\SS^N_+\right) = \Pi_{\cV}(\SS^N_+).
 \]
The dual cone can therefore be represented in the form 
\begin{equation} \label{e-spec-shadow}
 K^* = \left\{s \in \R^n : 
 \sum_{i=1}^n s_i H_i + \sum_{j=1}^k u_j U_j \succeq 0, 
 \mbox{ for some } u \in \R^k \right\},
\end{equation}
where $H_1, H_2, \ldots H_n, U_1, U_2, \ldots, U_k \in \SS^N$ are given.
This kind of semidefinite representation is called 
a \emph{lifted-LMI} or  \emph{spectrahedral shadow} representation 
(of $K^*$);
see \cite{HeltonVinnikov2007,Nem2007,ChuaTuncel2008,HeltonNie2010,GPT2013,Sch2018,Averkov2019,Fawzi2020} and the references therein.  
In a spectrahedral shadow representation the dual cone is expressed as the 
cone of positive semidefinite ``completable'' matrices (``completable'' 
by some element of $\cV^{\perp}$).  In our context, for the spectrahedral 
shadow representation~(\ref{e-spec-shadow}), $\{H_1, H_2, \ldots, H_n\}$
is a basis for $\cV$ and $\{U_1, U_2, \ldots, U_k\}$ can be taken as a 
basis for $\cV^{\perp}$.  If so, then $n+k = N(N+1)/2$.
 
Note that by our choices 
for these representations of $K$ and $K^*$ (i.e., for this choice of
inner product and the space), we always have $K \subseteq K^*$.

\subsection{Symmetric bilinear forms}
\begin{dfn} \label{def:K-bilinear-form}
Let $K$ be a homogeneous cone in a finite-dimensional real 
vector space $\cV$.
A \emph{homogeneous $K$-bilinear symmetric form}
$\cB(u,v)$ is a mapping from $\R^p\times\R^p$ to
$\cV$ that satisfies the following properties.
\begin{enumerate}[label = \arabic{enumi}.]
\item $\cB(\alpha_1 u^{(1)} +\alpha_2 u^{(2)},v) = 
\alpha_1 \cB(u^{(1)},v) + \alpha_2 \cB(u^{(2)},v)$ 
for all $u^{(1)}, u^{(2)}, v \in \R^p$ and $\alpha_1,\alpha_2\in{\R}$.
\item $\cB(u,v) = \cB(v,u)$ for all $u,v \in \R^p$.
\item $\cB(u,u)\in K$ for all $u \in \R^p$.
\item $\cB(u,u)=0$ implies $u=0$. 
\item
There exists a transitive subset
$G\subseteq \aut(K)$ such that for every $g\in G$, there
exists a linear transformation $\bar{g}$ on ${\R}^p$
which satisfies 
\begin{equation} \label{e-bilinear-prop5}
g(\cB(u,v)) = \cB(\bar{g}u,\bar{g}v) \quad
\mbox{for all  $u,v \in \R^p$}.
\end{equation}
\end{enumerate}
\end{dfn}
\noindent In this definition, $p=0$ is allowed.
When $p=0$, the mapping $\cB$ is the trivial bilinear form 
(a constant zero vector).

We now discuss some implications of the five properties in the 
definition.  We use the standard inner product $u^\top v$ in $\R^p$,
an inner product $\langle s, x\rangle$ in $\cV$, and denote by 
\[
K^* = \{ s \in \cV : \langle s, x \rangle \geq 0 \mbox{\ for all
$x\in K$}\} 
\]
the corresponding dual cone.
The trace inner product is used for symmetric matrices.

A function $\cB$ that satisfies properties~1 and~2 in 
Definition~\ref{def:K-bilinear-form} is called a \emph{symmetric bilinear 
form}.  
With every symmetric bilinear form $\cB$ one can associate a linear matrix 
function $\cH: \cV \rightarrow \SS^p$, defined by the identity
\begin{equation} \label{e-cB-cH}
u^\top \cH(s) v = v^\top \cH(s)u = \langle s, \cB(u,v)\rangle  
\quad \mbox{for all $u,v\in\R^p$, $s\in \cV$}.
\end{equation}
An explicit formula for the entries of $\cH(s)$ is 
\begin{equation} \label{e-cH-def}
 \cH(s)_{ij} = \langle s, \cB(e_i,e_j)\rangle, \quad i,j=1,\ldots,p,
\end{equation}
where $e_1 = (1,0,\ldots,0)$, $e_2 = (0,1,0,\ldots,0)$,
\ldots, $e_p = (0,\ldots,0,1)$ are the standard unit vectors in $\R^p$.
This expression follows from~(\ref{e-cB-cH}) if we use the bilinearity 
property~1 in the definition to expand $\cB(u,v)$ as
\begin{eqnarray} 
\cB(u,v) 
& = &
\cB(u_1 e_1 + \cdots + u_p e_p, v_1 e_1 + \cdots + v_p e_p) \nonumber \\
 & = & \sum_{i=1}^p \sum_{j=1}^p u_i v_j \cB(e_i,e_j).
\label{e-cB-explicit}
\end{eqnarray}
We will refer to $\cH$ as the \emph{dual representation}
of the bilinear form $\cB$.
The adjoint of $\cH$ (with respect to the inner product $\langle \cdot,
\cdot \rangle$ in $\cV$ and the trace inner product in $\SS^p$) 
is the linear mapping $\cH^*: \R^p \rightarrow \cV$ that maps a matrix 
$Y\in\SS^p$ to the vector 
\[
 \cH^*(Y) =  \sum_{i=1}^p \sum_{j=1}^p Y_{ij} \cB(e_i,e_j).
\]
Hence, from~(\ref{e-cB-explicit}), 
we have the following expression for $\cB$:
\begin{equation} \label{e-cB-cH*}
\cB(u,v) = \cH^*((uv^\top + vu^\top)/2).
\end{equation}
In particular, $\cB(u,u) = \cH^*(uu^\top)$.
If $\cB$ is the trivial bilinear form, we define $\cH^*$ as the constant
zero in $\cV$.

The formula $\cB(u,u) = \cH^*(uu^\top)$ has an important consequence 
for semidefinite programming applications. It implies that 
the ``sum of squares'' cone 
\begin{equation} \label{e-sos-cone}
C := \left\{ \sum_{i=1}^k \cB(u^{(i)}, u^{(i)}) :
  \mbox{for some $k$ and $u^{(1)}, \ldots, u^{(k)} \in \R^p$} \right\}
\end{equation}
of any symmetric bilinear form $\cB$ has a spectrahedral representation
$C = \{ \cH^*(Y) : Y \succeq 0\}$
\cite{Nes:00,Faybusovich2002,PaA:13}.
This follows from $\cB(u,u) = \cH^*(uu^\top)$ and linearity of $\cH^*$:
all elements in $C$ can be expressed as
\[
\sum_{i=1}^k \cB(u^{(i)}, u^{(i)}) 
= \sum_{i=1}^k \cH^*(u^{(i)} (u^{(i)})^\top) 
= \cH^*(Y)
\]
where $Y  = \sum_i u^{(i)} (u^{(i)})^\top \succeq 0$,
and, conversely, if 
$x = \cH^*(Y)$ with $Y\succeq 0$, then any decomposition
$Y = \sum_i u^{(i)} (u^{(i)})^\top$
gives an expression $x=\sum_i \cB(u^{(i)}, u^{(i)})$ that shows that
$x\in C$.

Properties 3--5 in Definition~\ref{def:K-bilinear-form} can 
be stated in equivalent forms involving the dual representation $\cH$
and its adjoint.
\begin{proposition}
\label{thm:cH}\ 
Let $K$ be a regular convex cone in a finite-dimensional real vector space
$\cV$.  Let $\cB:\R^p \times \R^p \rightarrow \cV$ be a symmetric 
bilinear form and $\cH$ its dual representation defined in~(\ref{e-cB-cH}).
\begin{enumerate}[label = \arabic{enumi}.]
\item 
Each of the following two statements is equivalent to 
the property that $\cB(u,u) \in K$ for all~$u$:
\begin{subequations} \label{e-cH-prop1}
\begin{equation}
\cH^*(Y) \in K \; \mbox{\ for all $Y\succeq 0$}, 
\end{equation}
\begin{equation}
\cH(s) \succeq 0 \; \mbox{\ for all $s\in K^*$}. 
\end{equation}
\end{subequations}

\item 
Each of the following two statements is equivalent to 
the property that $\cB(u,u) \in K\setminus\{0\}$ for all $u\neq 0$:
\begin{subequations} \label{e-cH-prop2}
\begin{equation}\label{e-cH-prop2-a}
\cH^*(Y) \in K\setminus\{0\} \; \mbox{ for all nonzero $Y\succeq 0$}, 
\end{equation}
\begin{equation}
\cH(s) \succ 0 \; \mbox{\ for all $s\in \inte(K^*)$}.
\end{equation}
\end{subequations}

\item Let $g:\cV\rightarrow\cV$ and $\bar g:\R^p\rightarrow\R^p$ be 
linear transformations.
Each of the following two statements is equivalent to the
property that $g(\cB(u,v)) = \cB(\bar g u, \bar g v)$ for all 
$u,v\in\R^p$:
\begin{subequations}\label{e-cH-prop3}
\begin{equation} 
g(\cH^*(Y)) = \cH^*(\bar g Y \bar g^\top) \; \mbox{\ for all $Y\in \SS^p$},
\end{equation}
\begin{equation} 
\cH(g^*(s)) = \bar g^\top \cH(s) \bar g \; \mbox{\ for all  $s\in \cV$}. 
\end{equation}
\end{subequations}
\end{enumerate}
\end{proposition}

Part 3 of the proposition follows directly 
from~(\ref{e-cB-cH})  and~(\ref{e-cB-cH*}).
The statements about $\cH^*$ in the first two parts follow
from $\cB(u,u) = \cH^*(uu^\top)$ and linearity of $\cH^*$.
The statements about $\cH$ follow from
$u^\top \cH(s) u = \langle s, \cB(u,u)\rangle$, and the equivalences
\begin{eqnarray*}
\cB(u,u) \in K
& \Longleftrightarrow &  
\langle s, \cB(u,u) \rangle \geq 0 \quad \mbox{for all $s\in K^*$} \\
& \Longleftrightarrow &  
u^\top \cH(s) u \geq 0 \quad \mbox{for all $s\in K^*$}
\end{eqnarray*}
and
\begin{eqnarray*}
\cB(u,u) \in K\setminus\{0\}
& \Longleftrightarrow &  
\langle s, \cB(u,u) \rangle > 0 \quad \mbox{for all $s\in \inte(K^*)$} \\
& \Longleftrightarrow &  
u^\top \cH(s) u > 0 \quad \mbox{for all $s\in \inte(K^*)$}.
\end{eqnarray*}

\begin{examp}\label{examp:6.1}
We take $\cV = \SS^3_E$ where $E$ is the homogeneous chordal pattern in
Figure~\ref{f-vinberg}, i.e., $\cV$ is the space of matrices of the form
\[
 X = \left[\begin{array}{ccc} 
   X_{11} & 0 & X_{31} \\
   0 & X_{22} & X_{32}  \\
   X_{31} & X_{32} & X_{33} \end{array}\right].
\]
We use the inner product $\langle S, X\rangle = \Tr(SX)$ on $\cV$
and define $K = \cV \cap \SS^3_+$. 
(This cone is known as the \emph{Vinberg cone}, the smallest dimensional
homogeneous cone which is not symmetric.)
Consider the following symmetric bilinear form 
$\cB : \R^6\times \R^6 \rightarrow \cV$:
\begin{eqnarray*}
\lefteqn{\cB(u,v)} \\
& = & 
\frac{1}{2} 
 \left[\begin{array}{ccc}
 u_1 & u_3 & 0 \\ 0 & 0 & u_5 \\ u_2 & u_4 & u_6 \end{array}\right]
 \left[\begin{array}{ccc}
 v_1 & v_3 & 0 \\ 0 & 0 & v_5 \\ v_2 & v_4 & v_6 \end{array}\right]^\top
+
 \left[\begin{array}{ccc}
 v_1 & v_3 & 0 \\ 0 & 0 & v_5 \\ v_2 & v_4 & v_6 \end{array}\right]
 \left[\begin{array}{ccc}
 u_1 & u_3 & 0 \\ 0 & 0 & u_5 \\ u_2 & u_4 & u_6 \end{array}\right]^\top
\\
& = & 
\frac{1}{2}
 \left[\begin{array}{ccc}
  2(u_1v_1 + u_3v_3) & 0 & u_1v_2+u_2v_1+ u_3v_4 + u_4v_3 \\
 0 & 2u_5v_5 & u_5v_6+u_6v_5 \\
 u_1v_2+u_2v_1+ u_3v_4 + u_4v_3 & 
 u_5v_6+u_6v_5 & 2(u_2 v_2 + u_4v_4 + u_6v_6)
 \end{array}\right].
\end{eqnarray*}
The dual representation $\cH:\cV \rightarrow \SS^6$ is
\[
 \cH(S) = \left[\begin{array}{cccccc}
 S_{11}  & S_{31} & 0      & 0      & 0  & 0 \\
 S_{31}  & S_{33} & 0      & 0      & 0  & 0 \\
 0       & 0      & S_{11} & S_{31} & 0  & 0 \\   
 0       & 0      & S_{31} & S_{33} & 0  & 0 \\
 0       & 0      & 0      & 0      & S_{22} & S_{32} \\
 0       & 0      & 0      & 0      & S_{32} & S_{33}
\end{array}\right], 
\]
and its adjoint
$\cH^*: \SS^6 \rightarrow \cV$ is
\[
 \cH^*(Y) = \left[\begin{array}{ccc}
 Y_{11}+Y_{33}  & 0       & Y_{21} + Y_{43} \\
 0              & Y_{55}  & Y_{65}  \\
 Y_{21} + Y_{43} & Y_{65} & Y_{22} + Y_{44} + Y_{66} \end{array}\right]. 
\]
This bilinear form $\cB$ satisfies the five properties
in Definition~\ref{def:K-bilinear-form}. 
It satisfies properties~3 and~4, as can be seen from 
\begin{eqnarray*}
\cB(u,u) & = &
\frac{1}{2} \left[\begin{array}{ccc}
 u_1 & u_3 & 0 \\ 0 & 0 & u_5 \\ u_2 & u_4 & u_6 \end{array}\right]
 \left[\begin{array}{ccc}
 u_1 & u_3 & 0 \\ 0 & 0 & u_5 \\ u_2 & u_4 & u_6 \end{array}\right]^\top
\\
& = & \left[\begin{array}{ccc}
  u_1^2 + u_3^2 & 0 & u_1u_2 + u_3u_4 \\
  0 & u_5^2 & u_5u_6 \\
  u_1u_2 + u_3 u_4 & u_5u_6 & u_2^2 + u_4^2 + u_6^2 \end{array}\right].
\end{eqnarray*}
The first expression shows that $\cB(u,u)$ is positive semidefinite
for all $u$; the second expression that $\cB(u,u)=0$ only if $u=0$.
For property~5 we use the transitive subset of $\aut(K)$ discussed 
in Section~\ref{s-primal-auto}. The automorphisms in $G$ are
the mappings $g = \cL:\cV \rightarrow \cV$ defined 
as $\cL(X) = LXL^T$, where $L$ is a nonsingular triangular matrix
\[
 L = \left[\begin{array}{ccc}
  L_{11} & 0 & 0 \\ 0 & L_{22} & 0 \\ L_{31} & L_{32} & L_{33}
 \end{array}\right].
\]
Then $\cL(\cB(u,v)) = \cB(\bar g u, \bar g v)$ where 
\[
\bar g = 
\left[\begin{array}{cccccc}
L_{11} & 0      & 0      & 0      & 0      & 0 \\
L_{31} & L_{33} & 0      & 0      & 0      & 0 \\ 
0      & 0      & L_{11} & 0      & 0      & 0 \\
0      & 0      & L_{31} & L_{33} & 0      & 0 \\
0      & 0      & 0      & 0      & L_{22} & 0 \\
0      & 0      & 0      & 0      & L_{32} & L_{33} 
\end{array}\right].
\]
\end{examp}

\begin{examp}\label{examp:6.2}
With the same choice of $\cV$ and $K$,
define $\cB : \R^2 \times \R^2 \rightarrow \cV$ as
\begin{eqnarray*}
\cB(u,v) & = & \frac{1}{2} 
\left[\begin{array}{c} u_1 \\ 0 \\ u_2  \end{array}\right]
\left[\begin{array}{c} v_1 \\ 0 \\ v_2  \end{array}\right]^\top
+ \frac{1}{2}
\left[\begin{array}{c} v_1 \\ 0 \\ v_2  \end{array}\right]
\left[\begin{array}{c} u_1 \\ 0 \\ u_2  \end{array}\right]^\top \\
& = & \frac{1}{2} \left[\begin{array}{ccc}
 2 u_1 v_1 & 0 &  u_1v_2 + u_2 v_1 \\ 
 0       & 0 & 0 \\
 u_1v_2 + u_2 v_1 & 0 & 2u_2 v_2
 \end{array}\right].
\end{eqnarray*}
The dual representation $\cH: \cV \rightarrow \SS^2$ and its adjoint
$\cH^*: \SS^2 \rightarrow \cV$ are
\[
\cH(S) = \left[\begin{array}{cc} 
 S_{11} & S_{31} \\ S_{31} & S_{33} 
\end{array}\right],
\qquad
\cH^*(Y) = \left[\begin{array}{ccc}
 Y_{11} & 0 & Y_{21} \\
 0      & 0 & 0      \\
 Y_{21} & 0 & Y_{22}
\end{array}\right].
\]
Here,
\[
\cB(u,u) = 
\left[\begin{array}{c} u_1 \\ 0 \\ u_2  \end{array}\right]
\left[\begin{array}{c} u_1 \\ 0 \\ u_2  \end{array}\right]^\top
= \left[\begin{array}{ccc}
 u_1^2  & 0 & u_1u_2 \\ 0       & 0 & 0 \\
u_1u_2 & 0 & u_2^2 \end{array}\right],
\]
which satisfies properties~3 and~4 in 
Definition~\ref{def:K-bilinear-form}.
Property~5 holds for the same transitive subset $G$ as
in the previous example and
\[
\bar g = \left[\begin{array}{cc}
 L_{11} & 0 \\ L_{31} & L_{33} \end{array}\right].
\]
Therefore $\cB$ is another homogeneous $K$-bilinear symmetric form
for the same cone $K$.
Note that, in contrast to the previous example,
$K$ is not equal to the sum-of-squares cone~(\ref{e-sos-cone}). 
Here, strict inclusions  
$K \supset \{ \cH^*(Y) : Y \succeq 0 \}$ and
$K^* \subset \{ S : \cH(S) \succeq 0 \}$ hold.
\end{examp}

\subsection{Siegel cone}
Let $K$ be a homogeneous cone and $\cB$ a homogeneous $K$-bilinear
symmetric form as defined in Definition~\ref{def:K-bilinear-form}.
We define the \emph{Siegel cone} associated with $K$ and $\cB$ as
\begin{eqnarray}
\lefteqn{\SC(K,\cB) }  \nonumber \\
& := & \left\{(x,u,\alpha) \in \cV \times \R^p \times \R :
  x\in K, \, \alpha \geq 0, \, \alpha x - \cB(u,u) \in K\right\}
 \nonumber \\
& = & 
  \left\{(x,u,\alpha) : \alpha > 0, \, 
 x - \frac{1}{\alpha} \cB(u,u) \in K \right\} \;
 \cup \; \left\{(x,0,0) : x\in K \right\}.  \label{e-SC-def} 
\end{eqnarray}
If $\cB$ is the trivial bilinear form ($p=0$), 
the Siegel cone is $\SC(K,\cB) = K \times \R_+$.

The following equivalent definition follows from the results in the 
previous section and makes it clear that $\SC(K,\cB)$ is convex: 
if $p\geq 1$,
\begin{eqnarray} 
\lefteqn{\SC(K,\cB)} \nonumber \\
 & = & \left\{(x,u,\alpha) : 
  x - \cH^*(Y) \in K, \left[\begin{array}{cc}
    \alpha & u^\top \\ u  & Y \end{array}\right] \succeq 0
 \mbox{\ for some $Y \in \SS^p$}\right\}.
\label{e-SC-alt}
\end{eqnarray}
This definition also shows that $\SC(K,\cB)$ has a spectrahedral shadow
representation if the cone $K$ has a spectrahedral shadow representation.
The equivalence of~(\ref{e-SC-def}) and~(\ref{e-SC-alt}) can be seen 
as follows. 
We first note that in both definitions the only elements 
$(x,u,\alpha)$ with $\alpha=0$ are the vectors $(x,0,0)$, $x\in K$.  
If $\alpha=0$, the matrix inequality in~(\ref{e-SC-alt}) requires
$u=0$ and $Y\succeq 0$.
Since $\cH^*(Y) \in K$ for all $Y\succeq 0$, the condition on $x$ then
reduces to $x\in K$.  
Next, suppose $\alpha > 0$ and $(x,u,\alpha)$ is in the 
cone~(\ref{e-SC-def}).  
Then $Y=(1/\alpha) uu^\top$ satisfies the conditions in~(\ref{e-SC-alt}), 
so $(x,u,\alpha)$ is in the cone~(\ref{e-SC-alt}).
Conversely, suppose $\alpha > 0$ and $(x,u,\alpha)$ satisfies the
conditions in~(\ref{e-SC-alt}) for some $Y\succeq 0$.
Then $Y\succeq (1/\alpha)uu^\top$ and therefore 
$\cH^*(Y) - (1/\alpha) \cH^*(uu^\top) \in K$.  Hence 
$x - \cB(u,u)/\alpha = x - \cH^*(uu^\top)/\alpha \in K$,
so $(x,u,\alpha)$ is an element of the cone~(\ref{e-SC-def}).

To establish the equivalence between~(\ref{e-SC-alt}) 
and~(\ref{e-SC-def}) we only used Properties 1--3 
of Definition~\ref{def:K-bilinear-form}. 
Clearly, $\SC(K,\cB)$ has nonempty interior in $\cV\times \R^p \times \R$, 
since $K$ has nonempty interior in $\cV$.
Property~4 further implies that $\SC(K,\cB)$ is closed and pointed,
so it is a regular cone.
It is closed because $\SC(K,\cB)$ can be expressed as the image
of a closed convex cone $\SS_+^{p+1} \times K$ under the linear 
transformation
\[
 \cA(W, w) = (\cH^*(W_{22}) + w, W_{21}, W_{11})
\]
where $W_{11}$ is scalar, $W_{21} \in \R^p$, and 
$W_{22}$ is the trailing $p$-by-$p$ submatrix in
\[
 W = \left[\begin{array}{cc} W_{11} & W_{21}^\top \\ W_{21} & W_{22}
 \end{array}\right].
\]
Property~4 in its form~(\ref{e-cH-prop2-a}) implies that
$\cA(W,w) =0$, $W\succeq 0$, $w\in K$ only holds for $W=0$, $w=0$.
Hence, by Theorem 9.1 in \cite{Roc:70}, the set $\SC(K,\cB)$ is closed.
By a similar argument, $\SC(K,\cB)$ is pointed.  Suppose
$(x,u,\alpha)\in \SC(K,\cB)$ and $-(x,u,\alpha) \in \SC(K,\cB)$, so
\[
(x,u,\alpha) = \mathcal A(W,w), \qquad
-(x,u,\alpha) = \mathcal A(\tilde W,\tilde w), \qquad
\]
for some $W,\tilde W\succeq 0$, $w,\tilde w\in K$.
Therefore $0 = \mathcal A(W+\tilde W, w + \tilde w)$, and by 
Property~4, $W=\tilde W= 0$ and $w=\tilde w=0$. Hence $(x,u,\alpha) = 0$.

The dual cone of $\SC(K,\cB)$, if we use the inner product
$\langle s, x \rangle + 2v^\top u + \beta\alpha$ between $(s,v,\beta)$
and $(x,u,\alpha)$, is given by
\begin{equation} \label{e-SC*}
 \SC(K,\cB)^* = \left\{(s,v,\beta) \in \cV \times \R^p \times \R :
  s \in K^*, \left[\begin{array}{cc}
    \beta  & v^\top \\ v& \cH(s) \end{array}\right] \succeq 0 \right\}
\end{equation}
if $p\geq 1$.
If $\cB$ is the trivial bilinear form ($p=0$), 
the dual Siegel cone is $\SC(K,\cB)^* = K^* \times \R_+$.
The dual Siegel cone is closed, convex, and pointed, and
Property~4 in Definition~\ref{def:K-bilinear-form}
implies that it has nonempty interior.

Note that the expression~(\ref{e-SC*}) shows that $\SC(K,\cB)^*$ has 
a spectrahedral representation if $K^*$ has a spectrahedral
representation.

So far we have only used Properties 1--4 of
Definition~\ref{def:K-bilinear-form}.  Property~5 further
imples that the Siegel cone is a homogeneous cone.
This result is due to Vinberg \cite{Vinberg2}.
To see this, it is sufficient to verify that the group generated by 
the following linear transformations on $\cV \times \R^p \times \R$ forms
a transitive subset of $\aut(\SC(K,\cB))$:
\begin{eqnarray}
  T_1(\gamma) 
 & : &  (x,u,\alpha) \mapsto (x, \sqrt{\gamma}u, \gamma \alpha),  
 \label{e-T1} \\ 
  T_2(w) 
 & : &  (x,u,\alpha) \mapsto (x+2\cB(w,u)+ \alpha \cB(w,w),\, 
   u+\alpha w, \, \alpha),
\label{e-T2} \\
  T_3(g) 
 & : & (x,u,\alpha) \mapsto (g(x), \bar{g}u,\alpha).
\label{e-T3}
\end{eqnarray}
Here,
$T_1$ is parametrized by a scalar $\gamma > 0$,
$T_2$ by a vector $w\in\R^p$, and $T_3$ by an automorphism
$g\in G$, where $G$ is the transitive subset of $\aut(K)$ mentioned
in property~5 of Definition~\ref{def:K-bilinear-form}.
The mapping $\bar g$ is the corresponding linear transformation
in $\R^p$ and satisfies~(\ref{e-bilinear-prop5}).
It is easy to check, using~(\ref{e-SC-def}) or~(\ref{e-SC-alt}), 
that these transformations are automorphisms of $\SC(K,\cB)$.
To verify that they form a transitive subset, consider an arbitrary
pair of points $(x,u,\alpha)$ and $(\hat x, \hat u, \hat \alpha)$
in the interior of $\SC(K,\cB)$.   Let $g\in G$ be an automorphism
that maps $x- \cB(u,u)/\alpha$ to $\hat x - \cB(\hat u, \hat u)/
\hat \alpha$.   Then the mapping
\[
  T_1(\hat\alpha) \, \circ \, T_2(\hat u/\hat\alpha^{1/2}) 
  \, \circ \, T_3(g)
  \, \circ \, T_2(-u/\alpha^{1/2}) \, \circ \, T_1(1/\alpha)
\]
is an automorphism of $\SC(K,\cB)$ that maps 
$(x,u,\alpha)$ to $(\hat x, \hat u, \hat\alpha)$.

By duality, the adjoints of the mappings $T_1(\gamma)$,
$T_2(w)$, $T_3(g)$ form a transitive subset of the automorphism
group of $\SC(K,\cB)^*$.  The adjoints are given by
\begin{eqnarray}
  T_1(\gamma)^*
 & : &  (s,v,\beta) \mapsto (s, \sqrt{\gamma}v, \gamma \beta),  
\label{e-T1*} \\ 
  T_2(w)^* 
 & : &  (s,v,\beta) \mapsto (s, \, \cH(s)w + v, \, w^\top \cH(s)w + 
    2 w^\top v + \beta),
 \label{e-T2*} \\
  T_3(g)^* 
 & : & (s,v,\beta) \mapsto (g^*(s), \bar{g}^\top v,\beta)
\label{e-T3*}
\end{eqnarray}
and are exploited in the work of Rothaus \cite{Rothaus1966}.
These mappings correspond to congruence operations
\[
\left[\begin{array}{cc} \sqrt\gamma & 0 \\ 0 & I \end{array}\right]
 \left[\begin{array}{ccc}
  \beta & v^\top \\ v & \cH(s) \end{array}\right]
\left[\begin{array}{cc} \sqrt\gamma & 0 \\ 0 & I \end{array}\right], 
\qquad
\left[\begin{array}{cc} 1 & w^\top \\ 0 & I \end{array}\right]
 \left[\begin{array}{ccc}
  \beta & v^\top \\ v & \cH(s) \end{array}\right]
\left[\begin{array}{cc} 1 & 0 \\ w & I \end{array}\right],
\]
and
\[
\left[\begin{array}{cc} 1 & 0 \\ 0 & \bar g^\top \end{array}\right]
 \left[\begin{array}{ccc}
  \beta & v^\top \\ v & \cH(s) \end{array}\right]
\left[\begin{array}{cc} 1 & 0 \\ 0 & \bar g \end{array}\right],
\]
respectively, where on the last line we use the identity 
$\bar g^\top \cH(s)\bar g= \cH(g^*(s))$.

Rothaus calls the mapping $\cH$, associated 
with a homogenous $K$-bilinear symmetric form $\cB$ via the 
definition~(\ref{e-cB-cH}), a \emph{representation} of $K^*$,
and he calls the dual Siegel cone $\SC(K,\cB)^*$ an \emph{extension}
of $K^*$ \emph{from the representation} $\cH$ \cite{Rothaus1966}.
We have used the term \emph{dual representation} 
for~$\cH$, to avoid confusion with general semidefinite representations 
of convex cones (i.e., spectrahedral representations or
spectrahedral shadow representations).

\begin{examp}\label{examp:6.3}
We take $K = \R_+$ and $\cV = \R$, and the trivial bilinear form
($p=0$ and $\cB = 0$).  The Siegel cone is
\begin{eqnarray*}
 \SC(K,\cB) & = & \{(x,\alpha) \in \R \times \R: x\geq 0, \alpha\geq 0,
  \alpha x \geq 0 \}\\
& = & \{ (x,\alpha) \in \R\times\R : \left[\begin{array}{cc}
    \alpha & 0 \\ 0 & x \end{array}\right] \succeq 0\}.
\end{eqnarray*}
This cone is linearly isomorphic to $\R^2_+$.

For the same $K$ and $\cV$, consider $\cB:\R\times \R
\rightarrow \R$ defined as $\cB(u,v) = uv$.
The symmetric form clearly satisfies properties 1--4 in the 
Definition~\ref{def:K-bilinear-form}, and property~5 with
automorphisms $g(x) = \gamma x$ for $\gamma > 0$, 
and linear transformations $\bar g = \sqrt \gamma$.
The Siegel cone 
\begin{eqnarray*}
 \SC(K,\cB) & = & \{(x,u,\alpha) \in \R \times \R \times \R: 
 x\geq 0, \, \alpha\geq 0, \, \alpha x \geq u^2 \}\\
& = & \{ (x,u,\alpha) \in \R\times\R \times \R : \left[\begin{array}{cc}
    \alpha & u \\ u & x \end{array}\right] \succeq 0\}
\end{eqnarray*}
is linearly isomorphic to $\SS^2_+$.

Finally, consider $\cB : \R^p \times \R^p$ with $p>1$, defined 
as $\cB(u,v) = u^\top v$.
This is another $K$-bilinear homogenous form.  In property~5, we
take automorphisms $g(x) = \gamma x$ for $\gamma > 0$ 
and $\bar g = \sqrt \gamma I$.
The dual representation $\cH:\R \rightarrow \SS^p$ and its 
adjoint $\cH^* : \SS^p \rightarrow \R$ are 
$\cH(s) = sI$ and $\cH^*(Y) = \Tr(Y)$.
With this choice of $\cB$, we obtain 
\begin{eqnarray*}
 \SC(K,\cB) 
& = & \{(x,u, \alpha) \in \R \times \R^p \times \R: 
 x\geq 0, \, \alpha\geq 0, \, \alpha x \geq u^\top u\} \\
& = & \{ (x,u,\alpha) \in \R\times\R^p \times \R : \left[\begin{array}{cc}
    \alpha I & u \\ u^\top & x \end{array}\right] \succeq 0\}.
\end{eqnarray*}
This cone is linearly isomorphic to the rotated quadratic 
cone~(\ref{e-rot-quad}).

The most general form of a homogeneous $K$-bilinear symmetric form
for $K=\R^+$ is $\cB(u,v) = u^TQ^{-1}v$, with $Q\in \SS^p_{++}$.
With this choice, 
\begin{eqnarray*}
 \SC(K,\cB) 
& = & \{(x,u, \alpha) \in \R \times \R^p \times \R: 
 x\geq 0, \, \alpha\geq 0, \, \alpha x \geq u^\top Q^{-1} u\} \\
& = & \{ (x,u,\alpha) \in \R\times\R^p \times \R : \left[\begin{array}{cc}
    \alpha Q & u \\ u^\top & x \end{array}\right] \succeq 0\}.
\end{eqnarray*}
This cone is a rotated quadratic cone after a linear transformation.

Hence, $\R^2_+$, $\SS^2_+$, and the rotated
quadratic cones are essentially the only types of cones that can be
constructed as Siegel cones of $\R_+$.

\end{examp}

\begin{examp}\label{examp:6.4}
We continue Examples~\ref{examp:6.1} and~\ref{examp:6.2}.
The Siegel cone for the bilinear form in Example~\ref{examp:6.1} is
\[
 \SC(K,\cB) = \{(x,u,\alpha) \in \SS^3_E \times \R^6 \times \R: 
 \left[\begin{array}{cccccc}
  \alpha & 0 & 0 & u_1 & 0 & u_2 \\
  0 & \alpha & 0 & u_3 & 0 & u_4 \\
  0 & 0 & \alpha & 0  & u_5 & u_6 \\
  u_1 & u_3 & 0 & X_{11} & 0 & X_{31} \\
  0   & 0   & u_5 & 0 & X_{22} & X_{32} \\
  u_2 & u_4 & u_6 & X_{31} & X_{32} & X_{33} \end{array}\right]
 \succeq 0 \}.
\]
This is an example of a \emph{sparse matrix norm cone} discussed in
Section~\ref{s-hom-cones-examples}.
For the bilinear form of Example~\ref{examp:6.2}, we obtain 
\[
 \SC(K,\cB) = \{(x,u,\alpha) \in \SS^3_E \times \R^2 \times \R: 
 \left[\begin{array}{cccc}
  \alpha & u_1 & 0 & u_2 \\
  u_1 & X_{11} & 0 & X_{31} \\
  0   & 0 & X_{22} & X_{32} \\
  u_2 & X_{31} & X_{32} & X_{33} \end{array}\right]
 \succeq 0 \}.
\]
This is a homogeneous sparse matrix cone (ordered using a trivially
perfect elimination ordering).
\end{examp}

\subsection{Siegel domain construction of homogeneous cones}
The Siegel cone is the key tool in Vinberg's recursive construction
of all homogeneous cones.  A homogeneous cone $K$
and a homogeneous $K$-bilinear symmetric form $\cB$ together
yield a Siegel cone $\SC(K,\cB)$ which is a homogeneous cone in a 
higher dimensional space.
The converse is also true.  For every homogeneous cone $\widehat K$
of dimension at least 2,
there exists a lower dimensional homogeneous cone $K$ and a 
homogeneous $K$-bilinear symmetric form $\cB$ such that $\widehat K$ 
is linearly isomorphic to the Siegel cone $\SC(K,\cB)$; 
see for example Rothaus~\cite{Rothaus1966} or Gindikin~\cite{Gindikin1992}.
This provides an inductive characterization (called 
\emph{Siegel domain construction}) of all homogeneous cones, 
starting with the ray $\R_+$ in $\R$ as the first homogenous cone. 
The construction may be viewed as an abstraction of the 
commonly used concept of \emph{Schur complement}.

The minimum number of steps required to construct
$K$ in this recursive way is called the \emph{Siegel rank of $K$}.
We denote this integer valued function of homogeneous cones by
$\Srank(K)$ and define $\Srank(\R_+):=1$.
Since $K$ is homogeneous if and only if  $K^*$ is,
Vinberg's classification theory described above also applies to $K^*$. 
Furthermore, the Siegel ranks of $K$ and $K^*$ are always the same.

\begin{examp}\label{examp:6.5}
The Vinberg cone
\begin{equation} \label{e-vinberg}
K:=\left\{x \in \R^5: 
\begin{bmatrix} x_1 &  0 & x_2 \\0 &  x_3 & x_4\\ x_2 & x_4 & x_5
\end{bmatrix} \succeq 0 \right\}
\end{equation}
has Siegel rank 3. Let us construct it via the recursive procedure. 

\begin{itemize}
\item Let $K_1=\R_+$. Define
$\cB_1: \R \times \R \to \R$ by $\cB(u,v):= uv$. Then, $K_1$ and $\cB_1$ 
satisfy the conditions in Definition~\ref{def:K-bilinear-form} and
their Siegel cone is linearly isomorphic to $\SS_2^+$;
see Example~\ref{examp:6.3}.
This shows that $\Srank(\SS^2_+) = 2$. 
\item Let $K_2:=\SS^2_+$. Define $\cB_2 : \R \times \R \to \SS^2$ by
\[
\cB_2(u,v) := \frac{1}{2}\begin{bmatrix} 0 \\u \end{bmatrix}
\begin{bmatrix} 0 & v \end{bmatrix}
+
\frac{1}{2}\begin{bmatrix} 0 \\v \end{bmatrix}
\begin{bmatrix} 0 & u \end{bmatrix}
=
\begin{bmatrix} 0 & 0 \\ 0 & uv \end{bmatrix}.
\]
Again, $K_2$ and $\cB_2$ satisfy the conditions
in Definition~\ref{def:K-bilinear-form}. 
To check the fifth condition, we can take as the transitive subset $G$ 
the set of linear maps $g(X) = L  X L^{\top}$,  
where $L$ is a nonsingular  2-by-2 lower-triangular matrix
\[
L = \begin{bmatrix} L_{11} & 0 \\ L_{21} & L_{22} \end{bmatrix}, 
\]
and define $\bar g= L_{22}$, so that  
$L(\cB(u,v)L^\top = \cB(L_{22}u, L_{22}v)$ as desired.
The Siegel cone for $K_2$ and the bilinear form $\cB_2$ is 
\[
\SC(K_2,\cB_2) = \{(X, u, \alpha) \in \SS^2 \times \R \times \R
 : \left[\begin{array}{ccc} \alpha  & 0 & u \\
      0  & X_{11} & X_{21} \\
      u  & X_{21} & X_{21} \end{array}\right] \succeq 0 \},
\]
which is linearly isomorphic to the Vinberg cone~(\ref{e-vinberg}).
Thus, we have derived the Vinberg cone as a homogeneous cone with 
$\Srank(K) = 3$.
\end{itemize}

\end{examp}

\subsection{Semidefinite representations of homogeneous cones}
Since every homogeneous cone of Siegel rank $r\geq 2$ arises from
a homogeneous cone of Siegel rank $r-1$, via the above construction, we can
establish many properties of homogeneous cones by induction on the
Siegel rank.   For example, from~(\ref{e-SC-alt}) we see that
$\SC(K,\cB)$ has a spectrahedral shadow representation if $K$ has
a spectrahedral shadow representation.
By induction, starting with $\R_+$, it follows that every 
homogeneous cone has a spectrahedral shadow or lifted-LMI representation.

From~(\ref{e-SC*}), we also see that if the dual of a homogenous 
cone $K$ has a  spectrahedral representation then so does 
$\SC(K,\cB)^*$.  
For example, if $\cV = \R^n$ and 
$K^* = \{ s \in \R^n: \mathcal A(s) \succeq 0\}$
is an LMI representation of $K^*$, then 
\[
\SC(K,\cB)^* = \left\{ (s,v,\beta) \in \R^n \times \R^p \times \R: 
   \left[\begin{array}{ccc}
   \beta & v^\top & 0 \\ v & \cH(s) & 0 \\ 0 & 0 & \mathcal A(s)
 \end{array}\right] \succeq 0 \right\}
\]
is an LMI representation of $\SC(K,\cB)^*$.
Since the set of 
homogeneous cones is closed under duality, every homogeneous cone of 
Siegel rank $r \geq 2$ must arise as $\SC(K, \cB)^*$ for some homogeneous 
cone $K$ of Siegel rank $r-1$ and some homogeneous $K$-bilinear symmetric 
form $\cB$.
Therefore, by induction on $r$, we can establish that every homogeneous 
cone has a spectrahedral representation
of the form $\cV \cap \SS_+^N $ for some $N \geq 1$ and 
a linear subspace $\cV$ of $\SS^N$.
Moreover, it can be assumed that there exists a transitive subset
of $\aut(\cV\cap \SS^N_+)$, consisting of congruences
$\mathcal R(U)= R^{\top}UR$.  This again follows by induction 
from~(\ref{e-SC*}) and the fact that the group generated by 
the mappings~(\ref{e-T1*})--(\ref{e-T3*}) forms a transitive subset
of $\aut(\SC(K,\cB)^*)$.

This high-level description of a recursive construction of a
spectrahedral representation does not necessarily lead to an
efficient representation.
In Section~\ref{s-hom-matrix-cones}, we saw a more structured and 
potentially more efficient canonical form for the spectahedral 
representation of a homogeneous cone, due to Ishi \cite{Ishi2015}.

In Example~\ref{examp:6.3} we enumerated the three types of cones
(up to linear isomorphisms) that can be constructed as Siegel cones
of $K=\R_+$ by choosing different bilinear forms $\cB$.
The possible homogeneous  $K$-bilinear symmetric forms are the trivial
symmetric form and the inner products $\cB(u,v) = v^\top u$ in $\R^p$,
for $p\geq 1$.  The corresponding dual representations are
$\cH(s) = sI$, where $I$ is a $p$-by-$p$ identity matrix.
To continue the inductive construction, one needs to 
characterize all homogenous $\widehat K$-bilinear symmetric forms 
$\widehat{\cB}$ for a Siegel cone $\widehat K = \SC(K,\cB)$, 
or, equivalently, the  dual representations $\widehat{\cH}$ of 
$\widehat{\cB}$, from the dual representations $\cH$ of $\cB$.
By definition, 
\[
\widehat{\cH}(s, v, \beta) = \left[\begin{array}{cc} 
  \beta & v^\top \\ v & \cH(s)
\end{array}\right]
\]
is a dual representation of one possible $\widehat{\cB}$. 
Rothaus has characterized the possible dual representations
of homogeneous $\widehat K$-bilinear symmetric 
forms~\cite{Rothaus1963,Rothaus1966,Rothaus1968}.
Using properties 1--5 of Definition~\ref{def:K-bilinear-form} 
he  proves \cite[Lemma 3.5 and Theorem 3.7]{Rothaus1966}
that they are all of the form
\[
\widehat{\cH}(s, v, \beta)
 = \left[\begin{array}{cc} \beta I & \cC(v)^\top \\ 
 \cC(v) & \widetilde \cH(s) \end{array}\right], 
\]
where $\widetilde \cH: \cV \rightarrow \SS^q$ is the dual representation
of a homogeneous $K$-bilinear form $\widetilde{\cB} : 
\R^q \times \R^q \rightarrow \cV$,
and $\cC: \R^p \rightarrow \R^{q \times m}$ is a linear matrix function.
Property~5 in Definition~\ref{def:K-bilinear-form} imposes
additional constraints on  $\mathcal C$ and $\widetilde H$
(specifically, the matrix equations 
\[
\cC(u)^\top \cC(v)  + \cC(v)^\top \cC(u) = 2(u^\top v) I,
\qquad  
\cC(\cH(s)v) = \widetilde{\cH}(s) \cC(v)
\]
are satisfied for all $s,u,v$).

\section{Primal--dual interior-point methods} \label{s-ipm}
Theorems~\ref{thm:asymmetric-factor} and \ref{thm:6.3} provide 
a technique for ``scaling'' primal and dual conic optimization problems 
over homogeneous matrix cones. 
These types of primal--dual scalings can be used to design and
analyze \emph{scale-invariant} and \emph{primal--dual symmetric} 
algorithms in the sense of \cite{Tun:98}. 

Consider a linear conic optimization problem in standard form (P) 
and its dual~(D),
with cones $K$ and $K^*$ in a finite-dimensional real vector space $\cW$: 
\begin{equation} \label{e-conic-LP}
\mbox{($\mathrm{P}$)} \quad
\begin{array}[t]{ll}
\mbox{minimize} & \iprod{c}{x} \\
\mbox{subject to} & \cA(x)=b \\ 
   &  x \in K 
\end{array}
\qquad\qquad
\mbox{($\mathrm{D}$)} \quad
\begin{array}[t]{ll}
\mbox{maximize} & b^\top y \\
\mbox{subject to} & \cA^*(y) +s = c \\
 &  s \in K^*.
\end{array}
\end{equation}
The linear mapping $\mathcal A: \cW\rightarrow \R^m$,
and the vectors $b\in\R^m$ and $c\in\cW$ are given. 
The primal optimization variable is $x\in\cW$, the dual variables
are $y\in\R^m$ and $s\in\cW$.
The linear conic optimization problem is a natural extension of linear 
programming (the special case with $\cW=\R^n$, $K=K^*=\R^n_+$), and
has been widely used in the development of interior-point 
methods; see the surveys in \cite{Renegar2001,NemTodd2008,BeN:01}.
Advances in conic optimization algorithms and software have also
enabled the creation of highly influential modeling software
for convex optimization  \cite{yalmip,cvx,cvxpy,cvxr,PICOS}.  
These modeling tools take advantage of the fact  that a few 
different types of convex cones 
(the positive semidefinite cone, the second order cone, and the
exponential cone) are sufficient to reformulate most convex optimization 
problems encountered in practice as conic optimization problems.

In this section we first review some recent literature on 
interior-point algorithms for conic optimization, 
and then comment on the special case of homogeneous cones $K$ and $K^*$.

Until 2001, the research and the literature on primal--dual interior-point 
methods with polynomial iteration complexity were dominated by 
a very high level of activity concentrated on 
semidefinite programming, i.e., the special case $\cW=\SS^N$, 
$K=K^*=\SS^N_+$.  We must note, however, that Nesterov
and Nemirovski's 1994 monograph~\cite{NN1994} already contained 
primal--dual interior-point algorithms, with polynomial iteration 
complexity, for general conic programming.  Their results were based on 
the primal--dual symmetric, generalized Tanabe--Todd--Ye potential 
function.  Generalized, because the original Tanabe--Todd--Ye potential 
function was proposed for linear programming and linear complementarity 
problems, whereas Nesterov and Nemirovski's generalization 
replaced the logarithmic barrier 
functions for the nonnegative orthant with 
\emph{self-concordant} barrier functions for $K$, 
and used the Legendre--Fenchel conjugates for the dual cone $K^*$.
Even though the algorithm and 
analysis were based on a primal--dual symmetric potential function,
the underlying algorithm was not primal--dual symmetric (the algorithm 
chose different kinds of search directions and steps in the primal 
and dual spaces).

In a breakthrough work, Nesterov and Todd~\cite{NT1997,NT1998} identified
properties of self-concordant barriers (they called such special 
self-concordant barriers \emph{self-scaled}) allowing the design and 
analysis of primal--dual symmetric interior-point algorithms with an 
outstanding number of desired properties and matching the best iteration
complexity bounds. Only \emph{symmetric cones} (those that are 
homogeneous and self-dual) admit self-scaled barriers. 
Therefore, the Nesterov--Todd algorithms only apply to symmetric cones,
i.e., they are limited to second order cone programming and semidefinite 
programming (over symmetric matrices with real entries, Hermitian 
matrices with complex entries or quaternion entries, and 
Hermitian $3$-by-$3$ matrices over the octonions).
It quickly became clear that generalizing all of the desired properties 
of Nesterov--Todd algorithms beyond symmetric cones was impossible 
(see a result of Nesterov in Theorem~7.2 of \cite{Guler1997}, \cite{Tun:98} and Lemma~6.4 of \cite{NT2016}).
However, as explained below, Theorem~\ref{thm:6.3} does 
allow some possibilities for homogeneous cones that are not self-dual.

Until recently, the literature on primal--dual symmetric interior-point 
algorithms for general conic optimization was quite sparse
(beyond the special case of symmetric cones), but papers on this subject
have been increasing in number and in their depth.
A general framework for primal--dual symmetric interior-point algorithms 
with polynomial iteration complexity was proposed in 
2001~\cite{Tuncel2001}. The same paper also showed how the theory of 
quasi-Newton methods and quasi-Newton-like updates can be applied
to the computation of a primal--dual scaling in interior-point methods.
Chares~\cite{Chares2009} considered $p$-norm and power cone optimization 
problems, and recently Roy and Xiao~\cite{RoyXiao2022} proved Chares's 
conjecture on self-concordance of a very efficient barrier function for generalized power cones.
Nesterov~\cite{Nesterov2012}
proposed primal--dual interior-point algorithms for general conic optimization
which are based on a primal--dual scaling 
that approximately satisfies
the conditions used in \cite{NT1997,NT1998,Tuncel2001}. Myklebust and Tun\c{c}el~\cite{MT2014}
streamlined the computation of the primal--dual scaling in 
\cite{Tuncel2001}, which involved a rank-four update to a symmetric 
positive definite matrix, by expressing it as a composition
of two rank-two quasi-Newton updates. Further, they proved that short-step
path following algorithms based on the framework of~\cite{Tuncel2001} 
achieve the same worst case polynomial iteration complexity as the 
current best interior-point algorithms for symmetric cone programming.
Skajaa and Ye~\cite{SkajaaYe2015} also used  the idea of quasi-Newton 
methods to design and analyze a primal--dual interior-point algorithm 
for general conic optimization. 
Their algorithm with some necessary modifications has been implemented 
by Papp and Yildiz~\cite{PappY2022}.  
Dahl and Andersen~\cite{DahlAndersen2022} followed the 
framework from~\cite{Tuncel2001, MT2014}, 
made a connection to Schnabel's work~\cite{Schnabel}, 
and designed and implemented primal--dual symmetric interior-point
algorithms for exponential cone programming problems. 

There are interior-point algorithms for hyperbolicity 
cones~\cite{RenegarS2014,NT2016}, and there is further potential for 
interesting primal--dual algorithms utilizing self concordant barriers 
for sophisticated matrix cones such as those arising from the quantum 
relative entropy~\cite{FaybusovichZ2022,KarimiT2019,FawziS2022}.

Chua~\cite{Chua2009} proposed a primal--dual interior-point algorithm for 
conic optimization with general homogeneous cones, 
based on Vinberg's axioms and exploiting the underlying structure,
including the transitive subset of the cone automorphism group.
Chua's algorithm achieves the current best iteration complexity bound for 
symmetric cone programming.
However, it is not primal--dual symmetric 
(in each iteration, the scaling is computed based on the
automorphism which maps the current dual iterate to identity).
Moreover, the search direction is well-defined only in a narrow
neighbourhood of the central path.

Through computational experiments several advantages of primal--dual
algorithms have been observed. In theoretical contexts, there are
additional justifications. See \cite{NesterovTodd2002} for a justification of the usage of
primal--dual central path setting; and see \cite{Todd2009} for a geometric
justification of the primal--dual scaling in the case of symmetric cones.

Next, based on the results and insights from the earlier sections,
we outline a new way of computing the primal--dual scaling
in new primal--dual symmetric interior-point algorithms for homogeneous
cones.

We assume the cone $K$ in~(\ref{e-conic-LP}) is a homogeneous matrix
cone of the types discussed in
Sections~\ref{s-cone}--\ref{s-hom-matrix-cones}, i.e.,
a homogeneous sparse matrix cone~(\ref{e-K-Kstar}) or the more 
general homogeneous matrix cone~(\ref{e-K-V}).    
In the first case we take $\cW = \SS^N_E$ in~(\ref{e-conic-LP});
in the second case $\cW = \cV$.  We denote by $F$ the logarithmic barrier 
function~(\ref{e-log-det-barrier}) and~(\ref{e-log-det-barrier-hom}).  
This is a \emph{$\vartheta$-logarithmically homogeneous self-concordant 
barrier} (or \emph{$\vartheta$-normal barrier}), with parameter $\vartheta=N$.
To simplify the notation and the application to homogeneous cones 
in other representations (e.g., with $\mathcal W=\R^p$),  
we continue to use lower-case symbols $x, s$ for the variables.

Let $x_k \in \inte(K)$, $s_k \in \inte(K^*)$ be current iterates in a
primal--dual algorithm for~(\ref{e-conic-LP}).
Theorems~\ref{thm:asymmetric-factor} 
and \ref{thm:6.3} state that there exists an automorphism 
$\mathcal{L}$ of $K$ that satisfies
\begin{equation} \label{e-L-scaling}
 \mathcal L^{-1}(x_k) = \mathcal L^{*}(s_k).
\end{equation}
Moreover $\mathcal L^{-*} \circ \mathcal L^{-1} = F''(w)$, where
$w$ is the \emph{primal--dual scaling point} defined by 
\begin{equation} \label{e-pd-scaling}
 F''(w; x_k) = s_k. 
\end{equation}
(If $\mathcal W = \R^p$, this equation is written more simply as 
$F''(w)x_k = s_k$.)
If we make a change of variables $\bar x = \cL^{-1}(x)$, 
$\bar s = \cL^*(s)$, problems (P) and (D) in~(\ref{e-conic-LP}) are 
transformed to 
$(\bar{\textup{P}})$ and $(\bar{\textup{D}})$ given below:
\[
\mbox{($\bar{\mathrm{P}}$)} \quad
\begin{array}[t]{ll}
\mbox{minimize} & \iprod{\bar{c}}{\bar x} \\
\mbox{subject to} & \bar{\cA}(\bar x)=b \\ 
   &  \bar x \in K 
\end{array}
\qquad\qquad
\mbox{($\bar{\mathrm{D}}$)} \quad
\begin{array}[t]{ll}
\mbox{maximize} & b^\top y \\
\mbox{subject to} & \bar{\cA}^* (y) +\bar s = \bar{c} \\
 &  \bar s \in K^*.
\end{array}
\quad
\]
Here, $\bar{\cA} := A \circ \mathcal{L}$ and 
$\bar{c} := \mathcal{L}^*(c)$. 
In the scaled problem the current iterates $x_k$, $s_k$ are mapped to the 
same point 
\begin{equation} \label{e-vk-def}
 v_k := \cL^{-1}(x_k) = \cL^*(s_k).
\end{equation}
The scaled problem generalizes the \emph{$v$-space formulation} from 
the interior-point literature for linear complementarity problems and 
linear optimization problems over symmetric cones
\cite{KMNY:91,JRT:96,StZ:99}.
Extending the definition in~\cite[Section 3]{Tuncel2001}, we can 
define a \emph{primal--dual affine scaling direction} at $x_k, s_k$ 
as the solution $(d_x, d_y, d_s)$ of the linear system
\[
\cA(d_x) = 0, \qquad
{\cA}^*(d_y) + d_s = 0, \qquad
\cL^{-1}(d_x) + \cL^*(d_s) = -v_k.
\]

If $F$ is a self-scaled barrier of a symmetric cone, the 
equation~(\ref{e-pd-scaling}) defines the Nesterov--Todd scaling point 
$w$~\cite{NT1997,NT1998}, and automatically implies 
\begin{equation} \label{e-pd-scaling-shadow}
F''(w; \tilde x_k) = \tilde s_k, 
\end{equation}
where $\tilde x_k := -F'_*(s_k)$ and $\tilde s_k := -F'(x_k)$.   
For general convex cones, the equation~(\ref{e-pd-scaling}) can still
be used to define primal--dual scaling points, and algorithms based on such
scalings have been studied in \cite{Tuncel2001,Nesterov2012}.
An important difference is that for general cones
(and for the non-self-dual homogeneous cones discussed in this paper) 
the equation~(\ref{e-pd-scaling})
does not imply~(\ref{e-pd-scaling-shadow}).
In the algorithms of \cite{Tuncel2001,MT2014} this difficulty is
addressed by making a rank-four update to $F''(w)$, or to an approximation
of $F''(w)$, to define a positive definite self-adjoint mapping 
$\mathcal H$ that satisfies both $\mathcal H(x_k) = s_k$ and 
$\mathcal H(\tilde x_k) = \tilde s_k$.  
In these algorithms,
a primal--dual search direction $(d_x,d_y,d_s)$ is computed from 
$\mathcal H$ by solving the equation
\begin{equation} \label{e-myklebust}
\cA(d_x) = 0, \; \quad
{\cA}^*(d_y) + d_s = 0, \; \quad
\mathcal H^{-1/2}(d_x) + \mathcal H^{1/2}(d_s) = 
-v_k + \gamma\mu \tilde v_k,
\end{equation}
where $\gamma \in[0,1]$ is a centering parameter,
$\mu = \langle s_k, x_k\rangle/\vartheta$, and
\[
v_k := \mathcal H^{1/2}(x_k) = \mathcal H^{-1/2}(s_k), \qquad
\tilde v_k := \mathcal H^{1/2}(\tilde x_k) = \mathcal H^{-1/2}(\tilde s_k).
\]

For homogeneous cones, different and simpler updates are possible, 
because the factor $\cL$ of the primal--dual scaling matrix 
$F''(w) =  \cL^{-*} \circ \cL^{-1}$ can be modified by a rank-one update 
to obtain a scaling ${\mathcal L}_+$ that satisfies the two conditions
\begin{equation} \label{e-Lbar-eqs}
{\mathcal L}_+^{-1}(x_k) = {\mathcal L}_+^{*}(s_k), \qquad 
{\mathcal L}_+^{-1}(\tilde x_k) = {\mathcal L}_+^{*}(\tilde s_k).
\end{equation}
The construction of $\mathcal L_+$ is similar to a quasi-Newton update,
with the difference that ${\mathcal L}_+$ must satisfy the 
two equations~(\ref{e-Lbar-eqs}), as opposed to one (secant) equation
in standard quasi-Newton updates.
A simplification of the updates in \cite{MT2014} that 
achieves this goal proceeds as follows.
Define 
\[
\delta_\mathrm p :=  x_k - \mu \tilde x_k, \qquad
\delta_\mathrm d :=  s_k - \mu \tilde s_k
\]
where $\mu = \langle s_k, x_k\rangle/\vartheta$.
General properties of $\vartheta$-logarithmically homogeneous barriers
($\langle F'(x), x\rangle = \langle s, F'_*(s)\rangle = -\vartheta$)
imply that
\[
\langle s_k, \delta_\mathrm p \rangle = 
\langle \delta_\mathrm d, x_k \rangle =  0.
\]
One can also show that
\[
\langle \delta_\mathrm d, \delta_\mathrm p\rangle \geq  0 
\]
with equality only if $\delta_\mathrm p$ and $\delta_\mathrm d$ are 
both zero \cite[Corollary 4.1]{Tuncel2001}.
We note the simple $v$-space expressions
\[
 \cL^{-1}(\delta_\mathrm p) = v_k + \mu F'_*(v_k), \qquad
 \cL^*(\delta_\mathrm d) = v_k + \mu F'(v_k), 
\]
which follow from~(\ref{e-vk-def}) and the composition 
properties~(\ref{e-composition}),~(\ref{e-composition-dual}).
The second terms on the two right-hand sides are not equal because
$F'_* \neq F'$, unless  the cone is symmetric.
The purpose of the update of $\mathcal L$ is to achieve
${\mathcal L}_+^{-1}(\delta_\mathrm p) = 
 {\mathcal L}_+^{*}(\delta_\mathrm d)$ 
while preserving ${\mathcal L}_+^{-1}(x_k) = {\mathcal L}_+^{*}(s_k)$.

If $\delta_\mathrm p = \delta_\mathrm d =0$, 
no update is needed and we take $\mathcal  L_+ = \cL$.
Otherwise,  $\langle \delta_\mathrm d, \delta_\mathrm p\rangle > 0$ 
and we use the (Broyden) rank-one  update 
\begin{equation} \label{e-cL-bfgs}
{\mathcal L}_+(\cdot) = \mathcal L(\cdot)
   + \frac{\langle \hat v_k, \cdot\rangle}{\|\hat v_k\|^2}
  (\delta_\mathrm p - \mathcal L(\hat v_k)), 
\end{equation}
where $\hat v_k$ is a multiple of $\cL^*(\delta_\mathrm d)$,
scaled to have norm $\|\hat v_k\| = \langle \delta_\mathrm  p,
\delta_\mathrm d\rangle^{1/2}$, i.e.,
\begin{equation} \label{e-hatvk-bfgs}
\hat v_k = \frac{1}{\alpha} \cL^*(\delta_\mathrm d), \qquad
\alpha = 
\frac {\|\cL^*(\delta_\mathrm d)\|} 
{\langle \delta_\mathrm d, \delta_\mathrm p\rangle^{1/2}}.
\end{equation}
The mapping $\mathcal L_+$ is invertible with inverse 
\[
{\mathcal L}_+^{-1}(\cdot) = \mathcal L^{-1}(\cdot) +  \alpha
\frac{\langle \mathcal L^{-*}(\hat v_k), \cdot\rangle}
{\|\hat v_k\|^2}
 (\hat v_k-\mathcal L^{-1}(\delta_\mathrm p)).
\]
We verify that the update $\mathcal L_+$ satisfies
\begin{equation} \label{e-scaling-rank-1}
{\mathcal L}_+^{-1}(x_k) = {\mathcal L}_+^{*}(s_k) = v_k, \qquad
{\mathcal L}_+^{-1}(\delta_\mathrm p) = 
 {\mathcal L}_+^{*}(\delta_\mathrm d) = \hat v_k. 
\end{equation}
This is equivalent to~(\ref{e-Lbar-eqs}) because
$\delta_\mathrm p$ and $\delta_\mathrm d$ are linear combinations
of $x_k, \tilde x_k$ and $s_k,\tilde s_k$, respectively.
To show~(\ref{e-scaling-rank-1}) we first note that
\[
\langle \hat v_k, v_k\rangle
= \frac{1}{\alpha} \langle \cL^*(\delta_\mathrm d), v_k\rangle
= \frac{1}{\alpha} \langle \delta_\mathrm d, \cL(v_k)\rangle
= \frac{1}{\alpha} \langle \delta_\mathrm d, x_k\rangle
= 0.
\]
Therefore, applying~(\ref{e-cL-bfgs}) to $v_k$ gives
\[
 \mathcal L_+(v_k) = \mathcal L(v_k)
   + \frac{\langle \hat v_k, v_k\rangle}{\|\hat v_k\|^2}
  (\delta_\mathrm p - \mathcal L(\hat v_k)) 
= \mathcal L(v_k) = x_k.
\]
Applying the adjoint to $s_k$ gives
\begin{eqnarray*}
\mathcal L_+^*(s_k)
& = & \mathcal L^*(s_k) + \frac{\langle s_k, \delta_\mathrm p - 
 \cL(\hat v_k)\rangle}{\|\hat v_k\|^2} \hat v_k \\
& = & v_k  - \frac{\langle \mathcal L^*(s_k), 
 \hat v_k\rangle}{\|\hat v_k\|^2} \hat v_k \\
& = & v_k  - \frac{\langle v_k,
 \hat v_k\rangle}{\|\hat v_k\|^2} \hat v_k \\
& = & v_k.
\end{eqnarray*}
This proves the first two equations in~(\ref{e-scaling-rank-1}).
The equation $\mathcal L_+(\hat v_k) = \delta_\mathrm p$
is immediate from~(\ref{e-cL-bfgs}).
The last equation $\mathcal L^*_+(\delta_\mathrm d) = \hat v_k$ 
follows from
\begin{eqnarray*}
\mathcal L_+^*(\delta_\mathrm d)
& = & \mathcal L^*(\delta_\mathrm d) + \frac{\langle \delta_\mathrm d, 
 \delta_\mathrm p\rangle - 
 \langle \mathcal L^*(\delta_\mathrm d), \hat v_k\rangle}{\|\hat v_k\|^2} 
 \hat v_k\\
& = & \alpha \hat v_k  + \frac{\|\hat v_k\|^2 -\alpha\|\hat v_k\|^2}
{\|\hat v_k\|^2} \hat v_k \\
& = & \hat v_k .
\end{eqnarray*}
The rank-one update~(\ref{e-cL-bfgs}) is the square-root form
of a Broyden--Fletcher--Goldfarb--Shanno (BFGS) update
(i.e., $\mathcal H_+ = \mathcal L_+ \circ \mathcal L_+^*$ is the BFGS 
update of $\mathcal H= \cL \circ \cL^*$ for the secant condition 
$\mathcal H_+(\delta_\mathrm d) = \delta_\mathrm p$).  
Many other rank-one updates will serve the same purpose.
For example, in \cite[p.140]{Sor:82} a closely related family of 
quasi-Newton updates is defined.  Sorensen's 
updates are parametrized by the vector $\hat v_k$ (in our notation). 
Instead of~(\ref{e-hatvk-bfgs}) one can choose for $\hat v_k$ 
any vector that satisfies 
\begin{equation} \label{e-hatv-conds}
 \|\hat v_k\|^2 = \langle \delta_\mathrm d, \delta_\mathrm p\rangle,
 \qquad
 \langle \hat v_k, v\rangle = 0, 
\end{equation}
and
$\langle \mathcal L^{*}(\delta_\mathrm d), \hat v_k\rangle
 \neq \langle \delta_\mathrm d, \delta_\mathrm p\rangle$, 
$\langle \hat v_k, \mathcal L^{-1}(\delta_\mathrm p) \rangle
 \neq \langle \delta_\mathrm d, \delta_\mathrm p\rangle$. 
Then the mapping ${\mathcal L}_+$ defined by
\[
{\mathcal L}_+(\cdot)
 = \mathcal L(\cdot)
  + \frac{\langle \mathcal L^{*}(\delta_\mathrm d) - \hat v, \cdot\rangle}
 {\langle \mathcal L^{*}(\delta_\mathrm d) - \hat v, \hat v\rangle}
   (\delta_\mathrm p - \mathcal L(\hat v))
\]
is invertible and satisfies~(\ref{e-scaling-rank-1}).
When $\alpha\neq 1$, the update~(\ref{e-cL-bfgs}) is a special
case if we choose the $\hat v_k$ given in~(\ref{e-hatvk-bfgs}).

Corresponding to the updated primal--dual scaling $\mathcal L_+$
that satisfies~(\ref{e-scaling-rank-1}), a 
primal--dual search direction at $x_k$, $s_k$ can be defined as the 
solution of the equation
\[
\cA(d_x) = 0, \qquad
{\cA}^*(d_y) + d_s = 0, \qquad
\mathcal L_+^{-1}(d_x) + \mathcal L_+^*(d_s) = 
-v_k + \gamma\mu \tilde v_k,
\]
where $\gamma \in[0,1]$ is a centering parameter,
$\mu = \langle s_k, x_k\rangle/\vartheta$, and
\[
\tilde v_k 
:= \mathcal L_+^{-1}(\tilde x_k) = \mathcal L_+^*(\tilde s_k)
= \frac{1}{\mu} (v_k - \hat v_k).
\]
This primal--dual search direction simplifies the algorithms developed and 
analyzed in the framework of \cite{Tuncel2001}, and are similar to 
the algorithms in \cite{MT2014},
based on the search direction defined in~(\ref{e-myklebust}).
However, the primal--dual scalings for homogeneous cones described
above have stronger properties than 
the primal--dual scalings for general convex cones discussed 
in \cite{MT2014}.
As we mentioned earlier,
instead of a rank-four update of the scaling matrix $F''(w)$, 
we perform a rank-one update to the factors (each of which is an automorphism
of the underlying cone) in the decomposition
$F''(w) = \cL^{-*} \circ \cL^{-1}$ 
to satisfy the second of the key equations~(\ref{e-Lbar-eqs}).
Moreover, if we apply a short-step strategy satisfying the assumptions
in~\cite{MT2014}, these new algorithms achieve the same polynomial
time iteration complexity as the current best primal--dual symmetric
interior-point algorithms for symmetric cone programming.

Our approach above offers more possibilities for the design and analysis of
algorithms which have a significant part operating in the $v$-space. In addition
to the references for $v$-space based algorithms we mentioned above, another
example is the algorithm for linear programming proposed by~\cite{Nesterov2008}.

\section{Conclusion}
\label{s-conclusions}
Special cases of homogeneous matrix cones have been studied in the
conic optimization literature.
Sparse SDPs with arrow patterns are quite common, and arise, for example,
in robust least squares and robust quadratic programming 
\cite{AVD:10,ElL:97,BEN:09}, and in structural optimization~\cite{Koc:21}.
They also appear in semidefinite relaxations of optimization problem 
with quadratic equality constraints,
when the constraints involve only squares $x_i^2$ of variables 
but not cross-products $x_ix_j$ with $i\neq j$ (for
example, Boolean constraints expressed as $x_i(x_i-1) = 0$).
Sparse matrix cones with block-arrow structure are often highlighted as 
an important example of chordal structure \cite{VaA:15,ZFP:21}.  As we have
seen, their homogenous cone property actually distinguishes them from 
general chordal sparse matrix cones.
The matrix norm cones described in Section~\ref{s-hom-cones-examples}
have also been studied separately, for their important role in 
optimization problems involving the matrix trace norm and spectral norm.
Except for these special cases, 
homogeneous matrix cones have been largely unexplored in modeling
convex optimization problems and in the development of scalable algorithms.

Our approach in this paper builds on fundamental results from various 
disciplines: abstract algebra, graph theory, sparse matrix computation 
and theory, convex conic optimization. An exciting next step
in research is the development of specialized algorithms and software 
which exploit the special structures we exposed here. There are many 
other interesting directions to be explored. Hyperbolicity cones are 
the next class of well-known convex cones which contain homogeneous cones 
as a strict subset. It may be fruitful to find a class of convex cones 
strictly between homogeneous cones and hyperbolicity cones providing a 
common generalization of homogeneous cones and cones of symmetric
positive semidefinite matrices with chordal sparsity.

\section*{Acknowledgments}
Research of the first author was supported in part by Discovery Grants 
from NSERC and by U.S.\ Office of Naval Research under award numbers
N00014-15-1-2171 and N00014-18-1-2078.
This work was started while the authors were visiting the Simons 
Institute for the Theory of Computing, supported in part by the 
DIMACS/Simons Collaboration on Bridging Continuous and Discrete 
Optimization through NSF grant \#CCF-1740425.

\appendix

\section{Background on homogeneous chordal graphs}
\label{s-app-graphs}
This appendix contains additional details for Section~\ref{s-graphs}.
We first review some results by Wolk \cite{Wol:62,Wol:65} and 
Golumbic \cite{Gol:78}, and then discuss the LBFS algorithm 
for recognizing and reordering homogeneous chordal graphs
\cite{Chu:08}.
We will use the term \emph{D-graph} when discussing Wolk's results 
in the next section, and the term \emph{trivially perfect graph}
when discussing Golumbic's results.  After that we use the term
\emph{homogeneous chordal graph} as in the rest of the paper.

\subsection{D-graphs}
Wolk \cite{Wol:62,Wol:65} defines a \emph{D-graph} or
\emph{graph with the diagonal property} as an undirected graph
that does not contain $P_4$ or $C_4$ as an induced subgraph.
He shows that this property characterizes the comparability graphs of
rooted forests.  

\begin{figure}
\hspace*{\fill}
\hspace*{\fill}
\begin{tikzpicture}[scale=1.5,baseline]
   \centering
   \tikzset{Vertex/.style = {shape = circle, draw, 
       minimum size = 12pt, inner sep = 1, fill=none}}
   \foreach \l / \k / \x / \y in {
       $u$/1/0.0/ 1.0,
       $x$/2/1.0/ 1.0,
       $w$/3/1.0/ 0.0,
       $v$/4/0.0/ 0.0}
       \node[Vertex, font=\small](\k) at (\x,\y){\l};
   \draw[>=Latex,->] (1) --(4);
   \draw[>=Latex,->] (3) --(4);
   \draw[>=Latex,->] (3) --(2);
   \draw[dashed] (1)--(2);
\end{tikzpicture}
\hspace*{\fill}
\begin{tikzpicture}[scale=1.5,baseline]
   \centering
   \tikzset{Vertex/.style = {shape = circle, draw, 
       minimum size = 12pt, inner sep = 1, fill=none}}
   \foreach \l / \k / \x / \y in {
       $u$/1/0.0/ 1.0,
       $x$/2/1.0/ 1.0,
       $w$/3/1.0/ 0.0,
       $v$/4/0.0/ 0.0}
       \node[Vertex, font=\small](\k) at (\x,\y){\l};
   \draw[>=Latex,<-] (1) --(4);
   \draw[>=Latex,<-] (3) --(4);
   \draw[>=Latex,<-] (3) --(2);
   \draw[dashed] (1)--(2);
\end{tikzpicture}
\hspace*{\fill}
\hspace*{\fill}
\caption{The vertices $u$, $v$, $w$, $x$ induce $C_4$ (if $u$ and $x$
are adjacent) or $P_4$ (otherwise). 
On the left, we assume $v>u$.  Transitivity of the partial ordering
and the fact that $\{u,w\}\not\in E$ and $\{v,x\} \not\in E$
imply that $v<w$ and $w<x$.  This ordering is incompatible with 
a tree structure because the vertex $w$  has two ancestors $v$ and $x$
that do not form an ancestor--descendant pair in the tree.
On the right we assume $u>v$.  Here, transitivity of the partial ordering
implies that $v<w$ and $x<w$.  Now $v$ has two ancestors $u$ and $w$ 
that do not form an ancestor--descendant pair.  }
\label{f-wolk}
\end{figure}
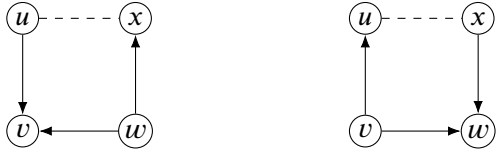

It is easy to show by contradiction that the absence 
of induced subgraphs $C_4$ or $P_4$ is a necessary condition for 
a graph $G = (V,E)$ to be the comparability graph of a rooted forest.
Suppose the vertices $u$, $v$, $w$, $x$ induce $C_4$ or $P_4$
(Figure~\ref{f-wolk}) and that there exists a rooted forest
with $G$ as its comparability graph.
We use the notation $a<b$ to denote that $b$ is an ancestor of $a$ in 
the forest ($b$ is on the unique path from $a$ to a root of the forest).
This defines a partial ordering: if $a<b$ and $b<c$, then $a<c$.   
There are two possible orientations for the edge $\{u,v\}$ in
Figure~\ref{f-wolk}, and for each orientation, there is only one 
possible orientation of the edges $\{v,w\}$ and $\{w,x\}$ that is 
compatible with the fact that $\{u,w\}\not \in E$ and $\{v,x\}\not\in E$.
For example, if $u<v$ as in the graph on the left, 
then necessarily $v>w$, because $v<w$ 
would imply that $u<w$  and therefore $\{u,w\}\in E$.
Now the two orientations in the figure are incompatible with
a tree structure because in each case, we find a vertex
($w$ in the graph on the left and $v$ on the right)
with two ancestors that are not mutually comparable (do not form
an ancestor--descendant pair in the tree).

For the second part of Wolk's result 
(every D-graph is the comparability graph of a rooted forest), 
we refer to Section~\ref{s-etree}, where we discuss how to
construct a rooted forest with comparability graph $G$.

Wolk also established the important 
property that every connected component of a D-graph has a 
\emph{universal vertex}, i.e., a vertex adjacent to all other vertices
\cite[lemma]{Wol:62}.
This can be seen as follows.  Without loss of generality we assume
that $G$ is connected. 
Let $v$ be the vertex with highest degree,
and denote its neigborhood by $\adj(v) = \{u_1, \ldots, u_k\}$ where $k$
is the degree of $v$.  We need to show that $v$ is a universal
vertex, i.e.,  $k=|V|-1$.
Assume that $k<|V|-1$.   Since the graph is connected, there exists
a vertex $w$ adjacent to one of the vertices $u_i$ and not adjacent
to $v$.  Thus,
$\{w,u_i\}\in E$, $\{u_i, v\} \in E$, and $\{w, v\}\not \in E$.
Consider any vertex $u_j$, $j\neq i$.  Since $\{v,u_j\} \in E$,
the vertices $v$, $u_i$, $w$, $u_j$ induce a $P_4$ or $C_4$
unless $u_i$ and $u_j$ are adjacent.  Therefore, if the graph is a D-graph,
$u_i$ must be adjacent to all $u_j$, $j\neq i$.  
However it is also adjacent to $v$ and to $w$,
so its degree is higher than the degree of $v$.  This contradicts our
assumption that $v$ is a vertex with maximum degree.

It was mentioned on page~\pageref{p-recursive-D-graphs} that this 
property leads to useful recursive characterization of D-graphs.
One consequence of this characterization is that D-graphs
are \emph{interval graphs} \cite{YCC:96}.
(In an interval graph the vertices represent intervals in $\R$;
two vertices are adjacent if and only if the corresponding intervals 
intersect.)
This follows from the construction method above, since clearly a disjoint
union of interval graphs is an interval graph, and the addition of
a universal vertex to an interval graph results in an interval graph.
The interval graphs are a subclass of the chordal graphs
\cite[chapter 8]{Gol:04}.

\subsection{Trivally perfect graphs} 
Golumbic \cite{Gol:78} defines a graph $G=(V,E)$ to be \emph{trivially 
perfect} if $\alpha(G_W) = m(G_W)$ holds for all $W\subseteq V$,
where $G_W$ denotes the subgraph induced by $W$, $\alpha(G_W)$ is
the stability number, and $m(G_W)$ the number of maximal cliques.
To motivate the name, recall that a graph is \emph{perfect} if
$\alpha(G_W) = \bar{\chi}(G_W)$ for all $W$, where $\bar{\chi}(G_W)$
is the clique cover number of $G_W$.   
Clearly, $\bar{\chi}(G_W) \leq m(G_W)$, so $\alpha(G_W) = m(G_W)$ 
immediately implies that $\alpha(G_W) = \bar{\chi}(G_W)$.

Golumbic gives the following simple proof to show that the 
trivially perfect graphs are exactly the graphs that do not contain
$C_4$ or $P_4$ as induced subgraphs.
First, we note that $\alpha(C_4) =2 < m(C_4) = 4$
and $\alpha(P_4) = 2 < m(P_4) = 3$, so a trivially perfect graph
cannot contain $C_4$ or $P_4$.
To show that the condition is sufficient, assume that $G$ does
not contain $C_4$ or $P_4$ as induced subgraphs.
Suppose $\alpha(G_W) < m(G_W)$ for some $W\subseteq V$.
Let $S$ be a maximum stable set of $G_W$.
Since $|S| = \alpha(G_W) < m(G_W)$, there exists a vertex
$s\in S$ that belongs to two different maximal cliques of $G_W$,
so we can find $x,y\in W$ with $\{s,x\}\in E$, $\{s,y\} \in E$,
$\{x,y\}\not\in E$.
Let $u$ be any element of $S\setminus \{s\}$ 
(note that $|S| = \alpha(G_W) \geq 2$ since $\{x,y\}\not\in E$). 
Therefore $\{s,u\}\not\in E$.  
If $\{x,u\}\in E$ and $\{u,y\}\in E$,
then the vertices $s,x,u,y$ induce a subgraph $C_4$.
If $\{x,u\}\in E$ and $\{u,y\}\not\in E$,
or $\{x,u\}\not\in E$ and $\{u,y\}\in E$,
then they induce a subgraph $P_4$.
We conclude that $\{x,u\}\not\in E$ and $\{u,y\}\not\in E$ for all 
$u\in S\setminus\{s\}$.
However this means that the set $(S\setminus \{s\}) \cup \{x,y\}$ is a 
stable set larger than $S$, 
contradicting the assumption that $S$ is a maximum stable set.

\subsection{Lexicographic Breadth First Search}
\label{s-lbfs}
We now discuss Chu's algorithm \cite{Chu:08} for recognizing
homogeneous chordal graphs and constructing a trivially perfect
elimination ordering $\sigma: \{1,2,\ldots, |V|\} \rightarrow V$.
The algorithm can be interpreted as reversing the recursive
construction of a homogeneous chordal graph via the operations
of disjoint union and addition of a universal vertex.
We number the vertices in the order $|V|, \ldots, 1$, i.e., select
$\sigma(|V|)$, \ldots, $\sigma(1)$ in that order.
At each step we find a universal vertex, give it the next available
number, and remove it from the graph.
Note that a universal vertex in a homogeneous chordal graph is easily 
found as a vertex with highest degree.   
 
Chu's algorithm  maintains a list  $L = (V_1, \ldots, V_K)$ 
of nonempty disjoint subsets of $V$.
The vertices in each set $V_i$ are ordered by nondecreasing 
degree (in $G$).

\begin{itemize}
\item 
Define $K=1$ and $L = (V_1)$,
with $V_1$ containing the elements of $V$ sorted in 
order of nondecreasing degree.  

\item For $i=|V|,\ldots,1$:  
\begin{enumerate}[label = \arabic{enumi}.]
\item Let $v$ be the last vertex in $V_K$.  
 Define $\sigma(i) = v$.
\item If $\adj(v) \cap V_j \neq \emptyset$ for some $j<K$,
terminate.  The graph is not a homogeneous chordal graph.
\item Otherwise, partition $V_K\setminus\{v\}$ in two sets 
\[
W' = V_K \cap \adj(v),\qquad 
W = (V_K\setminus \{v\}) \setminus W'.   
\]
The vertices in $W$ and $W'$ are kept in the order of nondecreasing 
degree (in $G$).
Replace the list $L$ by
\begin{equation} \label{e-L}
 L := (V_1, \ldots, V_{K-1}, W, W').
\end{equation}
If $W$ or $W'$ is empty, remove the empty sets from $L$.
Set $K$ equal to the length of the new list $L$.
\end{enumerate}
\end{itemize}
The complexity of the algorithm is $O(|E| + |V|)$. 

\begin{figure}
\hspace*{\fill}
\begin{minipage}{.4\linewidth}
\centering
\begin{tikzpicture}[scale=0.40]
   \tikzset{Bullet/.style = {
        shape = circle, minimum size = 3pt, inner sep = 0pt, fill=black, 
        draw=black, thick}}
   \foreach \i/\j/\k/\l in { 
        1/  5/ 2/12, 
        1/  4/ 2/ 6, 
        1/ 12/ 2/ 8, 
        2/  3/ 1/ 9, 
        2/  5/ 1/12,
        2/  4/ 1/ 6, 
        2/ 12/ 1/ 8, 
        3/  5/ 9/12,
        3/  4/ 9/ 6, 
        3/ 12/ 9/ 8, 
        5/  4/12/ 6, 
        5/ 12/12/ 8, 
        4/ 12/ 6/ 8, 
        6/ 11/11/ 5,  
        6/ 12/11/ 8, 
        7/  8/ 4/ 7, 
        7/ 11/ 4/ 5, 
        7/ 12/ 4/ 8, 
        8/ 11/ 7/ 5,
        8/ 12/ 7/ 8,
        9/ 10/10/ 3,
        9/ 11/10/ 5, 
        9/ 12/10/ 8, 
       10/ 11/ 3/ 5, 
       10/ 12/ 3/ 8, 
       11/ 12/ 5/ 8
   } { 
       \node[Bullet] at (\l-1,-\k+1){}; 
       \node[Bullet] at (\k-1,-\l+1){}; 
   }
   \foreach \i/\k in {
       1/2, 2/1, 3/9, 5/12, 4/6, 6/11, 7/4, 8/7, 9/10, 10/3, 11/5, 12/8} 
       \node[font=\small] at (\i-1,-\i+1) {$\i$};  
   \draw[thick] (-0.5,0.5) rectangle (11.5,-11.5);
\end{tikzpicture}
\end{minipage}
\hspace*{\fill}
\begin{minipage}{.3\linewidth}
\begin{center}
\begin{tabular}{cc} \toprule
Vertex & Degree\\ \midrule
1 & 4 \\ 
2 & 3 \\
3 & 3 \\
4 & 3 \\
5 & 6 \\
6 & 5 \\
7 & 3 \\
8 & 11 \\
9 & 4 \\ 
10 & 3 \\
11 & 2 \\ 
12 & 5 \\ \bottomrule 
\end{tabular}
\end{center}
\end{minipage}
\hspace*{\fill}
\caption{Undirected graph with vertex set $V=\{1,\ldots, 12\}$
and edges indicated by dots.  The table lists the degrees
of the 12 vertices.} 
\label{f-example-lbfs}
\end{figure}
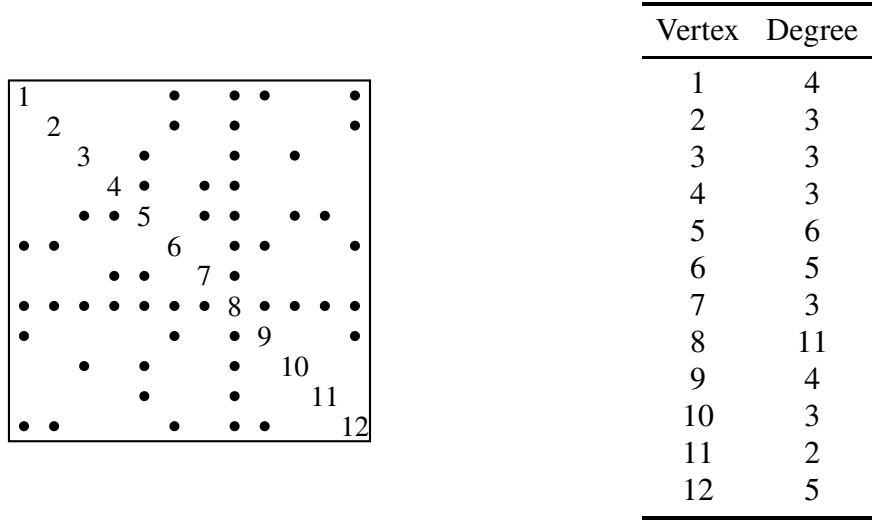

As an example we apply the algorithm to the graph
in Figure~\ref{f-example-lbfs}.
The sequence of partitions $L$ is shown in Figure~\ref{f-example-L}.
\begin{figure}
\begin{center}
\begin{tabular}{c@{\hskip 1em}l} \toprule 
 $i$ & \multicolumn{1}{c}{$L$} \\ \midrule
$12$ & 
\framebox{11, 2, 3, 4, 7, 10, 1, 9, 6, 12, 5, 8}  
\\*[1ex]
$11$ & 
\framebox{11, 2, 3, 4, 7, 10, 1, 9, 6, 12, 5} 
\\*[1ex]
$10$ & 
\framebox{2, 1, 9, 6, 12} \framebox{11, 3, 4, 7, 10} 
 \\*[1ex]
$9$ & 
\framebox{2, 1, 9, 6, 12} \framebox{11, 4, 7} \framebox{3} 
 \\*[1ex]
$8$ & 
\framebox{2, 1, 9, 6, 12} \framebox{11, 4, 7} 
\\*[1ex]
$7$ & 
\framebox{2, 1, 9, 6, 12} \framebox{11} \framebox{4}  
\\*[1ex]
$6$ & 
\framebox{2, 1, 9, 6, 12} \framebox{11} 
\\*[1ex]
$5$ & 
\framebox{2, 1, 9, 6, 12} 
\\*[1ex]
$4$ & 
\framebox{2, 1, 9, 6} 
\\*[1ex]
$3$ & 
\framebox{2, 1, 9} 
\\*[1ex]
$2$ & 
\framebox{2} \framebox{1} 
\\*[1ex]
$1$ & 
\framebox{2}   \\ \bottomrule
\end{tabular}
\end{center}
\caption{The partition $L$ at the start of each cycle
in the LBFS algorithm.
We start from the vertices $V$ sorted by degree.} 
\label{f-example-L}
\end{figure}
The ordering found by the algorithm is
\begin{subequations}\label{e-LBFS-sigma}
\begin{eqnarray} 
(\sigma(1), \ldots, \sigma(12))
& = & (2, 1, 9, 6, 12, 11, 4, 7, 3, 10, 5, 8) \\
(\sigma^{-1}(1), \ldots, \sigma^{-1}(12))
& = & (2, 1, 9, 7, 11, 4, 8, 12, 3, 10, 6, 5). 
\end{eqnarray}
\end{subequations}

We now verify that the algorithm recognizes homogeneous chordal graphs
\cite[theorem 3]{Chu:08}.  
First, assume that $G$ is a homogeneous chordal graph.
Let $L=(V_1,\ldots,V_K)$ be the partition at the start of a cycle
in the for-loop.
Assume that each set $V_j$ induces a homogeneous chordal subgraph 
$G_{V_j}$, disconnected from the other induced subgraphs $G_{V_k}$,
$k\neq j$.  
If $G$ is a homogeneous chordal graph, this assumption holds at the start 
of the algorithm.
Since $v\in V_K$, we have $\adj(v)\cap V_j =\emptyset$ for $j<K$,
so the algorithm does not terminate  in step 2.
Since $G_{V_K}$ is a homogeneous chordal graph, the sets $W$ and $W'$, 
which are subsets of $V_K$, also induce homogeneous chordal graphs.
Moreover, $v$ is a vertex with maximum degree in $V_K$, 
and therefore a universal vertex in the connected component of $G_{V_K}$ 
to which it belongs.  
This implies that $G_W$ is disconnected from $G_{W'}$.
We conclude that the sets in the new partition computed in step 3
of the algorithm define homogeneous chordal subgraphs that are mutually 
disconnected.
Therefore the algorithm completes the for-loop and does not terminate 
early.

Chu also shows that when the algorithm terminates early in step 2, 
a subgraph $P_4$ or $C_4$ that certifies that the graph is not a 
homogeneous chordal graphs, is easily obtained \cite[lemma 4]{Chu:08}.

Next we show that if the algorithm terminates successfully, the graph
$G$ is a homogeneous chordal graph.
Let $(V_1, \ldots, V_{K-1}, W, W')$ be the partition~(\ref{e-L}) 
at the end of cycle $i$ in the for-loop (with $W$ and $W'$ possibly
empty).
Assume that each set in this partition defines a homogeneous chordal 
graph, disconnected from the graphs induced by the other sets.
This is certainly true for $i=2$, since $L = (V_1, W, W')$
with $V_1 = \{\sigma(1)\}$ and $W=W'=\emptyset$. 
The set $V_K$ at the beginning of cycle  $i$
can be constructed by first adding a universal vertex $v$ to $W'$ and 
then making the disjoint union with the graph induced by $W$.  
Therefore $V_K$ induces a homogeneous chordal graph, disconnected
from the graphs induced by $V_1$, \ldots, $V_{K-1}$.
We conclude that the sets in the partition $L =(V_1,\ldots,V_K)$ at the 
beginning of cycle $i$ induce mutually disconnected homogeneous chordal 
subgraphs.  Therefore if the algorithm terminates the for-loop, the 
initial graph $G=(V,E)$ is a homogeneous chordal graph.

\subsection{Elimination tree} \label{s-etree}
We now discuss the ordering $\sigma$ produced by LBFS.
We use the notation 
\begin{eqnarray*}
 \madj(v) & = & \{w \in \adj(v) \,\, : \,\, \sigma^{-1}(w) > \sigma^{-1}(v)\} \\
 \ladj(v) & = & \{w \in \adj(v) \,\, : \,\, \sigma^{-1}(w) < \sigma^{-1}(v)\}
\end{eqnarray*}
for the \emph{higher} and \emph{lower neigborhoods} of $v$.   
We also define 
\[
 p(v) = \arg\min{\{\sigma^{-1}(w) \,\, : \,\, w \in\madj(v)\}}
\]
with the convention that $p(v) =v$ if $\madj(v)$ is empty.
The graph with vertex set $V$ and edges $\{v,p(v)\}$ for 
$p(v) \neq v$ is acyclic, since, by definition, 
$\sigma^{-1}(p(v)) > \sigma^{-1}(v)$.  
It is a rooted forest if we take the vertices with $p(v) = v$ as its roots.
The vertex $p(v)$ is the \emph{parent} of $v$ in the rooted forest.

\begin{figure}
\hspace*{\fill}
\begin{minipage}{.5\linewidth}
\centering
\begin{tikzpicture}[scale=0.4]
   \tikzset{Bullet/.style = {
       shape = circle, minimum size = 3pt, inner sep = 0pt, fill=black, 
       draw=black, thick}}
   \foreach \i/\j/\k/\l in { 
        1/  5/ 2/12, 
        1/  4/ 2/ 6, 
        1/ 12/ 2/ 8, 
        2/  3/ 1/ 9, 
        2/  5/ 1/12,
        2/  4/ 1/ 6, 
        2/ 12/ 1/ 8, 
        3/  5/ 9/12,
        3/  4/ 9/ 6, 
        3/ 12/ 9/ 8, 
        5/  4/12/ 6, 
        5/ 12/12/ 8, 
        4/ 12/ 6/ 8, 
        6/ 11/11/ 5,  
        6/ 12/11/ 8, 
        7/  8/ 4/ 7, 
        7/ 11/ 4/ 5, 
        7/ 12/ 4/ 8, 
        8/ 11/ 7/ 5,
        8/ 12/ 7/ 8,
       10/  9/10/ 3,
       10/ 11/10/ 5, 
       10/ 12/10/ 8, 
        9/ 11/ 3/ 5, 
        9/ 12/ 3/ 8, 
       11/ 12/ 5/ 8
   } { 
       \node[Bullet] at (\j-1,-\i+1){}; 
       \node[Bullet] at (\i-1,-\j+1){}; 
   }
   \foreach \i\k in {
       1/2, 2/1, 3/9, 5/12, 4/6, 6/11, 7/4, 8/7, 9/3, 10/10, 11/5, 12/8} 
       \node[font=\small] at (\i-1,-\i+1) {$\k$};  
   \draw[thick] (-0.5,0.5) rectangle (11.5,-11.5);
\end{tikzpicture}
\end{minipage}
\hspace*{\fill}
\begin{minipage}{.4\linewidth}
\centering
\begin{tikzpicture}[xscale=1.7, yscale=1.0, label distance = -.2em]
   \tikzset{VertexStyle/.style = {shape = circle, draw, 
       minimum size = 13pt, inner sep = 1, fill=none}}
   \foreach \v/\k/\x/\y in {
       2/1/-1.2/1,
       1/2/-0.4/0,
       9/3/-0.4/1,
       6/4/-0.8/2,
      12/5/-0.8/3,
      11/6/ 0.2/2,
       4/7/ 0.8/1,
       7/8/ 0.8/2,
       3/9/ 1.4/1,
      10/10/ 1.4/2, 
       5/11/ 0.8/3,
       8/12/ 0  /4}
      \node[VertexStyle, font=\small](\v) at (\x,\y)
          [label = right: \scriptsize $\k$]
          {\small $\v$};
   \foreach \i/\j in {2/6, 1/9, 9/6, 6/12, 12/8, 5/8, 11/5, 7/5, 4/7,
      10/5, 3/10} \draw (\i)--(\j);
\end{tikzpicture}
\end{minipage}
\hspace*{\fill}
\caption{The graph of Figure~\ref{f-example-lbfs} ordered using the
ordering~(\ref{e-LBFS-sigma}) and  the corresponding elimination tree.
The number next to node $v$ in the elimination tree is $\sigma^{-1}(v)$.} 
\label{f-reordered}
\end{figure}
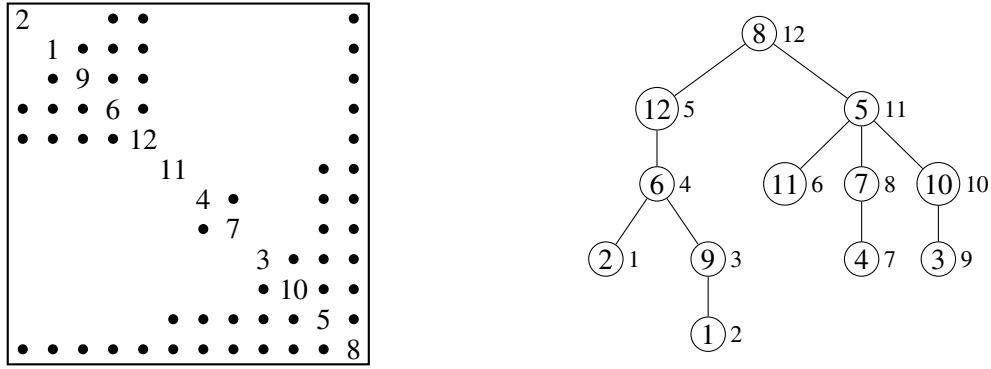

Figure~\ref{f-reordered} illustrates these definitions for the example.
In the array representation of the ordered graph, vertex $v$ appears on 
the diagonal of the array in position $\sigma^{-1}(v)$.  
The elements of $\madj(v)$ are found as the nonzeros below the diagonal 
in column $\sigma^{-1}(v)$.
The elements of $\ladj(v)$ are the elements to the left of the
diagonal in row $\sigma^{-1}(v)$.
The parent $p(v)$ of $v$ is the first nonzero below the diagonal.

The parent function can be computed by modifying the LBFS algorithm
as follows \cite[p.11]{Chu:08}.
Let $v = \sigma(i)$ be the vertex selected in step~1 of cycle 
$i$ of the algorithm.  This is called the \emph{pivot} \cite{Chu:08}.
The set $W'$ in step~3 is the lower neighborhood $\ladj(v)$, since
it contains the vertices adjacent to $v$ that will be numbered after $v$.
Since $w\in\ladj(v)$ if and only if $v \in \madj(w)$,
we find the parent $p(w)$ as the last pivot $v$ before 
$w$ is numbered for which $w\in\ladj(v)$.
To construct the parent function, we initialize $p(v) = v$ for all 
$v\in V$ at the start of the algorithm.
In step~3 of the algorithm we set $p(w) = v$ for all $w\in W'$.
Figure~\ref{f-example-tree} shows the value of $p(v)$ at the end
of each LBFS cycle in the example.
\begin{figure}
\hspace*{\fill}
\begin{minipage}{.5\linewidth}
\centering
\begin{tabular}{
c@{\hskip 1em}
c@{\hskip .5em}
c@{\hskip .5em}
c@{\hskip .5em}
c@{\hskip .5em}
c@{\hskip .5em}
c@{\hskip .5em}
c@{\hskip .5em}
c@{\hskip .5em}
c@{\hskip .5em}
c@{\hskip .5em}
c@{\hskip .5em}
c} \toprule 
& \multicolumn{12}{c}{Vertex} \\ \cmidrule{2-13} 
$i$ & 
$1$ & $2$ & $3$ & $4$ & $5$ & $6$ & $7$ & $8$ & $9$ & $10$ & $11$ & $12$ 
\\ \midrule
$12$ & 
$8$ & $8$ & $8$ & $8$ & $8$ & $8$ & $8$ & $8$ & $8$ & $8$ & $8$ & $8$ 
\\ 
$11$ & 
$8$ & $8$ & $5$ & $5$ & $8$ & $8$ & $5$ & $8$ & $8$ & $5$ & $5$ & $8$ 
\\ 
$10$ & 
$8$ & $8$ & $10$ & $5$ & $8$ & $8$ & $5$ & $8$ & $8$ & $5$ & $5$ & $8$ 
\\ 
$9$ & 
$8$ & $8$ & $10$ & $5$ & $8$ & $8$ & $5$ & $8$ & $8$ & $5$ & $5$ & $8$ 
\\
$8$ & 
$8$ & $8$ & $10$ & $7$ & $8$ & $8$ & $5$ & $8$ & $8$ & $5$ & $5$ & $8$ 
\\
$7$ & 
$8$ & $8$ & $10$ & $7$ & $8$ & $8$ & $5$ & $8$ & $8$ & $5$ & $5$ & $8$ 
\\
$6$ & 
$8$ & $8$ & $10$ & $7$ & $8$ & $8$ & $5$ & $8$ & $8$ & $5$ & $5$ & $8$ 
\\
$5$ & 
$12$ & $12$ & $10$ & $7$ & $8$ & $12$ & $5$ & $8$ & $12$ & $5$ & $5$ & $8$ 
\\
$4$ & 
$6$ & $6$ & $10$ & $7$ & $8$ & $12$ & $5$ & $8$ & $6$ & $5$ & $5$ & $8$ 
\\
$3$ & 
$9$ & $6$ & $10$ & $7$ & $8$ & $12$ & $5$ & $8$ & $6$ & $5$ & $5$ & $8$ 
\\
$2$ & 
$9$ & $6$ & $10$ & $7$ & $8$ & $12$ & $5$ & $8$ & $6$ & $5$ & $5$ & $8$ 
\\
$1$& 
$9$ & $6$ & $10$ & $7$ & $8$ & $12$ & $5$ & $8$ & $6$ & $5$ & $5$ & $8$ 
\\
\bottomrule
\end{tabular}
\end{minipage}
\hspace*{\fill}
\caption{Parent function $p(w)$ at the end of cycle $i=12,\ldots,1$.}
\label{f-example-tree}
\end{figure}

Assume now, without loss of generality, that the graph $G$ is connected,
so the rooted forest defined by the parent function $p(v)$ is 
a tree, called the \emph{elimination tree}.
Consider the partition~(\ref{e-L}) 
in cycle $i = \sigma^{-1}(v)$, when $v$ is the pivot. 
A vertex $w\in W' = \ladj(v)$ receives $p(w) = v$. 
This vertex is not adjacent to any of the elements of
$V_1$, \ldots, $V_{K-1}$, $W$.  If in subsequent cycles, the value
of $p(w)$ is updated, the new value can only be another element 
in $\ladj(v)$.
It follows that the vertices in $\ladj(v)$ form the subtree in
the elimination tree with root $v$.
Moreover, by definition of $W' = \ladj(v)$, the vertex $v$ is adjacent
to every element in $\ladj(v)$, i.e., all the descendants of $v$ in
the elimination tree.   
Equivalently, every vertex $w$
is adjacent to all its ancestors in the elimination tree (all vertices
on the unique path between $w$ and the root).
Finally, if two vertices $w$, $z$ do not form an ancestor--descendant
pair, they were placed in different sets of the partition
when their least common ancestor was the pivot.  Therefore $w$ and
$z$ are not adjacent.
In summary, two vertices in $G$ are adjacent if and only if 
they are comparable (form an ancestor--descendant pair)
in the elimination tree.
In other words, $G$ is the comparability graph of the elimination tree.
It also follows that $\sigma$ is a \emph{trivially perfect elimination 
ordering}, i.e., $\madj(v)$ induces a complete subgraph of $G$ 
and $\madj(v)$ contains the vertices on the path from $v$
to the root.

Finally, we note that placing $W'$ last in the list~(\ref{e-L}) 
ensures that the computed ordering is a \emph{postordering,} 
i.e., if $\sigma^{-1}(v) = i$ and $v$ has $k$
descendants in the elimination tree, then the descendants $w$ will have
consecutive positions $\sigma^{-1}(w) = i-k, \ldots, i-1$ in the ordering.

\section{Matrix algorithms for homogeneous 
chordal sparsity} \label{s-lin-alg}
In this appendix we outline algorithms for the basic matrix operations
discussed in Sections~\ref{s-cone} and~\ref{s-barriers}.
The algorithms are similar to the multifrontal algorithms for
matrices with chordal sparsity patterns described 
in \cite{AndersonDahlVandenberghe2013,VaA:15}, with additional
simplifications to exploit homogeneous chordal sparsity.

We consider a sparsity pattern described by a homogeneous
chordal graph $G = (V,E)$ with $V=\{1,2,\ldots,N\}$ and  assume the 
numerical order $1,2,\ldots,N$ is a trivially perfect 
elimination ordering of $V$. 
We denote by $\alpha_i$ the set of row indices of the lower-triangular
nonzeros in column~$i$, and by 
$\bar\alpha_i$ the set $\{i\} \cup \alpha_i$.
The parent of a non-root vertex $i$ in the elimination tree is denoted
by $p(i)$.  By definition, this is the first element of $\alpha_i$.

If $1,\ldots, N$ is a perfect elimination ordering of a chordal pattern, 
we have the important property 
\begin{equation} \label{e-col-chordal}
 \alpha_i \subseteq \bar\alpha_{p(i)} 
\end{equation}
for all non-root vertices $i$.
By applying this recursively, we see that the vertices indexed by 
$\alpha_i$ are on the path from vertex $i$ to the root, i.e.,
$\alpha_i \subseteq \{p(i), p^2(i), \ldots, p^k(i)\}$, if $k$ is the 
depth (distance to the root) of vertex $i$.
If $1,\ldots, N$ is a trivially perfect elimination ordering of a 
homogeneous chordal pattern, we have equality:
\begin{equation} \label{e-col-tp-app}
 \alpha_i = \bar\alpha_{p(i)}. 
\end{equation}
Therefore $\alpha_i = \{p(i), p^2(i),\ldots, p^k(i)\}$, the set of 
ancestors of vertex $i$ in the elimination tree.

The algorithms presented in the rest of this section use a recursion
on the elimination tree.  A recursion in \emph{topological order}
visits each node of the elimination tree before its parent.
A recursion in \emph{inverse topological order} visits each node before 
its children.
We also use the notation $\ch(i)$ for the set of children of node $i$
in the elimination tree.
Supernodal elimination trees can be used to 
formulate faster supernodal or blocked versions, but this extension
will not be discussed in detail.

\subsection{Cone automorphisms}
Our main interest in this section is the evaluation of the linear mappings 
$\cL$ and $\cL^*$ 
defined in~(\ref{e-mL}) and~(\ref{e-mL*}).  We first consider two simpler 
operations, matrix--matrix products  $L\tilde L$ and $L\tilde L^\top$,
where $L, \tilde L\in\TT^N_E$.

\subsubsection{Products of lower-triangular matrices}
Let $L, \tilde L \in\TT^N_E$.
Column $k$ of the product $Y=L\tilde L$ can be computed by 
initializing the column as zero, and running the iteration
\[
 \left[\begin{array}{c} Y_{jk} \\ Y_{\alpha_jk} \end{array}\right]
:= 
 \left[\begin{array}{c} Y_{jk} \\ Y_{\alpha_jk} \end{array}\right]
+ \tilde L_{jk}  
 \left[\begin{array}{c} L_{jj} \\ L_{\alpha_jj} \end{array}\right],
 \quad j\in \bar\alpha_k.
\]
Here we rely on the fact that the nonzero elements in column $k$ of 
$\tilde L$ are in the rows indexed by $\bar \alpha_k$, and the
nonzeros in column $j$ of $L$ are in the rows indexed by $\bar\alpha_j$.
Now, the property~(\ref{e-col-tp}) implies that for a trivially perfect
elimination ordering, $\alpha_j \subset \alpha_k$ for $j\in \alpha_k$.
Therefore the nonzeros in the $k$th column of $Y$ are in 
the rows indexed by $\bar\alpha_k$.
This again shows that $Y=L\tilde L\in\TT_E^N$, as already noted in 
Theorem~\ref{thm:1.1}.

Next we consider products $Y = L\tilde L^\top$, 
where $L, \tilde L \in \TT_E^N$. The matrix $Y$
is not symmetric, but has a symmetric sparsity pattern,
and if $E$ is chordal and $1,\ldots,N$ is a perfect elimination order,
then the sparsity pattern of $Y$ is $E$.
To see this consider the formula for the $ij$ element of $Y$:
\[
 Y_{ij}  = \sum_{k=1}^{\min\{i,j\}} L_{ik} \tilde L_{jk}.
\]
For a perfect elimination ordering of a chordal graph, 
$\{i,k\}\in E$, $\{j,k\}\in E$ implies that $\{i,j\}\in E$.  
So if $\sum_k L_{ik}\tilde L_{jk}$ is nonzero then $\{i,j\} \in E$.  

An efficient method for computing $Y$ can be formulated as a recursion
on the elimination tree, using ideas from the multifrontal
Cholesky factorization (see Section~\ref{s-chol}).  
As in the multifrontal method, we start from the equation
for the $\bar\alpha_i$-by-$ \bar\alpha_i$ block:
\begin{eqnarray} 
\lefteqn{\left[\begin{array}{cc}
 Y_{ii} & Y_{i\alpha_i} \\ Y_{\alpha_ii} & Y_{\alpha_i\alpha_i}
 \end{array}\right]} \nonumber \\
& = & 
 \left[\begin{array}{@{\hskip .05em}c@{\hskip .05em}} 
  L_{ii} \\ L_{\alpha_ii} \end{array}\right]
 \left[\begin{array}{@{\hskip .05em}c@{\hskip .05em}} 
  \tilde L_{ii} \\ \tilde L_{\alpha_ii} 
  \end{array}\right]^\top
+ \sum_{k<i} 
 \left[\begin{array}{@{\hskip .05em}c@{\hskip .05em}} 
   L_{ik} \\ L_{\alpha_ik} \end{array}\right]
 \left[\begin{array}{@{\hskip .05em}c@{\hskip .05em}} 
  \tilde L_{ik} \\ \tilde L_{\alpha_ik} 
  \end{array}\right]^\top
 +
 \sum_{k>i} 
 \left[\begin{array}{@{\hskip .05em}c@{\hskip .05em}} 
   0 \\ L_{\alpha_ik} \end{array}\right]
 \left[\begin{array}{@{\hskip .05em}c@{\hskip .05em}} 
   0  \\ \tilde L_{\alpha_ik} \end{array}\right]^\top.
\label{e-llt-prod}
\end{eqnarray}
We define for each vertex $j$ a nonsymmetric \emph{update matrix}
\begin{equation} \label{e-update-matrix}
 U_j = -\sum_{k\in T_j} L_{\alpha_jk} \tilde L_{\alpha_jk}^\top 
\end{equation}
where $T_j$ is the subtree of the elimination tree rooted at node $j$.
With this notation, the first column and row of 
equation~(\ref{e-llt-prod}), and the definition of the update matrix
$U_i$ using~(\ref{e-update-matrix}), can be combined in the equation
\[
\left[\begin{array}{cc}
 Y_{ii} & Y_{i\alpha_i} \\ Y_{\alpha_ii} & -U_i\end{array}\right]
 =  \left[\begin{array}{c} L_{ii} \\ L_{\alpha_ii} \end{array}\right]
  \left[\begin{array}{c} \tilde L_{ii} \\ \tilde L_{\alpha_ii} 
 \end{array}\right]^\top - \sum_{j\in \ch(i)} U_j,
\]
where $\ch(i)$ is the set of children of node $i$  in the elimination tree.
This recursion allows us to compute $Y$.  We enumerate the vertices
$i$ of the elimination tree in topological order (i.e., visiting each node
before its parent, for example, in the order $1,\ldots,N$).  
For each $i$ we compute $Y_{ii}$, $Y_{i\alpha_i}$, $Y_{\alpha_ii}$,
and $U_i$ from column $i$ of $L$ and $\tilde L$, and 
from the update matrices of the children of $i$.
After the update at vertex $i$ the matrices $U_j$ for $j\in \ch(i)$
can be discarded.

\subsubsection{Triangular scaling of symmetric matrix}
We now turn to the computation of $\cL(X)=LXL^\top$ where $X\in\SS^N_E$ 
and $L\in\TT^N_E$.  The operation can be reduced
to a combination of the previous cases by
splitting $X$ as $X = \tilde L + \tilde L^\top$, where $\tilde L$
is lower-triangular with nonzero elements
$\tilde L_{ii} = X_{ii}/2$, $\tilde L_{\alpha_ii} = X_{\alpha_ii}$ 
for $i=1,\ldots,n$.  The other entries of $\tilde L$ are zero.
Then $\cL(X)$ can be written as 
\[
 \cL(X) = L(\tilde L + \tilde L^\top) L^\top 
  = (L\tilde L)L^\top + L(L\tilde L)^\top.
\]
We first compute $\hat L = L\tilde L$ column by column using
\[
 \left[\begin{array}{c} \hat L_{ii} \\ \hat L_{\alpha_ii} 
 \end{array}\right]
 = \left[\begin{array}{cc} 
 L_{ii} & 0  \\ L_{\alpha_ii} & L_{\alpha_i\alpha_i} \end{array}\right]
 \left[\begin{array}{c} X_{ii}/2 \\ X_{\alpha_ii} \end{array}\right], 
 \quad i=1,\ldots,N.
\]
Then $Y = \cL(X)$ can be computed via 
\[
\left[\begin{array}{cc}
 Y_{ii} & Y_{i\alpha_i} \\ Y_{\alpha_ii} & -U_i\end{array}\right]
 =  \left[\begin{array}{c} \hat L_{ii} \\ \hat L_{\alpha_ii} 
    \end{array}\right]
  \left[\begin{array}{c} L_{ii} \\ L_{\alpha_ii} \end{array}\right]^\top 
 +  \left[\begin{array}{c} L_{ii} \\ L_{\alpha_ii} 
  \end{array}\right]
  \left[\begin{array}{c} \hat L_{ii} \\ \hat L_{\alpha_ii} 
 \end{array}\right]^\top - \sum_{j\in \ch(i)} U_j
\]
in topological order.
Combining the two steps gives the formula
\begin{eqnarray}
\lefteqn{
\left[\begin{array}{cc}
 Y_{ii} & Y_{i\alpha_i} \\ Y_{\alpha_ii} & -U_i\end{array}\right] }
\nonumber \\
 & = & \left[\begin{array}{cc} 
 L_{ii} & 0  \\ L_{\alpha_ii} & L_{\alpha_i\alpha_i} \end{array}\right]
 \left[\begin{array}{c} X_{ii}/2 \\ X_{\alpha_ii} \end{array}\right]
  \left[\begin{array}{c} L_{ii} \\ L_{\alpha_ii} \end{array}\right]^\top 
 \nonumber \\
& & \mbox{}
 +  \left[\begin{array}{c} L_{ii} \\ L_{\alpha_ii} 
  \end{array}\right]
  \left[\begin{array}{c} X_{ii}/2 \\ X_{\alpha_ii} \end{array}\right]^\top 
  \left[\begin{array}{cc} L_{ii} & L_{\alpha_ii}^\top \\   
  0 & L_{\alpha_i\alpha_i}^\top \end{array}\right] 
- \sum_{j\in \ch(i)} U_j  \nonumber \\
 & = & \left[\begin{array}{@{\hskip.1em}cc@{\hskip.1em}} 
 L_{ii} & 0  \\ L_{\alpha_ii} & L_{\alpha_i\alpha_i} \end{array}\right]
 \left[\begin{array}{@{\hskip.1em}cc@{\hskip.1em}} 
  X_{ii} & X_{\alpha_ii}^\top \\
   X_{\alpha_ii} & 0 \end{array}\right]
  \left[\begin{array}{@{\hskip.1em}cc@{\hskip .1em}} 
   L_{ii} & L_{\alpha_ii}^\top \\   
  0 & L_{\alpha_i\alpha_i}^\top \end{array}\right] 
- \sum_{j\in \ch(i)} U_j, \label{e-F-eq}
\end{eqnarray}
which can be evaluated by a recursion in topological order.
The algorithm is summarized as follows.

\begin{algdesc}{Forward mapping $\cL$} \label{a-F}
\begin{list}{}{}
\item[\textbf{Input.}]
A matrix $X \in \SS^N_E$ with homogeneous chordal 
sparsity pattern and trivially perfect elimination ordering 
$\sigma = (1,2,\ldots,N)$, a lower-triangular matrix $L\in\TT^N_E$, 
and the elimination tree for $\sigma$.
\item[\textbf{Output.}]  The matrix $Y = LXL^\top$.
\item[\textbf{Algorithm.}]
\begin{enumerate}[label = \arabic{enumi}.]
\item Define a lower-triangular matrix $W\in\TT^N_E$ with 
\[
 W_{ii} = X_{ii}, \qquad 
 W_{\alpha_ii} = L_{\alpha_i\alpha_i} X_{\alpha_ii}, \qquad
 i=1,\ldots,N.
\]
\item Enumerate the vertices $i=1, 2, \ldots,N$ of the
elimination tree in topological order.  
For each $i$, compute $U_i$, $Y_{ii}$, $Y_{\alpha_ii}$ using
the formula
\begin{eqnarray*}
\lefteqn{
\left[\begin{array}{cc}
 Y_{ii} & Y_{i\alpha_i} \\ Y_{\alpha_ii} & -U_i\end{array}\right]} \\
 & = & \left[\begin{array}{cc} 
 L_{ii} & 0  \\ L_{\alpha_ii} & I \end{array}\right]
 \left[\begin{array}{cc} W_{ii} & W_{\alpha_ii}^\top \\
   W_{\alpha_ii} & 0 \end{array}\right]
  \left[\begin{array}{cc} L_{ii} & L_{\alpha_ii}^\top \\   
  0 & I \end{array}\right] - \sum_{j\in \ch(i)} U_j.
\end{eqnarray*}
\end{enumerate}
\end{list}
\end{algdesc} 
The intermediate variable $W$ and the computation in step~1 are 
introduced to make the adjoint relation with the algorithm
for $\cL^*$ in the following paragraph clearer.

\subsubsection{Adjoint triangular scaling of symmetric matrix}
The next operation is $\cL^*(S) = \Pi_E(L^{\top} SL)$.
The $\bar\alpha_i $-by-$ \bar\alpha_i$ block of $Y=L^{\top}SL$ is
\begin{equation} \label{e-Fadj-eq}
\left[\begin{array}{cc}
 Y_{ii} & Y_{i\alpha_i} \\ Y_{\alpha_ii} & Y_{\alpha_i\alpha_i} 
 \end{array}\right]
 = \left[\begin{array}{cc} 
 L_{ii} & L_{\alpha_ii}^\top   \\ 0 & L_{\alpha_i\alpha_i}^\top
  \end{array}\right]
 \left[\begin{array}{cc} S_{ii} & S_{\alpha_ii}^\top \\
   S_{\alpha_ii} & S_{\alpha_i\alpha_i} \end{array}\right]
  \left[\begin{array}{cc} L_{ii} & 0 \\   
  L_{\alpha_ii} & L_{\alpha_i\alpha_i}  \end{array}\right]. 
\end{equation}
This follows from the fact that the block column of $L$ indexed
by $\alpha_i$ has no nonzeros outside the rows $\alpha_i$.
This is not true for a general chordal pattern. For a general
chordal pattern
the expression~(\ref{e-Fadj-eq}) gives the wrong value 
for the $22$ block $Y_{\alpha_i\alpha_i}$, although the expressions for 
$Y_{ii}$ and $Y_{\alpha_ii}$ are correct.
Even for a homogeneous chordal pattern, we actually do not use 
the $22$ block, since these elements are part of other columns 
and we need to compute them only once.
A possible implementation is as follows.
The intermediate variable $V_i$ in this algorithm is 
simply $S_{\alpha_i\alpha_i}$.   Passing this dense matrix from nodes
to their children is more efficient than retrieving $S_{\alpha_i\alpha_i}$
from a sparse matrix data structure \cite{AndersonDahlVandenberghe2013}.

\begin{algdesc}{Adjoint mapping $\cL^*$} \label{a-Fadj}
\begin{list}{}{}
\item[\textbf{Input.}]
A matrix $S \in \SS^N_E$ with a homogeneous chordal
sparsity pattern and trivially perfect elimination ordering 
$\sigma = (1,2,\ldots,N)$, a lower-triangular matrix
$L\in\TT^N_E$, and the elimination tree for $\sigma$.
\item[\textbf{Output.}]  The matrix $Y = \Pi_E(L^{\top}SL)$.
\item[\textbf{Algorithm.}]
\begin{enumerate}[label = \arabic{enumi}.]
\item Compute a lower-triangular matrix $W\in\TT^N_E$
by running the following recursion in reverse topological order.
For each $i$, compute $W_{ii}$ and $W_{i\alpha_i}$ from
\[
\left[\begin{array}{cc}
 W_{ii} & W_{i\alpha_i}^\top \\ W_{i\alpha_i} & \times \end{array}\right]
 = \left[\begin{array}{cc} 
 L_{ii} & L_{\alpha_ii}^\top   \\ 0 & I \end{array}\right]
 \left[\begin{array}{cc} S_{ii} & S_{\alpha_ii}^\top \\
   S_{\alpha_ii} & V_i \end{array}\right]
  \left[\begin{array}{cc} L_{ii} & 0 \\   
  L_{\alpha_ii} & I \end{array}\right], 
\]
and define
\[
 V_j = \left[\begin{array}{cc} S_{ii} & S_{\alpha_ii}^\top \\
   S_{\alpha_ii} & V_i \end{array}\right], \quad
 j\in\ch(i).
\]
\item For $i=1,\ldots,N$, set 
\[
 Y_{ii} = W_{ii}, \qquad
 Y_{\alpha_ii} = L_{\alpha_i\alpha_i}^\top W_{\alpha_ii}.
\]
\end{enumerate}
\end{list}
\end{algdesc}

\subsection{Inverse cone automorphisms} \label{s-prod-inv}
Next we consider the inverses of the mappings $\cL$ and $\cL^*$.
Again we start with some observations about simpler operations with the 
inverse of a sparse lower-triangular matrix.

\subsubsection{Products with inverse of lower-triangular matrix}
To solve $Lx=b$, we set $x :=b$ and run the iteration
\begin{equation} \label{e-forward-subs}
 \left[\begin{array}{c} x_j \\ x_{\alpha_j} \end{array}\right]
 := \left[\begin{array}{cc} 1/L_{jj} & 0 \\ -L_{\alpha_jj}/L_{jj}
   & I \end{array}\right] \left[\begin{array}{c} x_j \\ x_{\alpha_j}
 \end{array}\right], \qquad j=1,\ldots,N.
\end{equation}
The algorithm does not require chordality or homogeneous chordality,
but the order of the recursion matters.  
If the pattern is chordal and the right-hand side $b$ is sparse, 
we can simplify the iteration and iterate over a ``pruned'' elimination 
tree, defined by the vertices $k$ with $b_k\neq 0$ and their ancestors.
This follows from~(\ref{e-col-chordal}):
all the elements of $\bar\alpha_j$ are on the path
from vertex $j$ to the root, so the iteration~(\ref{e-forward-subs}) 
does not change entries outside this pruned elimination tree.
In particular, if $b$ has only one nonzero entry $b_k$, 
then in~(\ref{e-forward-subs}) we can iterate over the vertices 
$j = k, p(k), p^2(k), \ldots,$ on the path from 
$k$ to the root of the elimination tree. 

The product $X= L^{-1}\tilde L$ can be computed column by column, by
forward substitution.
Set $X= \tilde L$.  For each $k=1,\ldots,n$, run the iteration
\[
 \left[\begin{array}{c} X_{jk} \\ X_{\alpha_jk} \end{array}\right]
  := 
 \left[\begin{array}{cc}
   1/L_{jj} & 0 \\ -L_{\alpha_jj}/L_{jj} & I \end{array}\right]
 \left[\begin{array}{c} X_{jk} \\ X_{\alpha_jk} \end{array}\right],
 \qquad j = k, p(k), p^2(k), \ldots.
\]
This works for any chordal sparsity pattern.  In general, however, the sets
$\alpha_j$ for $j= p(k), p^2(k), \ldots$ are not subsets of 
$\alpha_k$, so the final sparsity pattern of $X_k$ can include nonzeros 
outside~$\alpha_k$.
For a homogeneous chordal pattern and trivially perfect elimination
ordering, the property~(\ref{e-col-tp})
implies that the indices of all lower-triangular nonzeros of $X_k$ are 
in $\alpha_k$.  Therefore $X=L^{-1}\tilde L$ has the same sparsity 
pattern as $L$ and $\tilde L$.

Applying this with $\tilde L=I$, we see that
the inverse $L^{-1}$ has the same sparsity pattern as $L$:
$L^{-1} \in \TT^N_E$;  see Theorem~\ref{thm:6.3}.  
This property 
does not hold for general chordal sparsity pattern.
As a consequence, the identities
\begin{equation} \label{e-Linv}
 (L^{-1})_{\bar\alpha_j\bar\alpha_j} = L_{\bar\alpha_j\bar\alpha_j}^{-1}
 = \left[\begin{array}{@{\hskip .1em}cc@{\hskip .1em}}
  1/L_{jj} & 0 \\ -(1/L_{jj})
 L_{\alpha_j\alpha_j}^{-1}L_{\alpha_jj} & 
   L_{\alpha_j\alpha_j}^{-1} \end{array}\right], \quad 
 j=1,\ldots,N,
\end{equation}
(which hold for any nonsingular triangular matrix and any index set 
$\bar\alpha_j$) characterize all the nonzero elements in $L^{-1}$.

\subsubsection{Inverse of triangular scaling}
The inverse of the mapping $\cL(X) = LXL^\top$ is $L^{-1}XL^{-\top}$. 
It can be evaluated via the formula~(\ref{e-F-eq}) applied to the 
inverse of $L$:
\[
\left[\begin{array}{cc}
 Y_{ii} & Y_{i\alpha_i} \\ Y_{\alpha_ii} & -U_i 
\end{array}\right]
 = 
L_{\bar\alpha_i\bar\alpha_i}^{-1}
\left[\begin{array}{cc}
   X_{ii} & X_{\alpha_ii}^\top \\ X_{\alpha_ii} & 0 \end{array}\right]
L_{\bar\alpha_i\bar\alpha_i}^{-\top}
 - \sum_{j\in\ch(i)} U_j.
\]
This can be simplified if we define update matrices 
$V_i = L_{\alpha_i\alpha_i} U_i L_{\alpha_i\alpha_i}^\top$ instead of~$U_i$:
\[
\left[\begin{array}{cc}
 Y_{ii} & Y_{i\alpha_i} \\ Y_{\alpha_ii} & 
 - L_{\alpha_i\alpha_i}^{-1}V_i L_{\alpha_i\alpha_i}^{-\top}
\end{array}\right]
 = L_{\bar\alpha_i\bar\alpha_i}^{-1} \left(
\left[\begin{array}{cc}
   X_{ii} & X_{\alpha_ii}^\top \\ X_{\alpha_ii} & 0 \end{array}\right]
 - \sum_{j\in\ch(i)} V_j \right) L_{\bar\alpha_i\bar\alpha_i}^{-\top},
\]
and, using~(\ref{e-Linv}),
\begin{eqnarray*}
\lefteqn{
\left[\begin{array}{cc}
 Y_{ii} & Y_{i\alpha_i}L_{\alpha_i\alpha_i}^\top \\ 
 L_{\alpha_i\alpha_i}Y_{\alpha_ii}  & - V_i 
\end{array}\right]} \\
 & = &
\left[\begin{array}{cc}
1/L_{ii} & 0 \\ - L_{\alpha_ii}/L_{ii} & I \end{array}\right]
\left( \left[\begin{array}{cc}
   X_{ii} & X_{\alpha_ii}^\top \\ X_{\alpha_ii} & 0 \end{array}\right]
 - \sum_{j\in\ch(i)} V_j \right) 
\left[\begin{array}{cc}
1/L_{ii} & -L_{\alpha_ii}^\top/L_{ii}  \\
 0 & I \end{array}\right].
\end{eqnarray*}
This is summarized in the following outline.

\begin{algdesc}{Inverse forward mapping $\cL^{-1}$} \label{a-F-inv}
\begin{list}{}{}
\item[\textbf{Input.}]
A matrix $X \in \SS^N_E$ with homogeneous chordal
sparsity pattern and trivially perfect elimination ordering 
$\sigma = (1,2,\ldots,N)$, a nonsingular
lower-triangular matrix $L\in\TT_E^N$, 
and the elimination tree for $\sigma$.
\item[\textbf{Output.}]  The matrix $Y = L^{-1}XL^{-\top}$.
\item[\textbf{Algorithm.}]
\begin{enumerate}[label = \arabic{enumi}.]
\item Enumerate the vertices $i=1, 2, \ldots,N$ of the
elimination tree in topological order.  
For each $i$, compute $V_i$, $W_{ii}$, $W_{\alpha_ii}$ using
the formula
\begin{eqnarray*}
\lefteqn{\left[\begin{array}{cc}
 W_{ii} & W_{i\alpha_i} \\ W_{\alpha_ii} & -V_i\end{array}\right]} \\
 &   = & 
 \left[\begin{array}{@{\hskip .1em}c@{\hskip .3em}c@{\hskip .1em}} 
 1/L_{ii} & 0  \\ -L_{\alpha_ii}/L_{ii} & I \end{array}\right]
 \left(
 \left[\begin{array}{@{\hskip .1em}c@{\hskip .3em}c@{\hskip .1em}} 
  X_{ii} & X_{\alpha_ii}^\top \\ X_{\alpha_ii} & 0 \end{array}\right]
  - \sum_{j\in \ch(i)} V_j \right) 
  \left[\begin{array}{@{\hskip .1em}c@{\hskip .3em}c@{\hskip .1em}} 
 1/L_{ii} & -L_{\alpha_ii}^\top/L_{ii} \\   
  0 & I \end{array}\right].
\end{eqnarray*}
\item For $i=1,\ldots,N$, compute
\[
 Y_{ii} = W_{ii}, \qquad 
 Y_{\alpha_ii} = L_{\alpha_i\alpha_i}^{-1} W_{\alpha_ii}. 
\]
\end{enumerate}
\end{list}
\end{algdesc} 

\subsubsection{Inverse of adjoint triangular scaling}
Applying~(\ref{e-Fadj-eq}) with $L^{-1}$ shows that
$\cL^{-*}(S) = \Pi_E(L^{-\top}SL^{-1})$ and
that the $\bar\alpha_i $-by-$ \bar\alpha_i$ block of 
$Y= (\cL^*)^{-1}(S)$ is given by
\begin{eqnarray*}
\lefteqn{
\left[\begin{array}{cc}
 Y_{ii} & Y_{i\alpha_i} \\ Y_{\alpha_ii} & Y_{\alpha_i\alpha_i} 
 \end{array}\right]} \\
 & =  &
 \left[\begin{array}{cc} 1/L_{ii} & -L_{\alpha_ii}^{\top} 
 L_{\alpha_i\alpha_i}^{-\top}/L_{ii}
 \\ 0 & L_{\alpha_i\alpha_i}^{-\top} \end{array}\right]
\left[\begin{array}{cc}
  S_{ii} & S_{\alpha_ii}^\top \\ S_{\alpha_ii} & S_{\alpha_i\alpha_i} 
 \end{array}\right]
 \left[\begin{array}{cc} 1/L_{ii} & 0 \\
 -L_{\alpha_i\alpha_i}^{-1}L_{\alpha_ii}/L_{ii} & 
 L_{\alpha_i\alpha_i}^{-1} 
 \end{array}\right].
\end{eqnarray*}
It can be computed as follows. Here the ``update matrices'' $V_j$ are 
defined as 
\[
V_i = Y_{\alpha_i\alpha_i} = 
L_{\alpha_i\alpha_i}^{-\top} S_{\alpha_i\alpha_i} 
L_{\alpha_i\alpha_i}^{-1}.
\]

\begin{algdesc}{Inverse adjoint mapping $({\cL}^*)^{-1}$} 
\label{a-Fadj-inv}
\begin{list}{}{}
\item[\textbf{Input.}]
A matrix $S \in \SS^N_E$ with homogeneous chordal sparsity pattern
and trivially perfect elimination ordering 
$\sigma = (1,2,\ldots,N)$, a nonsingular lower-triangular matrix 
$L\in\TT^N_E$, 
and the elimination tree for $\sigma$.
\item[\textbf{Output.}]  The matrix $Y = \Pi_E(L^{-\top}SL^{-1})$.
\item[\textbf{Algorithm.}]
\begin{enumerate}
\item For $i=1,\ldots,N$, set 
\[
 W_{ii} = S_{ii}, \qquad
 W_{\alpha_ii} = L_{\alpha_i\alpha_i}^{-\top} S_{\alpha_ii}.
\]
\item Enumerate the vertices $i=1,\ldots,N$ in reverse topological order.
 For each $i$, compute $Y_{ii}$ and $Y_{i\alpha_i}$ from
\begin{eqnarray*}
\lefteqn{
\left[\begin{array}{cc}
 Y_{ii} & Y_{i\alpha_i}^\top \\ Y_{\alpha_ii} & \times \end{array}\right]}
\\
 & = & \left[\begin{array}{cc} 
 1/L_{ii} & -L_{\alpha_ii}^\top/L_{ii}   \\ 0 & I \end{array}\right]
 \left[\begin{array}{cc} W_{ii} & W_{\alpha_ii}^\top \\
   W_{\alpha_ii} & V_i \end{array}\right]
  \left[\begin{array}{cc} 1/L_{ii} & 0 \\   
  -L_{\alpha_ii}/L_{ii} & I \end{array}\right], 
\end{eqnarray*}
and define
\[
 V_j = \left[\begin{array}{cc} Y_{ii} & Y_{\alpha_ii}^\top \\
   Y_{\alpha_ii} & V_i \end{array}\right], \quad
 j\in\ch(i).
\]
\end{enumerate}
\end{list}
\end{algdesc}

\subsection{Cholesky factorization} \label{s-chol}
Assume $X$ is positive definite with sparsity pattern $E$.
We define the Cholesky factorization as  a factorization 
$X = LL^\top$ with $L$ lower-triangular with positive diagonal elements.
If $\sigma = (1,2,\ldots,N)$ is a perfect elimination order, then
$L$ has the same sparsity pattern as $X$, i.e., $L\in\TT^N_E$.
In this section we specialize the multifrontal Cholesky factorization 
algorithm \cite{DuR:83,Liu:90,Liu:92}
to homogeneous chordal sparsity patterns.

Consider the $\bar\alpha_i$-by-$ \bar\alpha_i$ block of the 
factorization: 
\begin{eqnarray*}
\lefteqn{
\left[\begin{array}{cc}
 X_{ii} & X_{\alpha_ii}^\top \\ X_{\alpha_ii} & X_{\alpha_i\alpha_i}
 \end{array}\right]} \\
& =  &
 \left[\begin{array}{c} L_{ii} \\ L_{\alpha_ii} \end{array}\right]
 \left[\begin{array}{c} L_{ii} \\ L_{\alpha_ii} 
  \end{array}\right]^\top
+ \sum_{k<i} 
 \left[\begin{array}{c} L_{ik} \\ L_{\alpha_ik} \end{array}\right]
 \left[\begin{array}{c} L_{ik} \\ L_{\alpha_ik} 
  \end{array}\right]^\top
 +
 \sum_{k>i} 
 \left[\begin{array}{c} 0 \\ L_{\alpha_ik} \end{array}\right]
 \left[\begin{array}{c} 0  \\ L_{\alpha_ik} 
  \end{array}\right]^\top.
\end{eqnarray*}
If we consider only the first row and column in this equation, 
we can drop the last term on the right-hand side.  
In the second term, we can limit the sum to the vertices $k$ that
are proper descendants of $i$ in the elimination tree: 
\begin{equation}
\left[\begin{array}{cc}
 X_{ii} & X_{\alpha_ii}^\top \\ X_{\alpha_ii} & \times
 \end{array}\right] =
 \left[\begin{array}{c} L_{ii} \\ L_{\alpha_ii} \end{array}\right]
 \left[\begin{array}{c} L_{ii} \\ L_{\alpha_ii} 
  \end{array}\right]^\top
+ \sum_{j\in \ch(i)} \sum_{k\in T_j} 
 \left[\begin{array}{c} L_{ik} \\ L_{\alpha_ik} \end{array}\right]
 \left[\begin{array}{c} L_{ik} \\ L_{\alpha_ik} 
  \end{array}\right]^\top. \label{e-chol-update-eq}
\end{equation}
Here $T_j$ denotes the subtree of the elimination tree rooted at $j$. 
In the multifrontal method, one  defines for each node $j$ in the 
elimination tree an \emph{update matrix} 
\[
 U_j = -\sum_{k\in T_j} L_{\alpha_jk}L_{\alpha_jk}^\top.
\]
For a trivially perfect elimination ordering, 
$\alpha_j = \bar \alpha_i$  if $j\in \ch(i)$. 
The last term in~(\ref{e-chol-update-eq}) is therefore equal
to $-\sum_{j\in\ch(i)} U_j$, and the $22$ block of the entire 
right-hand side is $-U_i$.
Therefore
\[
\left[\begin{array}{cc}
 X_{ii} & X_{\alpha_ii}^\top \\ X_{\alpha_ii} & -U_i \end{array} \right]
 = 
 \left[\begin{array}{c} L_{ii} \\ L_{\alpha_ii} \end{array}\right]
 \left[\begin{array}{c} L_{ii} \\ L_{\alpha_ii} 
  \end{array}\right]^\top
- \sum_{j\in \ch(i)} U_j. \label{e-chol-update-eq2}
\]
Re-arranging this as
\begin{equation}
 \left[\begin{array}{cc}
   X_{ii} & X_{\alpha_ii}^\top \\ X_{\alpha_ii} & 0 \end{array}\right]
  + \sum_{j\in \ch(i)} U_j
 = \left[\begin{array}{c} L_{ii} \\ L_{\alpha_ii} \end{array}\right]
 \left[\begin{array}{c} L_{ii} \\ L_{\alpha_ii} \end{array}\right]^\top
 + \left[\begin{array}{cc} 0 & 0 \\ 0 & U_i \end{array}\right] 
\label{e-chol-recursion} 
\end{equation}
suggests a recursive algorithm for computing the factorization.

\begin{algdesc}{Cholesky factorization} \label{a-chol}
\begin{list}{}{}
\item[\textbf{Input.}]
A matrix $X \in \SS^N_E \cap \SS^N_{++}$,
with homogeneous chordal sparsity pattern and 
trivially perfect elimination ordering $\sigma = (1,2,\ldots,N)$, 
and the elimination tree for $\sigma$.
\item[\textbf{Output.}]  The Cholesky factorization $X = LL^\top$.
\item[\textbf{Algorithm.}]
Enumerate the vertices $i=1, 2, \ldots,N$ of the
elimination tree in topological order.  
For each $i$, form the \emph{frontal matrix}
\[
 \left[\begin{array}{cc}
  F_{11} & F_{21}^\top \\ F_{21} & F_{22} \end{array}\right] =
  \left[\begin{array}{cc}
     X_{ii}    & X_{\alpha_ii}^\top \\
     X_{\alpha_ii}    & 0 \end{array}\right] + \sum_{j\in\ch(i)} U_j.
\]
and calculate $L_{ii}$, $L_{\alpha_ii}$, and the update matrix $U_i$ from
\[
 L_{ii} = \sqrt{F_{11}}, \qquad 
 L_{\alpha_ii} = \frac{1}{L_{ii}}F_{21}, \qquad
 U_i = F_{22} -  L_{\alpha_ii} L_{\alpha_ii}^\top.
\]
\end{list}
\end{algdesc} 

\subsection{Maximum-determinant positive definite completion.} 
A matrix with a chordal sparsity pattern has a positive definite
completion if and only all completely specified principal submatrices
are positive definite \cite{GroneJohnsonSa1984}.
In our notation, $S\in\Pi_E(\SS^N_{++})$ if and only if
$S_{\bar\alpha_i\bar\alpha_i} \succ 0$ for all~$i$.
The positive definite completion with maximum determinant is the
inverse of a matrix $X \in \SS^N_E \cap \SS^N_{++}$.
If we parameterize $X=LL^\top$ by its Cholesky factor $L$, 
then $L\in\TT^N_E$ is the solution of the nonlinear equation
\[
\Pi_E(L^{-\top}L^{-1}) = S.
\]
The solution can be computed as follows 
\cite{AndersonDahlVandenberghe2013}.
Consider the $\bar\alpha_i $-by-$ \bar\alpha_i$ block of the 
equation $X^{-1}L = L^{-\top}$,
\begin{equation} \label{e-compl-eq}
 S_{\bar\alpha_i\bar\alpha_i} L_{\bar\alpha_i\bar\alpha_i}
   = L_{\bar\alpha_i\bar\alpha_i}^{-\top}.
\end{equation}
On the right-hand side we use $(L^{-1})_{\bar\alpha_i\bar\alpha_i}
 = L_{\bar\alpha_i\bar\alpha_i}^{-1}$, which holds for any nonsingular
lower-triangular matrix and any index set $\bar\alpha_i$.
On the left-hand side we use the fact that
the block column $\bar\alpha_i$ of $L$ has no zeros outside the rows 
indexed by $\bar\alpha_i$, since $1,\ldots,N$ is
a trivially perfect elimination ordering,
An algorithm for computing the Cholesky factor $L$ follows from the
first column of the equation~(\ref{e-compl-eq}):
\[
\left[\begin{array}{cc}
 S_{ii} & S_{\alpha_ii}^\top \\ S_{\alpha_ii} & S_{\alpha_i\alpha_i}
 \end{array}\right]
\left[\begin{array}{c} L_{ii} \\ L_{\alpha_ii} \end{array}\right]
= \left[\begin{array}{c} 1/L_{ii} \\ 0 \end{array}\right].
\]
The subvector $L_{\alpha_ii}/L_{ii}$ satisfies
\[
 \frac{1}{L_{ii}} L_{\alpha_ii} 
 = -S_{\alpha_i\alpha_i}^{-1} S_{\alpha_ii}  \\
 = -L_{\alpha_i\alpha_i} L_{\alpha_i\alpha_i}^\top S_{\alpha_ii}. 
\]
Substituting this in the first equation gives an expression for $L_{ii}$:
\[
 L_{ii} 
 =  
  (S_{ii} + S_{\alpha_ii}^\top (L_{\alpha_ii}/L_{ii}))^{-1/2} 
 =   
  (S_{ii} - S_{\alpha_ii}^\top L_{\alpha_i\alpha_i}
   L_{\alpha_i\alpha_i}^\top S_{\alpha_ii})^{-1/2}.
\]
In other words, if we define $u = L_{\alpha_i\alpha_i}^\top S_{\alpha_ii}$,
then  
\[
L_{ii} = (S_{ii} - \|u\|^2)^{-1/2}, \qquad L_{\alpha_ii} = -L_{ii}
   L_{\alpha_i\alpha_i}u.
\]
In the following outline we define $V_i = L_{\alpha_i\alpha_i}$.

\newpage

\begin{algdesc}{Maximum-determinant positive definite completion} 
\label{a-compl}
\begin{list}{}{}
\item[\textbf{Input.}]
A matrix $S \in \Pi_E(\SS^N_{++})$
with homogeneous chordal sparsity pattern and 
trivially perfect elimination ordering $\sigma = (1,2,\ldots,N)$, 
and the elimination tree for $\sigma$.
\item[\textbf{Output.}]  The nonsingular matrix $L\in \TT^N_E$ 
that satisfies $\Pi(L^{-\top}L^{-1}) = S$.
\item[\textbf{Algorithm.}]
Enumerate the vertices $i=1, 2, \ldots,N$ of the
elimination tree in inverse topological order.  
For each $i$, compute 
\[
 u = V_i^\top S_{\alpha_ii}, \qquad
 L_{ii} = \frac{1}{(S_{ii} - u^\top u)^{1/2}}, \qquad
 L_{\alpha_ii} = - L_{ii} V_i u.
\]
Then set
\[
 V_j = \left[\begin{array}{cc}  L_{ii} & 0 \\ L_{\alpha_ii} & V_i
 \end{array}\right], \quad j \in \ch(i).
\]
\end{list}
\end{algdesc}

\subsection{Gradient and Hessian of primal barrier}
In Section~\ref{s-barriers} we introduced the 
function $F(X) = -\ln\det X$ as logarithmic barrier function
for the cone $K = \SS^N_+ \cap \SS^N_E$.
Define $\cL(Y) = LYL^\top$, where $L$ is the Cholesky factor 
of $X$.
Then the gradient of $F$ at $X$, which is given by 
$F'(X) = -\Pi_E(X^{-1})$,
can be computed as
\[
F'(X) = -\Pi_E(L^{-\top}L^{-1}) = -(\cL^*)^{-1}(I).
\]
The following algorithm is Algorithm~\ref{a-Fadj-inv} with $X=I$.
It is also easily derived directly by considering the 
$\bar\alpha_i $-by-$ i$ block of the equation $X^{-1}L = L^{-\top}$, i.e.,
\begin{equation} \label{e-gradient-equation}
  \left[\begin{array}{cc}
    Y_{ii} & Y_{\alpha_ii}^\top \\ Y_{\alpha_ii} & Y_{\alpha_i\alpha_i}
 \end{array}\right] 
 \left[\begin{array}{c} L_{ii} \\ L_{\alpha_ii} \end{array}\right] 
 = \left[\begin{array}{c} 1/L_{ii} \\ 0 \end{array} \right].
\end{equation}
We define $V_i = (X^{-1})_{\alpha_i\alpha_i}$.

\begin{algdesc}{Projected inverse}   \label{a-projinv}
\begin{list}{}{}
\item[\textbf{Input}.]
The Cholesky factor $L$ of a positive definite matrix 
$X\in \SS^N_{++} \cap \SS^N_E$, with a homogeneous chordal sparsity
pattern and trivially perfect elimination ordering $1,\ldots, N$,
and an elimination tree for $\sigma$.
\item[\textbf{Output}.]
The projected inverse $Y = \Pi_E(L^{-\top}L^{-1})$.
\item[\textbf{Algorithm.}]
Enumerate the vertices $i= 1,2,\ldots,N$ in inverse topological 
order.  For each $i$, calculate 
\[
Y_{\alpha_ii} = -\frac{1}{L_{ii}} V_iL_{\alpha_ii}, \qquad
Y_{ii} = \frac{1}{L_{ii}} (\frac{1}{L_{ii}} - 
  L_{\alpha_ii}^\top Y_{\alpha_ii})
\]
and define the update matrices  
\[
 V_j =  \left[\begin{array}{cc} Y_{ii} & Y_{\alpha_ii}^\top \\ 
 Y_{\alpha_ii} & V_i \end{array}\right], \qquad j\in\ch(i). 
\]
\end{list}
\end{algdesc}
The Hessian of $F$ at $X$ is the linear mapping
\[
F''(X; Y) 
 = (\cL \circ \cL^*)^{-1}(Y),
\]
(see~(\ref{e-primal-Hessian-fact}))
and can be evaluated by calling algorithms~\ref{a-F-inv}
and~\ref{a-Fadj-inv}.

\subsection{Gradient and Hessian of dual barrier}
The barrier for the cone $\Pi_E(\SS^N_+)$ is 
\[
F_*(S) = \sup_X{(-\tr(SX) - F(X))} = N - F(\hat X)
\]
where $\hat X$ is the maximizer in the definition, i.e., the solution
of the equation $\Pi_E(X^{-1}) = S$.
Define $\cL(Y) = LYL^\top$, where $L$ is the Cholesky factor of 
$\hat X$, which can be computed by algorithm~\ref{a-compl}.
The gradient of $F_*$ at $S$ is 
\[
F'_*(S) = -\hat X = -\cL(I)
\]
and can be computed by applying algorithm~\ref{a-F} with $X=I$.

\begin{algdesc}{Dual gradient}
\begin{list}{}{}
\item[\textbf{Input.}]
The Cholesky factor $L$ of the inverse of the maximum-determinant
positive definite completion of a matrix $S \in \Pi_E(\SS^N_{++})$,
with a homogeneous chordal sparsity pattern and trivially perfect
elimination ordering $\sigma = (1,2,\ldots,N)$,
and the elimination tree for $\sigma$.
\item[\textbf{Output.}]  The matrix $Y = LL^\top$.
\item[\textbf{Algorithm.}]
Enumerate the vertices $i=1, 2, \ldots,N$ of the
elimination tree in topological order.  
For each $i$, compute $U_i$, $Y_{ii}$, $Y_{\alpha_ii}$ using
the formula
\[
\left[\begin{array}{cc}
 Y_{ii} & Y_{i\alpha_i} \\ Y_{\alpha_ii} & -U_i\end{array}\right]
 = \left[\begin{array}{c} L_{ii} \\ L_{\alpha_ii} \end{array}\right]
   \left[\begin{array}{c} L_{ii} \\ L_{\alpha_ii} \end{array}\right]^\top
 - \sum_{j\in \ch(i)} U_j.
\]
\end{list}
\end{algdesc} 
The Hessian of $F_*$ is given by 
$F''_*(S) = F''(\hat X)^{-1} = \cL\circ \cL^*$
and can be evaluated 
via Algorithms~\ref{a-F} and~\ref{a-Fadj}.

\bibliographystyle{plain}
\bibliography{TuncelVandenberghe}
\end{document}